\documentclass[12pt]{article}

\usepackage{epsfig}

\usepackage{amssymb} 
\usepackage{amsmath} 

\def\ZZ{\mbox{{\sl Z\hspace{ -.3em}Z}}}

\def\RR{\mbox{{\sl l\hspace{ -.15em}R}}}
\def\CC{\mbox{{\sl l\hspace{ -.47em}C}}}
\def\HH{\mbox{{\sl l\hspace{ -.15em}H}}}
\def\PP{\mbox{{\sl l\hspace{ -.15em}P}}}

\def\iint{\int\!\!\int}

\newtheorem{theo}{Theorem}[subsection]
\newtheorem{lemm}[theo]{Lemma}
\newtheorem{coro}[theo]{Corollary}
\newtheorem{defi}[theo]{Definition}
\newtheorem{prop}[theo]{Proposition}
\newtheorem{clai}[theo]{Claim}
\newtheorem{subl}[theo]{Sublemma}
\newtheorem{conj}[theo]{Conjecture}
\newtheorem{ques}[theo]{Question}
\newtheorem{exam}[theo]{Example}
\newenvironment{rema}{{\noindent \bf Remark:}}{}
\newenvironment{demo}{{\noindent \bf Proof: }}{\hfill$\square$}

\font\b=cmr10 scaled \magstep 3
\def\bigzerou{\smash{\lower1.7ex\hbox{\b 0}}} 
\def\R{\Re e}
\def\I{\Im m}

\title{\bf Conformally invariant energies of knots}
\author{R. ~Langevin and J. ~O'Hara\\
{\small Laboratoire de Topologie, CNRS-UMR 5584, Universit\'e de Bourgogne,} \\
{\small Department of Mathematics, Tokyo Metropolitan University}}

\begin{document}

\maketitle

\vskip 5mm
\noindent {\bf Abstract}

Conformally invariant functionals on the space of knots 
are introduced via extrinsic conformal geometry of the knot 
and integral geometry on the space of spheres. 
Our functionals are expressed in terms of a complex-valued 2-form 
which can be considered as the cross-ratio 
of a pair of infinitesimal segments of the knot. 
We show that our functionals detect the unknot 
as the total curvature does, and that their values 
explode if a knot degenerates to a singular knot 
with double points.

\tableofcontents

\vskip 5mm

\noindent {\bf keywords:} knot, energy, conformal, cross-ratio, M\"obius transformations.


\vskip 5mm
\section{Introduction. }
In 1929 Fenchel 
\cite{Fe1} 
proved that the total curvature of a closed curve is greater than 
or equal to $2\pi$. 
Between 1949 and 1951, Fary 
\cite{Far}
Fenchel 
\cite{Fe2}
and Milnor 
\cite{Mi}
proved independently 
that the total curvature of a non-trivial knot is 
greater than $4\pi$. 
This means that non-triviality of knots imposes 
a jump on the total curvature functional. 
For knots in $S^3$ one has to consider the sum of the length and the 
total curvature to get such a jump for non-trivial knots 
(\cite{Ban}, \cite{La-Ro}). 

In 1989 the second author introduced 
the {\it energy} $E^{(2)}$ of a knot 
\cite{Oh1} 
as the regularization of modified electrostatic energy of charged knots. 
This energy $E^{(2)}$ is a functional on the space of knots 
which explodes if a knot degenerates to a singular knot 
with double points. 
In 1994 Freedman, He, and Wang 
\cite{F-H-W} 
showed that the energy is conformally invariant 
and that 
there is a lower bound for the energy of non-trivial knots, 
which can detect the unknot. 

In this paper we introduce several conformally invariant 
functionals on the space of knots 
that satisfy the two properties above, namely, 
the jump imposed by the non-triviality and the explosion 
for the singular knots. 
These functionals can be expressed in terms of a complex-valued 
2-form on the configuration space $K\times K\setminus\triangle$ 
of a knot $K$, 
which is obtained as the cross-ratio of 
$dx$ and $dy$ 
on the 2-sphere $\Sigma\cong\CC\cup\{\infty\}$ 
that is twice tangent to the knot at $x$ and $y$.   

\medskip
This paper is arranged as follows:

In section 2 we introduce the geometry of the space of $n$-spheres. 
We realize it as a hypersurface in 
Minkowski space 
and give an invariant measure on it. 
A 4-tuple of points of $\RR^3$ or $S^3$ is contained in a sphere $\Sigma$ (unique if the four points do not belong to a circle) which allows us to 
define the cross-ratio of four points. 

In section 3 we make a brief review of the energy $E^{(2)}$, and introduce 
the cosine formula given by Doyle and Schramm. 

In section 4 we show that the integrand of $E^{(2)}$ 
can be interpreted in terms of the {\it infinitesimal cross-ratio}, 
which is the cross-ratio of four points $(x, x+dx, y, y+dy)$ 
on the knot, 
by using the 2-sphere which is twice tangent to the knot at $x$ and $y$. 
We show that the real part of this infinitesimal cross-ratio 
coincides with the pull-back of the canonical symplectic form 
of the cotangent bundle of $S^3$, which allows us to deduce 
the original definition of $E^{(2)}$ from its cosine formula. 

In section 5 we define another functional in terms of this 
infinitesimal cross-ratio 
and show that it has the jump and the explosion properties. 

In section 6 we introduce a {\it non-trivial sphere} 
for a knot, which is a sphere that intersects the knot 
in at least 4 points. 
We give a formula of the measure of the set of the non-trivial spheres 
for a knot (counted with multiplicity), 
which will be called the {\it measure of acyclicity}, 
in terms of the infinitesimal cross-ratio. 

In section 7 
we then study the position of the knot with respect to a {\it zone} which is a region bounded by 
two disjoint 2-spheres. 
Recall that a zone is characterised up to conformal transformation 
by one number which we call the {\it modulus}. 
One way to compute it is to send the two boundary spheres by 
a conformal transformation into concentric position. 
The modulus is then the ratio of the radii. 
By showing that too ``thin" zones can not catch the topology 
of a given non-trivial knot, 
we prove that the measure of acyclicity has 
the jump and the explosion properties. 

We put an appendix as section 8. 
In subsection 8.1 we show that the integrand of $E^{(2)}$ can be 
interpreted as the maximal modulus of an infinitesimal zone 
whose boundary spheres pass through pairs of points 
$\{x,y\}$ and $\{x+dx, y+dy\}$ respectively. 
In subsection 8.2 we introduce a conformal geometric generalization of 
the Gauss' integral formula for the linking number. 
In subsection 8.3 we express the radius of global curvature, 
which was introduced by Menger
and Gonzalez and Maddocks \cite{Go-Ma} 
in terms of the {\it 4-tuple map} introduced in section 2. 

\medskip
Throughout the paper, for simplicity's sake, 
when we consider knots in $\RR^3$ we assume that 
they are parametrized by arc-length unless otherwise mentioned. 
In what follows a knot $K$ always means the image of an embedding $f$ 
from $S^1$ into $\RR^3$ or $S^3$. 

{\it 
 The authors thank J.P. Troalen for making the pictures,
and J. Cantarella, Y. Kanda, R. Kusner, J. Maddocks, K. Ono, and D.
Rolfsen for
useful conversations.}


\section{The space of spheres and the cross-ratio of 4 points. }

\subsection{The space of $(n-1)$-spheres in $S^{n}$. }

\begin{defi}\rm 
(1) Define an indefinite quadratic form $L$ on $\RR^{n+2}$ and 
the associated {\it pseudo inner product} $L(\cdot ,\cdot )$ 
$$
\begin{array}{l}
L(x_1,\cdots ,x_{n+2})=(x_1)^2+\cdots +(x_{n+1})^2-(x_{n+2})^2,\nonumber\\
L(u,v)=u_1v_2+\cdots u_{n+1}v_{n+1}-u_{n+2}v_{n+2}.\nonumber
\end{array}
$$
This quadratic form with signature $(n+1,1)$ is called the 
{\it Lorentz quadratic form}. 
The Euclidean space $\RR^{n+2}$ equipped with this pseudo inner product $L$ 
is called the {\it Minkowski space} and denoted by $\RR^{n+1,1}$. 

\medskip
(2) The set ${\mathcal G}=O(n+1, 1)$ of linear isomorphism of $\RR^{n+2}$ 
which preserves this pseudo inner product $L$ is called the {\it 
(full homogeneous) Lorentz group}. We also call it the {\it conformal group} 
for short. 

\medskip
(3) A vector $v$ in $\RR^{n+2}$ is called {\it space-like} if $L(v)>0$, 
{\it light-like} if $L(v)=0$, and 
{\it time-like} if $L(v)<0$. 
A line is called space-like (or time-like) if it contains a space-like
 (or respectively, time-like) vector.
The isotropy cone $\{v\in\RR^{n+2} \, |\, L(v)=0\}$ of $L$ 
is called the {\it light cone}. 

\medskip 
(4) The points at infinity of the light cone in the upper half space 
$\{x_{n+2} >0\}$ form an $n$-dimensional sphere. Let it be denoted by 
$S_{\infty} ^{n}$. 
Since it can be considered as the set of lines through the origin 
in the light cone, 
it is identified with the intersection $S_1 ^{n}$ of 
the upper half light cone and the hyperplane $\{x_{n+2}=1\}$, 
which is given by 
$S_1 ^{n}=\{(x_1, \cdots , x_{n+1}\, |\, (x_1)^2+\cdots +(x_{n+1})^2-1=0\}$. 

\medskip
(5) An element of the Lorentz group ${\mathcal G}$ acts on $S_{\infty} ^{n}$ 
as a {\it conformal transformation} which is a composition of 
reflections with respect to $(n-1)$-spheres in $S_{\infty} ^{n}$. 
The set of conformal transformations of 
$S_{\infty} ^{n}\cong\RR^n\cup\{\infty\}$ 
is called the {\it conformal group} 
and is also denoted by ${\mathcal G}$. 

\end{defi}

\begin{clai} 
Let ${\mathcal S}$ be the set of oriented $(n-1)$-spheres in $S^{n}$, and 
let $\Lambda=\{x\in\RR^{n+2}\, |\, L(x)=1\}$ be the hyperbolic quadric 
hypersurface of one sheet. 
(This $\Lambda$ is called de Sitter space.) 
Then there is a canonical bijection between ${\mathcal S}$ and $\Lambda$. 
\end{clai} 

\begin{figure}[ht]
\centerline{\input{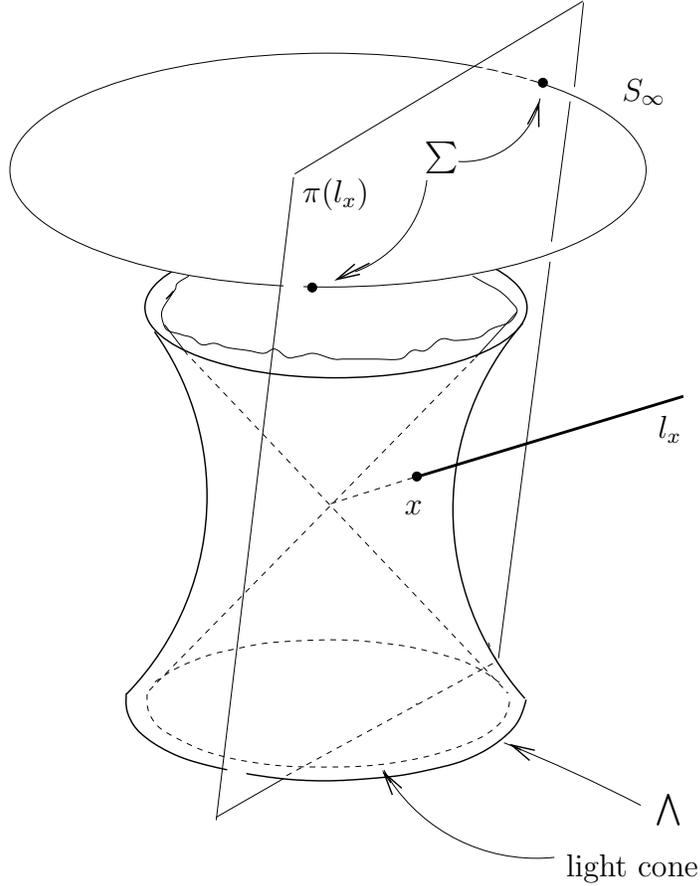}}
\caption{$S^n_{\infty}$ and the correspondence between $\Lambda$ and ${\mathcal S}$.}    
\end{figure}

\begin{demo} 
A point $x$ in $\Lambda$ determines an oriented half line 
$l_x=\stackrel{\longrightarrow}{Ox}$ from the origin. 
Let $\Pi (l_x)$ be the oriented hyperplane passing through 
the origin 
that is orthogonal to $l_x$ with respect to Lorentz quadratic form $L$. 
Since $l_x$ is space-like $\Pi (l_x)$ intersects the light cone transversely 
and therefore $\Pi (l_x)$ intersects $S_{\infty}^3$ in an oriented sphere 
$\Sigma_{\Pi (l_x)}$. The map
$
\Lambda\ni x\mapsto\Sigma_{\Pi(l_x)}\in{\mathcal S}
$
defines the bijection from $\Lambda$ to ${\mathcal S}$. 
\end{demo}

Let us identify $\Lambda$ with ${\mathcal S}$ through this bijection.

Let $S_1$ and $S_2$ be oriented 2-spheres. 
Considered as points in $\Lambda$, $S_1$ and $S_2$ satisfy:

\begin{equation}
\begin{array}{l}\label{intersect_or_not}
 {\rm (1)} S_1\cap S_2\ne\phi \ \mbox{if and only if}\ |L(S_1, S_2)|\le 1\\ 
 {\rm (2)} S_1\cap S_2=\phi \ \mbox{ if and only if}\ |L(S_1, S_2)|>1. 
\end{array}
 \end{equation}
Two spheres $S_1$ and $S_2$ are said to be {\it nested} if $S_1\cap S_2=\emptyset$ 
and {\it intersecting} if $S_1\cap S_2\ne\emptyset$ and $S_1\neq S_2$.

Dimension 2 planes through the origin intersect the quadric 
$\Lambda \subset \RR ^{n+2}$ in curves. Those curves are geodesics 
in the sense that they are critical (but not minimal) points 
of the length function which assigns the length $\int_a ^b \sqrt{|L (c'(t))|} dt$ 
to any given arc $c:[a,b]\rightarrow \Lambda$. 
There are three cases. 

\smallskip 
(1a) When the plane contains a time-like vector, the intersection consists of 
two non-compact curves. The pencil of spheres has limit (or Poncelet) points 
which are $P\cap S_{\infty}^n$; it is called a {\it Poncelet pencil}. 

In case (1) of the formula \ref{intersect_or_not}, the length of one of the arcs joining $S_1$ to $S_2$ in the pencil is the angle $\theta_0$ between the two spheres.
This angle satisfies: $L(S_1 ,S_2 ) = \cos(\theta_0 )$.

\smallskip 
(1b) When the plane is tangent to the light cone, the intersection is again two 
non-compact curves, but the pencil is made of spheres all tangent at the 
point $P \cap S_{\infty}^n$. 

\smallskip
(2) When all the vector of the plane $P$ are space-like, its intersection with 
$\Lambda$ is connected and closed. Any tangent vector to this intersection 
is space-like. The corresponding set of spheres is a {\it pencil} of 
spheres with a common $(n-2)$-sphere 
\newline $\Gamma =P^{\bot} \cap S_{\infty}^n$  
($P^{\bot}$ means the subspace $L-orthogonal$ to $P$, that is 
\newline $P^{\bot}= \{ w |  L(v,w) = 0 \ \forall v\in P \}$) .

In case (2) of the formula \ref{intersect_or_not}, we can get an interesting invariant of the annulus bounded by $S_1$ and $S_2$ from the number $|L(S_1 ,S_2 )|$ (see \ref{rem_module}): the length $t_0$ of the arc of the pencil joining $S_1$ to $S_2$ satisfies $|L(S_1 ,S_2 )|= Ch(t_0)$.

\begin{figure}[ht]
\centerline{\input{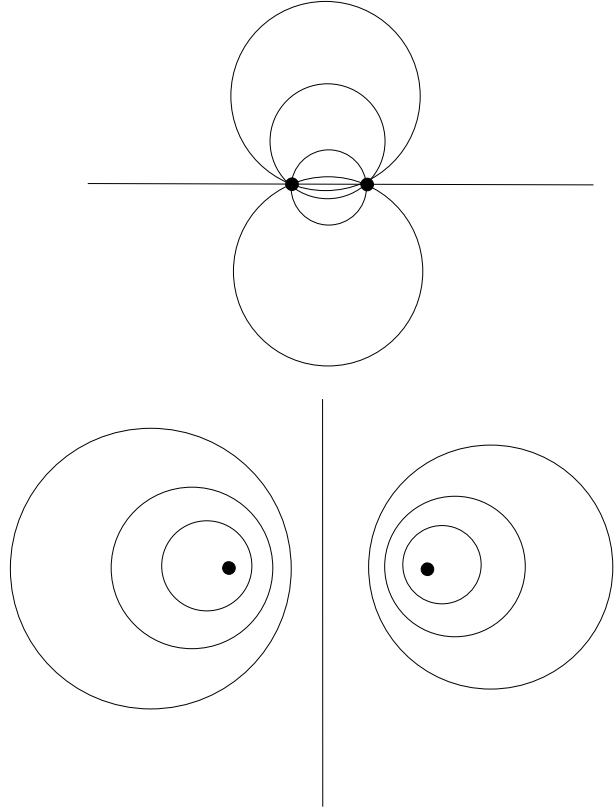}}
\caption{A pencil with a base circle and a Poncelet pencil.}    
\end{figure}

Working in $\Lambda \subset \RR^4$ we get the ``usual" theory of pencils of 
circles. We will be mainly interested in $\Lambda \subset \RR^5$ 
which is the set of oriented 2-spheres in $S^3$. 
We do not study the case when $n>3$. 

\begin{rema}
Let $(v_1 ,v_2 ,... ,v_{n+1} )$ be $n+1$ vectors in $T_x \Lambda $. 
The volume of the parallelepiped  
constructed on these vectors is
$$|\det (x ,v_1 ,v_2 ,... ,v_{n+1} )|= \sqrt{-\det({L}(v_i ,v_j))}$$
\end{rema}
This result is a corollary of the lemma \ref{parallel}.

\medskip
Let $\Delta=\{x\in\RR^{n+2}\, |\, L(x)=-1, \, x_{n+2}>0\}$ be 
the upper half of the hyperbolic quadric hypersurface of two sheets. 
The restriction of the quadratic form $L$ to each tangent space of 
$\Delta$ is positive definite. 
Then $\Delta$ is a model for the $n+1$-dimensional hyperbolic space 
$\HH=\HH^{n+1}$. 
Each sphere $\Sigma$ in $S^{n}$ is the  ``boundary at infinity" of a 
totally geodesic subspace  $h$ of $\HH$.
The restriction of the full homogeneous Lorentz group 
${\mathcal G}=O(n+1, 1)$ to $\HH$ is the group of isometries 
of the hyperbolic space $\HH$.

A choice of a point $z$ in $\HH$ determines a metric on the sphere $S^{n}$ 
at infinity of $\HH$ by the projection on $S^{n}$ of the metric 
of $T_z(\HH)$ using the geodesic rays which starts from $z$. 
Another choice of the point $z$ determines 
another conformally equivalent metric on the sphere $S^{n}$. 
Although the sphere does not admit a measure which is invariant 
under the action of the conformal group ${\mathcal G}$, 
the set ${\mathcal S}\cong \Lambda$ of spheres in $S^{n}$ do. 
Namely, $\Lambda$ is endowed with a ${\mathcal G}$-invariant measure 
$d\Sigma$. We define $d\Sigma$ in several ways in what follows. 

\medskip 
(1) The Lebesgue measure $|dx_1 \wedge \cdots \wedge dx_{n+2}|$ 
is invariant under the action of $\mathcal G$. 
It inherits from $\RR^{n+2}$ a ${\mathcal G}$-invariant measure 
defined on each tangent space at a point $x\in \Lambda$ by 
$d\Sigma = |\iota _x (dx_1  \wedge\cdots \wedge  dx_{n+2})|$. 

\medskip 
(2) This measure is given by the volume $(n+1)$-form 
$$
\omega _{\Lambda} 
=\frac{\sum_{i=1}^{n+2}(-1)^ix_idx_1\wedge \cdots 
\widehat{\wedge dx_i}\cdots \wedge dx_{n+2}}
{({x_1}^2+\cdots+{x_{n+1}}^2-{x_{n+2}}^2)^{\frac{n+2}2}}
$$ 
of $\Lambda$ associated with the Lorentz form $L$, 
where $\widehat{\wedge dx_i}$ means that $\wedge dx_i$ is removed. 
Then $\omega _{\Lambda}$ is invariant under the action of 
${\mathcal G}_+=SO(n+1,1)$. 
We remark that $\omega _{\Lambda}$ is expressed as 
$$
\omega _{\Lambda _+}=\frac1{x_{n+2}}dx_1\wedge 
\cdots  \wedge dx_{n+1}
$$
on $\Lambda _+=\{(x_1,\cdots ,x_{n+2})\in\Lambda |x_{n+2}>0\}$ 
if we take $(x_1,\cdots, x_{n+1})$ as its local coordinate. 

\medskip 
(3) This measure can also be regarded as a measure on the set of 
totally geodesic 
hyperplanes of the hyperbolic space $\HH$ which is invariant by the action 
of the hyperbolic isometries \cite{San}. 

\medskip 
(4) Let us now project the sphere $S^{n}$ stereographically on 
an affine space 
$\RR^{n}$ with Euclidean coordinates $(x_1 ,x_2 ,...,x_{n})$. 
There, a sphere $\Sigma$ is given by its center $(x_1,x_2,...,x_{n})$ and 
its radius r.
The measure $d\Sigma$ is expressed by 
$$
d\Sigma\, =\, \frac1{r^{n+1}}
|dx_1\wedge dx_2\wedge ...\wedge dx_{n}\wedge dr|
$$

The quadric $\Lambda$ is unbounded and the total measure of $\Lambda$ 
is infinite. 
The measure of the set of spheres intersecting a given arc of any size 
is infinite.
Since the set of spheres with geodesical radius larger than $r$ 
is compact in $\Lambda$, its volume is finite.

\begin{rema}
In what follows we assume that circles and spheres include 
lines and planes respectively as we work in conformal geometry. 
\end{rema}

\subsection{The 4-tuple map and the cross-ratio of 4 points. \label{ss4tuple}}

In this subsection we define the cross-ratio of ordered four points 
in $S^3$ or $\RR^3$ 
by means of the oriented 2-sphere that passes through them, 
which is determined uniquely unless the four points are concircular. 

The {\it configuration space} ${\rm Conf}_n(X)$ of a space $X$ is given by 
$$
\begin{array}{rcl}  
{\rm Conf}_n(X)&=&\{(x_1, \cdots , x_n)\}\in X^n | 
x_{i}\ne x_{j} \hskip 0.4cm {\rm if} \hskip 0.4cm i\ne j\}\\[1mm]
 &=&X\times \cdots \times X\setminus \triangle , 
\end{array} 
$$
where $\triangle$ is called the {\it big diagonal set}. 
Let ${\mathcal C}c(S^3)$ be the subset of ${\rm Conf}_4(S^3)$ 
consisting of concircular points. 

In what follows we define a map 
$\Sigma$ from ${\rm Conf}_4(S^3)\setminus {\mathcal C}c(S^3)$ 
to the set of oriented 2-spheres, ${\mathcal S}\cong\Lambda$, 
which assigns to ordered 4 points $(x, y, z, w)$ the oriented 2-sphere 
$\Sigma (x, y, z, w)$ that passes through them. 

\begin{defi}\rm 
The {\it (4,1)-exterior product} $v^1\wedge \cdots \wedge v^4$ of 
ordered four vectors 
$v^i=(v^i_1, \cdots , v^i_5)$ in $\RR^5$ $(i=1, 2, 3, 4)$ 
is the uniquely determined vector $u$ that satisfies 
$$ 
L(u,w)= det(v^1 , v^2 , v^3 , v^4 , w) \in\RR
\hskip 0.7cm (\forall \>\> w \in\RR^5). 
$$
Then it is given 
by 
$$
v^1\wedge \cdots \wedge v^4
=(\tilde\nu_1, \tilde\nu_2, \tilde\nu_3, \tilde\nu_4, -\tilde\nu_5)
$$
with
$$
\tilde\nu_j =(-1)^{j+1}
\left|
\begin{array}{ccccc}
v^1_1 & \cdots & \hat v^1_j & \cdots & v^1_5 \\
\vdots &       &  \vdots  &          &  \vdots  \\
v^4_1 & \cdots & \hat v^4_j & \cdots & v^4_5 \\
\end{array}
\right|,
$$
where $\hat {}$ means that the $j$-th column is removed. 
\end{defi}
We remark that this is a generalization of the vector product 
of two vectors in $\RR^3$. 

\begin{lemm}\label{parallel}
\noindent {\rm (1)} $v^1\wedge v^2\wedge v^3\wedge v^4\ne 0$ if and only if 
$v^1, v^2, v^3, v^4$ are linearly independent. 

\medskip
\noindent {\rm (2)} $L(v^1\wedge v^2\wedge v^3\wedge v^4 , v^j)=0$ 
for $j=1, 2, 3, 4$, namely, $v^1\wedge v^2\wedge v^3\wedge v^4$ is 
$L$-orthogonal to $Span\langle v^1, v^2, v^3, v^4\rangle$. 

\medskip
\noindent {\rm (3)} The norm of $v^1\wedge v^2\wedge v^3\wedge v^4$ is equal to 
the absolute value of the volume of the  parallelepiped  
spanned by $v^1, v^2, v^3$ and $v^4$ associated with 
the Lorentz quadratic form $L$: 
$$
\sqrt{|L(v^1\wedge v^2\wedge v^3\wedge v^4)|}
=\sqrt{|\det (L(v^i, v^j))|}. 
$$
\end{lemm}

\begin{demo}
(3) The last equality is a consequence of a formula in linear algebra. 
$$
\begin{array}{rcl}  
\det (L(v^i, v^j))
&=& \det \left\{\left(
\begin{array}{ccc}
v^1_1 & \cdots & v^1_5 \\
\vdots & & \vdots \\
v^4_1 & \cdots & v^4_5
\end{array}
\right)
\left(
\begin{array}{ccc}
v^1_1 & \cdots & v^4_1 \\
\vdots & & \vdots \\
v^1_4 & \cdots & v^4_4\\
-v^1_5 & \cdots & -v^4_5
\end{array}
\right)\right\}
\\[2mm]
&=&  -\tilde\nu_1^2-\tilde\nu_2^2-\tilde\nu_3^2-\tilde\nu_4^2+\tilde\nu_5^2
\end{array} 
$$
\end{demo}

\begin{clai}
 
\noindent {\rm (1)} Let $u=(u_1, \cdots , u_5)$ and $v=(v_1, \cdots , v_5)$ 
be linearly independent light-like vectors with $u_5, v_5>0$. 
Then $L(u,v)<0.$

Therefore, any pair of linearly independent light-like vectors can not be $L$-orthogonal. 

\medskip
\noindent {\rm (2)} Let 
$u=(u_1, \cdots , u_5)$ be a non-zero light-like vector. Then the hyperplane $u^{\bot}$ is tangent to the light-cone along the generatrix $\RR \cdot u$.

\medskip 
\noindent {\rm (3)}
If $u=(u_1, \cdots , u_5)$ is a non-zero time-like vector, then the hyperplane $u^{\bot}$ is space-like. 
\label{lorentzform}
\end{clai}

\begin{coro} 
Suppose $v^1, v^2, v^3, v^4$ are linearly independent. 

\noindent {\rm (1)} If $v^1\wedge v^2\wedge v^3\wedge v^4$ is time-like, then  all of $v^1, v^2, v^3, v^4$ are space-like. 

\medskip 
\noindent {\rm (2)} $v^1\wedge v^2\wedge v^3\wedge v^4$ is light-like if and only if 
$Span\langle v^1, v^2, v^3, v^4\rangle$ is the tangent space of 
the light cone at $v^1\wedge v^2\wedge v^3\wedge v^4$. 
\end{coro}

Let us identify the 3-sphere $S^3$ with 
the intersection $S^3_1$ of the light-cone $\{L=0\}$ 
and the hyperplane in $\RR^5$ defined by $\{(x_1, \cdots , x_5)|x_5=1\}$. 
Let 
$x^i=(x^i_1, x^i_2, x^i_3, x^i_4)\in S^3$ $(i=1,2,3,4)$ 
and 
$\tilde x^i=(x^i_1, x^i_2, x^i_3, x^i_4, 1)\in S^3_1\subset\RR^5$. 
Then $\tilde x^i$'s are linearly dependent if and only if 
$x^i$'s are concircular. 

Since an oriented 2-sphere $\Sigma (x^1,x^2,x^3,x^4)$ 
that passes through $\{x^1,x^2,x^3,x^4\}$ is obtained 
as the intersection of $S^3_1$ and an oriented hyperplane in $\RR^5$ 
that passes through $x^1,x^2,x^3,x^4$, and the origin, we obtain:  

\begin{prop} Define the  ``4-tuple map'' of the 3-sphere by 
$$
\Sigma :{\rm Conf}_4(S^3)\setminus {\mathcal C}c(S^3)\ni (x^1,x^2,x^3,x^4)
\mapsto 
\frac{\tilde x^1\wedge \tilde x^2\wedge \tilde x^3\wedge \tilde x^4}
{\sqrt{L(\tilde x^1\wedge \tilde x^2\wedge \tilde x^3\wedge \tilde x^4)}}
\in\Lambda\cong {\mathcal S}. 
$$
Then $\Sigma (x^1,x^2,x^3,x^4)$ is 
the oriented 2-sphere that passes through $\{x^1,x^2,x^3,x^4\}$. 
\end{prop}

Assume that the complex plane $\CC$ has the standard orientation. 

\begin{defi}\rm  
Let $(x^1,x^2,x^3,x^4)\in {\rm Conf}_4(S^3)$. 
Let $p:\Sigma (x^1,x^2,x^3,x^4)\to\CC$ 
be any stereographic projection that preserves the orientation, 
where, when $x^1,x^2,x^3$, and $x^4$ are concircular, 
we can choose any oriented 2-sphere 
that passes through them 
as 
$\Sigma (x^1,x^2,x^3,x^4)$. 

The {\it cross-ratio of ordered four points} 
$(x^1,x^2,x^3,x^4)$ is defined as 
the cross-ratio 
\begin{equation} \label{cross}
(p(\tilde x^2), p(\tilde x^3); p(\tilde x^1), p(\tilde x^4))
=\frac{p(\tilde x^2)-p(\tilde x^1)}{p(\tilde x^2)-p(\tilde x^4)}:
\frac{p(\tilde x^3)-p(\tilde x^1)}{p(\tilde x^3)-p(\tilde x^4)}. 
\end{equation}
We denote it by  $(x^2,x^3;x^1,x^4)$. 
When $x^1,x^2,x^3$ and $x^4$ are concircular, 
their cross-ratio is real and independent of the choice 
of $\Sigma (x^1,x^2,x^3,x^4)$. 
\end{defi}

\begin{lemm}\label{sign_cross}
The image of ${\rm Conf}_4(S^3)\setminus {\mathcal C}c(S^3)$ by the cross-ratio map of formula \ref{cross} is contained in a component of $\CC \setminus \RR$. 
\end{lemm}

\begin{demo} The set $\triangle \subset (S^3 )^4$ is of codimension $3$; the set ${\mathcal C}c(S^3)$ is of codimension $2$ in ${\rm Conf}_4(S^3)$ and the cross-ratio map is continuous on the arcwise connected set ${\rm Conf}_4(S^3)\setminus {\mathcal C}c(S^3)$. 
\end{demo}

The {\it 4-tuple map $\Sigma$ of} $\RR^3$ and the {\it cross-ratio of 
four points in} $\RR^3$ are defined through 
an orientation preserving 
stereographic projection : $p:S^3\to \RR^3$. 
When $x^1,x^2,x^3$, and $x^4$ are concircular in $\RR^3$, we define 
$\Sigma (x^1,x^2,x^3,x^4)$ to be any oriented sphere that passes through them.



\section{Review of $r^{-2}$-modified potential energy. }

\subsection{What is an energy functional for knots? }

In the following three sections we study knots in $\RR^3$. 
By an {\it open knot} 
we mean an embedding $\tilde f$ from a line into $\RR^3$ 
or its image $\tilde K=\tilde f(\RR)$ 
that approaches 
asymptotically 
a straight line at both ends.
We assume that knots are of class $C^2$. 

The {\it energy of knots} was proposed by 
Fukuhara \cite{Fuk} and Sakuma \cite{Sak} motivated by the following problem: 
Define a suitable functional on the space of knots, 
which we call ``{\it energy}'', and 
define a good-looking ``{\it canonical knot}'' for each isotopy class 
as one of the embeddings that attain the minimum value of this ``energy'' 
within its isotopy class. 

For this purpose we try to deform a given embedding 
along the negative gradient flow of the ``energy''  
until it comes to a critical point without changing its knot type. 
Hence the crossing changes should be avoided. 

Thus we are lead to the notion of the {\it energy functional for knots} 
which is a functional that explodes  
if a knot degenerates to an immersion with double points. More precisely; 

\begin{defi} \rm  
A functional $e$ on the space of knots is an {\it energy functional for knots} 
if for any real numbers $b$ and $\delta$ with $0<\delta \le\frac12$, 
there exists a positive constant $C=C(b,\delta)$ which satisfies 
the following condition: 

If a knot $f$ with length $l(f)$ contains a pair of points $x$ and $y$ 
which satisfy 
that the shorter arc-length between them is equal to $\delta l(f)$ 
and that $|x-y|\le Cl(f)$, 
then $e(f)\ge b$. 
\label{defofenergyfunctional}
\end{defi}
We call a functional an {\it energy functional in the weak sense} 
when it satisfies the same condition as above 
if the knot $f$ satisfies an additional geometric condition; 
for example, if the curvature of $f$ is bounded above 
by a constant. 

\medskip
\subsection{The regularization of $r^{-2}$-modified potential energy, $E^{(2)}$.}

One of the most natural and na\"{\i}ve candidates for an energy functional 
for knots would be the electrostatic energy of charged knots. 

The first attempt to consider the electrostatic energy of charged knots 
was carried out in the finite dimensional category by Fukuhara \cite{Fuk}. 
He considered the space of the polygonal knots whose vertices are charged, 
and studied their {\it modified} electrostatic energy 
under the assumption that Coulomb's repelling force between 
a pair of point charges of distance $r$ is proportional to 
$r^{-m}$ $(m=3,4,5, \cdots)$. 

The first example of an energy functional for smooth knots 
was defined by the second author in \cite{Oh1}. 
(The reader is referred to a survey article \cite{Oh2} for details.) 

Let $K=f(S^1)$ with $f:S^1=[0,1]/\sim \to \RR^3$ 
be a knot of class $C^2$ that is parametrized by the arc-length. 
Suppose that the knot is uniformly electrically charged. 
In order to obtain an energy functional for knots, we have to 
make the non-realistic assumption that Coulomb's repelling force 
between a pair of point charges of distance $r$ is proportional to 
$r^{-(\alpha +1)}$ 
and hence the potential is proportional to $r^{-\alpha }$ 
where $\alpha \ge 2$. 
Put $\alpha =2$ in what follows. 
Then its ``{\it $r^{-2}$-modified potential energy}'' is given by 
$$
\iint _{K\times K}\frac{dxdy}{|x-y|^2}, 
$$
which blows up at the diagonal set $\triangle$ for any knot. 
We normalize this blowing up in the following way: 

Let $\delta _K(x,y)$ denote the shorter arc-length 
between $x$ and $y$ along the knot $K$. 
Define the ``{\it $r^{-2}$-modified $\epsilon$-self avoiding voltage}'' 
at a point $x$, $V_{\epsilon}^{(2)}(K;x)$, and 
the ``{\it $r^{-2}$-modified $\epsilon$-off-diagonal potential energy}'' 
$E_{\epsilon}^{(2)}(K)$ by 
\begin{eqnarray*}
V_{\epsilon}^{(2)}(K;x) &=&\displaystyle{
\int_{y\in K, \delta_K(x,y)\ge\epsilon}
\frac{dy}{|x-y|^{2}}, }\\[2mm]
E_{\epsilon}^{(2)}(K)& =&\displaystyle{\int_{K}V_{\epsilon}^{(2)}(K;x)dx
=\iint _{K\times K, \delta_K(x,y)\ge\epsilon}\frac{dxdy}{|x-y|^{2}}.}
\end{eqnarray*}

Then the order of the blowing up of $V_{\epsilon}^{(2)}(K;x)$ and 
$E_{\epsilon}^{(2)}(K)$, as $\epsilon$ goes down to 0, 
is independent of the knot $K$ and the point $x$. 
Therefore we have the limits: 

\begin{eqnarray}
V^{(2)}(K;x) &=&\displaystyle{
\lim_{\epsilon\to 0}
\left(\int_{y\in K, \delta_K(x,y)\ge\epsilon}
\frac{dy}{|x-y|^{2}}-\frac 2{\epsilon }\right)+4, }
\nonumber \\[3mm]
E^{(2)}(K) &=&\displaystyle{\lim_{\epsilon\to 0}\left(
\iint_{K\times K, \delta_K(x,y)\ge\epsilon}\frac{dxdy}{|x-y|^{2}}
-\frac 2{\epsilon }\right)+4. }
\label{defofe}
\end{eqnarray}

Then $E^{(2)}(K)$ is equal to 
$$
E^{(2)}(K)=\iint _{K\times K}\left(\frac1{|x-y|^{2}}
-\frac1{\delta _K(x,y)^{2}}\right)dxdy. 
$$
We remark that the integrand converges near the diagonal set 
for knots of class $C^4$. 
This expression
\footnote{The idea of using the arc-length as the subtracting counter term 
to cancel the blowing up of the integral at the diagonal set 
was introduced by Nakauchi \cite{Na}. }
 works for open knots and (open) knots which are not necessarily parametrized 
by arc-length. 

Then $E^{(2)}$ is an energy functional for knots. In fact 
we can take a constant $C(b)$ which is independent of $\delta$ as 
$C(b,\delta)$ in Definition \ref{defofenergyfunctional}. 

\medskip
Freedman, He and Wang showed that $E^{(2)}$ is conformally invariant in 
\cite{F-H-W}. 
After their proof, 
Doyle and Schramm obtained a new formula for $E^{(2)}$ 
from a geometric interpretation, 
which implies a simpler proof. 
We introduce their formula in the next subsection. 

\begin{rema} 
If the power ``2'' in the formula of $E^{(2)}$ is replaced by 
$\alpha$ with $\alpha \ge 1$, 
then $E^{(\alpha)}$ can be defined similarly if $\alpha <3$. 
$E^{(\alpha)}$ thus defined is an energy functional for knots if and only if 
$\alpha \ge 2$ and is conformally invariant if and only if $\alpha =2$. 
\end{rema}

\medskip
The jump of $E^{(2)}$ for non-trivial knots is shown as follows. 

Freedman and He defined the {\it average crossing number} $ac(K)$ 
of a knot 
as the average of the numbers of the crossing points 
of the projected knot diagrams, where the average is 
taken over all the directions of the projection (\cite{Fr-He}). 

Then Freedman, He, and Wang showed that 
$E^{(2)}$ bounds the average crossing number from above (\cite{F-H-W}), 
which implies the infimum of the value of $E^{(2)}$ for non-trivial knots, 
because the average crossing number of a  non-trivial knot 
is greater than or equal to 3. 

\medskip
\subsection{The cosine formula for $E^{(2)}$. \label{subsectioncosineformula}}

In this subsection we introduce the cosine formula for $E^{(2)}$ 
invented by Doyle and Schramm (see \cite{Au-Sa}). 

Let us first recall a couple of formulae in the conformal geometry. 
Let $T$ be a conformal transformation of $\RR^3\cup\{\infty\}$. 
Put 
$$
|T^{\prime}(p)| =|\det dT(p)| ^{\frac16}
$$
for $P\in\RR^3$. 
In particular, if $T$ is a reflection with respect to a sphere of radius 
$r$ with center $0$ then $|T^{\prime}(p)| =\frac r{|p|}$. 
Then 
$$
|T(p)-T(q)|=|T^{\prime}(p)||T^{\prime}(q)||p-q|
$$
for $p,q\in\RR^3$ and 
$$
|(T\circ f)^{\prime}(s)|=|T^{\prime}(f(s))|^2|f^{\prime}(s)|
$$
for $f:S^1 \, \mbox{or} \, \RR\to\RR^3$. 

Let $I_x$ be the reflection 
with respect to the sphere 
with center $x$ and radius 1, 
$\tilde{K}_x=I_x(K)$, 
$\tilde{y}=I_x(y)$, $p_{\pm}(\epsilon )=f(s\pm \epsilon )$, and 
$\tilde{p}_{\pm}(\epsilon )=I_x(p_{\pm}(\epsilon ))$ for $0<\epsilon \ll 1$. 
We call $\tilde{K}_x$ the {\it inverted open knot} of $K$ at $x$. 

\begin{figure}[ht]
\centerline{\input{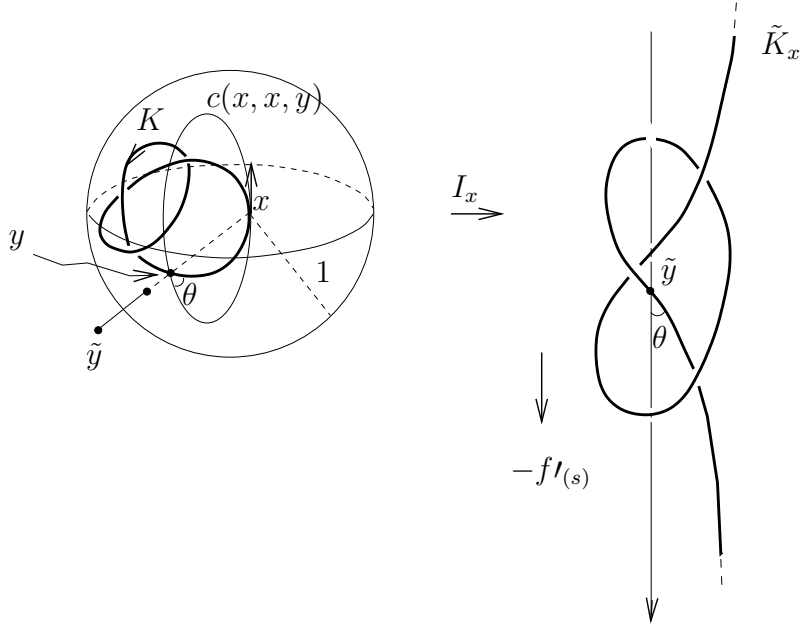}}
\caption{The inverted open knot $\tilde{K}_x$ and the angle $\theta$. }    
\end{figure}

\noindent
Then 
$$
|{I_x}^{\prime}(y)|=\frac1{|x-y|}
$$
and 
$$
|d\tilde{y}| =|d(I_x(y))|
 =|{I_x}^{\prime}(y)|^2|dy|
 =\frac{|dy|}{|x-y|^2}.
$$
For us $dx$ and $dy$ are infinitesimal arcs on $K\subset  \RR^3 \, \mbox{or} \, S^3$, and $|dx|$, $|dy|$, their lengths. We may, later, when the context is clear, use also $dx$ and $dy$ as two one-forms on $K\times K \setminus \Delta$.

Therefore 
$$
V_{\epsilon}^{(2)}(K;x)=
\int_{y\in K, \delta (x,y)\ge \epsilon}\frac{dy}{|x-y|^2}
=\int_{\tilde{p}_+(\epsilon )}^{\tilde{p}_-(\epsilon )}|d\tilde{y}|
$$
is the arc-length between $\tilde{p}_+(\epsilon )$ 
and $\tilde{p}_-(\epsilon )$ 
along the inverted knot $\tilde{K}_x$. 

Let $\bar t=f^{\prime}(s)$ be the unit tangent vector, 
$n=f^{\prime\prime}(s)/|f^{\prime\prime}(s)|$ 
be the unit principal normal vector, 
and $k=|f^{\prime\prime}(s)|$ 
be the curvature of $K$ at $x$. 
Put $n=0$ when $|f^{\prime\prime}(s)|=0$. 
Then 
$$
p_{\pm}(\epsilon ) =x\pm \epsilon \bar t+
\frac{\epsilon ^2k}2n+O(\epsilon ^3), 
$$
which implies 
$$
\tilde{p}_{\pm}(\epsilon ) 
=x\pm \frac1{\epsilon} \bar t+\frac{k}2n+O(\epsilon ).
$$
Hence 
$$
|\tilde{p}_+(\epsilon )-\tilde{p}_-(\epsilon )|=\frac2{\epsilon}+O(\epsilon ). 
$$
Therefore
$$
\int_{y\in K, \delta (x,y)\ge \epsilon}\frac{dy}{|x-y|^2}-\frac2{\epsilon}
=\int_{\tilde{p}_+(\epsilon )}^{\tilde{p}_-(\epsilon )}|d\tilde{y}|-
|\tilde{p}_+(\epsilon )-\tilde{p}_-(\epsilon )|+O(\epsilon )
$$
is equal to the difference of 
the arc-length along $\tilde{K}_x$ and the distance 
between $\tilde{p}_+(\epsilon )$ and $\tilde{p}_-(\epsilon )$ 
up to $O(\epsilon )$. 
(This was called the {\it ``wasted length"} of the inverted open knot.) 

\begin{defi} \rm \label{x_y_z_circle}
Let $C(x,x,y)$ denote the circle 
\footnote{In general $C(x,y,z)$ will denote the uniquely determined 
circle that passes through $x, y$ and $z$ after the notation of \cite{Go-Ma}.}
tangent to a knot $K$ at $x$ 
that passes through $y$	
with the natural orientation 
derived from that of $K$ at point $x$, 
and let $\theta$ $(0\le\theta\le\pi)$ be the angle between $C(x,x,y)$ 
and $C(y,y,x)$ at $x$ or at $y$. 
\end{defi}

We call $\theta =\theta_K (x,y)$ the {\it conformal angle}.

Put $\tilde C(\infty, \infty, \tilde y)
=I_x(C(x,x,y))$. 
Then $\tilde C(\infty, \infty, \tilde y)$ is the line 
which passes through $\tilde y$ and has the tangent vector $-f^{\prime}(s)$. 
The line passing through $\tilde{p}_+(\epsilon)$ and $\tilde{p}_-(\epsilon)$ 
approaches parallel as $\epsilon $ goes down to 0 to the line $\tilde C(\infty, \infty, \tilde y)$. 
Since $I_x$ is a conformal map the angle between the line $\tilde C(\infty, \infty, \tilde y)$
and the tangent line of $\tilde{K}_x$ at $\tilde{y}$ is equal to $\theta$. 
Therefore
\begin{eqnarray}
V^{(2)}(K;x)-4 &=&\displaystyle{\lim_{\epsilon \to 0}\left(\int_{y\in K, \delta (x,y)\ge \epsilon}
\frac{dy}{|x-y|^2}-\frac2{\epsilon}\right)}\nonumber \\[3mm]
 &=&\displaystyle{\lim_{\epsilon \to 0}\int_{\tilde{p}_+(\epsilon )}^{\tilde{p}_-(\epsilon )}
d\tilde{y}(1-\cos \theta )
 =\lim_{\epsilon \to 0}\int _{y\in K}\frac{(1-\cos \theta)dy}{|x-y|^2}.} \nonumber
\end{eqnarray}
Hence 
\begin{equation}
E^{(2)}(K)= \iint  _{K\times K}\frac{1-\cos \theta}{|x-y|^2}dxdy+4, 
\label{cosineformula}
\end{equation}
which is called the {\it cosine formula} for $E^{(2)}$ by Doyle and Schramm.

\section{The infinitesimal cross-ratio. }

In this section we introduce an infinitesimal interpretation 
of the 2-form 
$$
\frac{dxdy}{|x-y|^2}
$$
in the integrand of $E^{(2)}$ in terms of the cross-ratio. 

\medskip
\subsection{The infinitesimal cross-ratio as a ``bilocal'' function.}

In what follows in this section, 
we use the following notation: 
$$
K=f(S^1), \, x=f(s), \, x+dx=f(s+ds), \, 
y=f(t), \, \>\mbox{{\rm and}}\> \, y+dy=f(t+dt). 
$$
As is stated in subsection \ref{ss4tuple}, 
we can define the cross-ratio $(x+dx, y; x, y+dy)$ of 
the ordered four points $x, x+dx, y, y+dy$ 
via the oriented 2-sphere $\Sigma (x, x+dx, y, y+dy)$ that 
passes through them which is  uniquely defined unless the four points are concircular. 

\begin{defi}
\rm 
A {\it twice tangent sphere} $\Sigma_K(x,y)$ is an oriented 2-sphere that is 
tangent to the knot $K$ at the points $x$ and $y$ obtained as the limit of 
the sphere $\Sigma (x, x+dx, y, y+dy)$ as $dx$ and $dy$ go to 0. 
\end{defi}

\begin{figure}[ht] \label{twicetangent}
\centerline{\input{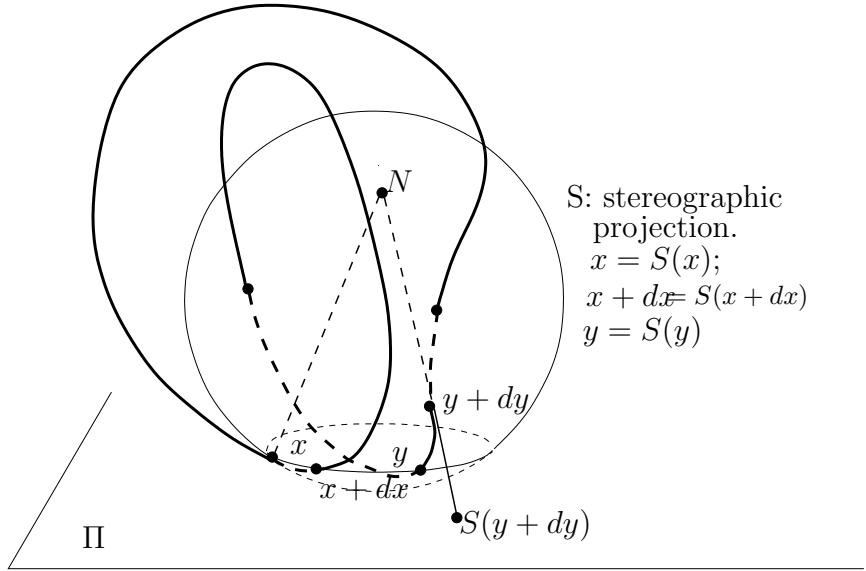}}
\caption{The twice tangent sphere and the infinitesimal cross-ratio.}    
\end{figure}

It is uniquely determined unless $f(s)$, 
$f^{\prime}(s)$, $f(t)$,  and $f^{\prime}(t)$ are concircular. 
We will give a formula for the twice tangent spheres in the next subsection. 

Let 
$\Pi =\Pi_x(y)$ be the plane which contains $C(x,x,y)$
with an identification with the complex plane $\CC $, and let 
$$
S:\Sigma_K(x,y)\to\Pi_x(y)\cup\{\infty\}\cong\CC \cup\{\infty\}
$$
be the orientation preserving stereographic projection. 
When $f(s)$, $f^{\prime}(s)$, $f(t)$,  and $f^{\prime}(t)$ are concircular, 
we can take $\Pi$ with any orientation as $\Sigma _K$. 
Let 
$\tilde x=S(x)$, $\tilde x+\widetilde{dx}=S (x+dx)$, 
$\tilde y=S(y)$, and $\tilde y+\widetilde{dy}=S(y+dy)$
denote the corresponding complex numbers. 
Then the cross ratio $(x+dx, y; x, y+dy)$ is given by\footnote{R.~Kusner already mentionned such an infinitesimal cross-ratio in joint AMS-SMM
meeting in Oaxaca, Mexico in the fall of 1997 and in a talk on quadrupoles
in Illinois (1998).}  
$$
\frac{(\tilde x+\widetilde{dx})-\tilde x}{(\tilde x+\widetilde{dx})
-(\tilde y+\widetilde{dy})}
:\frac{\tilde y-\tilde x}
{\tilde y-(\tilde y+\widetilde{dy})}
\sim\frac{\widetilde{dx}\widetilde{dy}}
{(\tilde x-\tilde y)^2}. 
$$
Let $v_x(y)$ be the unit tangent vector of $C(x,x,y)$ at $y$. 

\begin{defi} \rm  
We call $v_x(y)$ the {\it unit tangent vector at $x$ 
conformally translated to $y$}. 
\end{defi} 

Let $\theta _K(x,y)$ $(-\pi\le\theta _K(x,y)\le\pi)$ 
be the angle from $v_x(y)$ to $f^{\prime}(t)$, 
where the sign of $\theta _K(x,y)$ is given according to 
the orientation of $\Sigma _K(x,y)$. 
Then $|\theta _K(x,y)|$ is equal to the angle $\theta$ 
between $C(x,x,y)$ and $C(y,y,x)$ in the cosine formula. 
As $x, x+dx$ and $y$ are in $\Sigma_K(x,y)\cap\Pi$, 
they are invariant by $p$.

\begin{figure}[ht] \label{infcr}
\centerline{\input{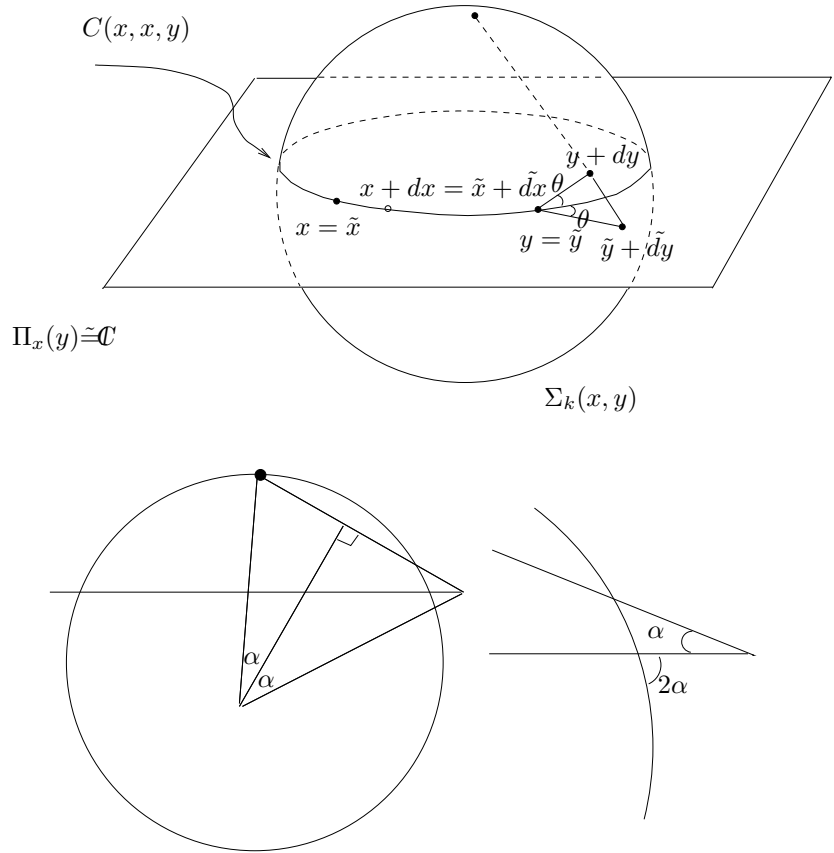}}
\caption{The twice tangent sphere, stereographic projection and computation of the infinitesimal cross-ratio.}    
\end{figure}

Since $|\widetilde{dy}|=|dy|$ as illustrated on the Figures \ref{infcr}, which consider the case when the argument of the cross-ratio is negative, 
the cross-ratio, $(x+dx, y; x, y+dy)$ 
has the absolute value : 
$$
\frac{|dx||dy|}{|x-y|^2},
$$
and the argument: $\theta _K(x,y)$. 

\begin{defi} \rm  
We call $(x+dx, y; x, y+dy)$
the {\it infinitesimal cross-ratio} of $dx$ and $dy$ and denote it 
by $\Omega_{CR}$. 
\label{defofinfcr}
\end{defi} 

\medskip

A stereographic projection from a sphere 
is a restriction of a conformal transformation of $\RR^3\cup\{\infty\}$. 
Hence 
the infinitesimal cross-ratio (or its is complex conjugate class) is 
the unique invariant of $dx$ and $dy$ 
under the action of the orientation preserving conformal group ${\mathcal G}_+$ 
(or respectively, the action of the conformal group 
${\mathcal G}$). 

Since the absolute value of the 2-form 
$$
\frac{|dx||dy|}{|x-y|^2}
$$
on $K\times K\setminus\triangle$ 
is the absolute value of $\Omega_{CR}$ 
and the angle $\theta$ $(0\le\theta\le\pi)$ is 
the absolute value of the argument of $\Omega_{CR}$, 
the cosine formula (\ref{cosineformula}) 
implies 
\begin{equation}
E^{(2)}(K)=\iint_{K\times K\setminus\triangle}
\{|\Omega_{CR}|-\R \Omega_{CR}\}+4. 
\label{cosineformulawithinfinitesimalcrossratio}
\end{equation}
Therefore it is conformally invariant. 

We close this subsection with the following claim which we will use later. 

\begin{clai} The argument of cross-ratio of a knot $K$ is identically 0 
if and only if $K$ is the standard circle. 
\label{arg0}
\end{clai}

\begin{demo} Fix a point $x\in K$ and consider the inverted open knot 
at $x$, $\tilde K_x$. 
Since $\theta (x,y)=0$ for any $y\in K$, 
the tangent vector of $\tilde K_x$ at $\tilde y$ is equal to $-f^{\prime}(s)$ 
for any $\tilde y\in\tilde K_x$. 
Therefore $\tilde K_x$ is the straight line, which means that $K$ is the standard circle. 
\end{demo}

\medskip 
\subsection{The 4-tuple map for a knot and the twice tangent spheres.}


In this subsection we define a map that assigns 
the oriented 2-sphere that passes through 
a given 4-tuples of points on a knot, and then, 
a map that assigns the twice tangent sphere. 

\begin{defi} \rm 
The {\it open concircular point set} ${\mathcal C}c(K)$ {\it of a knot} $K$ 
is given by ${\mathcal C}c(K)=(f^4)^{-1}({\mathcal C}c(S^3))$, 
where $f^4=f\times f\times f\times f:{\rm Conf}_4(S^1)\to {\rm Conf}_4(S^3)$ 
is the
natural map.
\end{defi}

Suppose $K$ is a non-trivial knot. 
G. Kuperberg showed that $K$ contains 4 points 
which are collinear (\cite{Ku}). 
As a corollary, ${\mathcal C}c(K)$ is not empty. 

\begin{ques} 
Is ${\mathcal C}c(K)$ not empty for any trivial knot $K$?
\end{ques}

Moreover we will show later that 
${\mathcal C}c(K)$ is generically of dimension 2. 
This is indicated by the following observation: 

Let $K_p=I_{\Sigma_p}(K)$ be the inverted closed knot 
with respect to a sphere with center $p\in\RR^3\setminus K$. 
Then G. Kuperberg's theorem implies that there is a line $L_p$ 
satisfying that $L_p\cap K_p$, and hence $I_p(L_p)\cap K$, 
contains 4 points. 
This means that for any point $p\in\RR^3\setminus K$ there is a circle 
$C_p$ that intersects $K$ in at least 4 points. 
This implies that, by considering the normal plane to $C_p$ at $p$, 
the set of such circles, and hence ${\mathcal C}c(K)$ has dimension 2. 

Assume that 
$$
f(t_i)=(f_1(t_i), f_2(t_i), f_3(t_i), f_4(t_i), 1)\in S^3_1=\{L=0\}\cap\{x_5=1\}\subset \RR^5 
$$
for $i=1,2,3,4$.  
Put 
$$
\nu _f 
=\nu _f(t_1, t_2, t_3, t_4)
=(\nu_1, \cdots , \nu_5)
=f(t_1)\wedge f(t_2)\wedge f(t_3)\wedge  f(t_4). 
$$ 

We remark that $\nu_5(t_1, t_2, t_3\, t_4)= 0$ 
if and oly if 
$f(t_i)$'s are linearly dependent as vectors in $\RR^5$, 
which means that 
$f(t_i)$'s are concircular as vectors in $S^3_1\cong S^3$. 

\begin{defi}\rm 
The {\it 4-tuple map $\psi_K$ for a knot $K$} 
is a map which assigns 
the oriented 2-sphere that passes through $f(t_i)$ $(i=1,2,3,4)$ 
to a 4-tuple 
$(t_1, t_2, t_3, t_4)\in {\rm Conf}_4(S^1)\setminus {\mathcal C}c(K)$, 
which is given by 
$$
\psi_K=
\Sigma \circ f^4 :{\rm Conf}_4(S^1)\setminus {\mathcal C}c(K)
\ni (t_1, t_2, t_3, t_4)\mapsto \hskip 1cm {}
$$
$$
{} \hskip 1cm 
\frac{\nu _f}{\sqrt{L(\nu _f)}}=
\frac{f(t_1)\wedge f(t_2)\wedge f(t_3) \wedge f(t_4)}
{\sqrt{L(f(t_1)\wedge f(t_2)\wedge f(t_3) \wedge f(t_4))}}
\in\Lambda\cong {\mathcal S}. 
$$
\end{defi}

\begin{clai} 
The open concircular point set 
${\mathcal C}c(K)$ is of dimension more than or equal to 2 
if it is not an empty set, especially if $K$ is 
a non-trivial knot. 
\end{clai}

\begin{demo} Assume that $(t_1, t_2, t_3, t_4)\in{\mathcal C}c(K)\in {\rm Conf}_4$. 
Then, considered as vectors in $\RR^5$, 
$f(t_1), f(t_2), f(t_3)$, and $f(t_4)$ are linearly dependent, 
but three of them, say, $f(t_1), f(t_2), \, and \,  f(t_3)$ 
are linearly independent. 
Suppose 
$$
f(t_4)=af(t_1)+bf(t_2)+cf(t_3).
$$
Then 
$$
\frac{\partial \nu_f}{\partial t_1}\, =-af(t_1)\wedge f(t_2)\wedge f(t_3)
\wedge f^{\prime}(t_1)
$$
$$
\frac{\partial \nu_f}{\partial t_2}\, =-bf(t_1)\wedge f(t_2)\wedge f(t_3)
\wedge f^{\prime}(t_2)
$$
$$
\frac{\partial \nu_f}{\partial t_3}\, =-cf(t_1)\wedge f(t_2)\wedge f(t_3)
\wedge f^{\prime}(t_3)
$$
$$
\frac{\partial \nu_f}{\partial t_4}\, =\hskip 0.4cm 
f(t_1)\wedge f(t_2)\wedge f(t_3)
\wedge f^{\prime}(t_4).
$$
Since $f(t_1), f(t_2), \,  and \,  f(t_3)$ are linearly inde\-pendent, 
$$
Span\left\langle\frac{\partial \nu_f}{\partial t_1}, 
\frac{\partial \nu_f}{\partial t_2}, 
\frac{\partial \nu_f}{\partial t_3}, 
\frac{\partial \nu_f}{\partial t_4}
\right\rangle
$$
has dimension at most 2. 
Therefore the kernel of $d\nu_f$, which is equal to 
the tangent space of ${\mathcal C}c(K)$, has dimension at least 2. 
\end{demo}

\medskip
Let us consider the behavior of $\psi_K$ when 
$t_2$ and $t_4$ approach $t_1$ and $t_3$ respectively. 
Then Taylor's expansion formula implies
\begin{eqnarray*}
f(t_1)\wedge f(t_2)\wedge f(t_3)\wedge f(t_4) 
&=&(t_2-t_1)(t_4-t_3)
f(t_1)\wedge f(t^{\prime}_1)\wedge f(t_3)\wedge f^{\prime}(t_3)\\[2mm]
&&+\> \mbox{higher order terms}. 
\end{eqnarray*}
Put 
$\nu _f^{(2,2)} (t_1, t_3)
= f(t_1)\wedge f(t^{\prime}_1)\wedge f(t_3)\wedge f^{\prime}(t_3).$ 
Then $\nu _f^{(2,2)} (t_1, t_3)=0$ if and only if 
$f(t_1)$, $f(t^{\prime}_1)$, $f(t_3)$, and $f^{\prime}(t_3)$ are 
linearly dependent as vectors in $\RR^5$, 
which occurs if and only if 
$f^{\prime}(t_3)$ can be expressed as 
a linear combination of $f(t_1)-f(t_3)$ and $f(t^{\prime}_1)$. 
This is equivalent to the condition that $f^{\prime}(t_3)$, considered as 
a tangent vector at $f(t_3)\in S^3_1$, lies in 
the circle which can be obtained as 
the intersection of $S^3_1$ and the 3-dimensional vector subspace 
spanned by $f(t_1)-f(t_3)$ and $f(t^{\prime}_1)$. 
Therefore $\nu _f^{(2,2)} (t_1, t_3)=0$ if and only if 
$f(t_1)$, $f^{\prime}(t_1)$, $f(t_3)$, and $f^{\prime}(t_3)$ are concircular. 

Put 
$$
{\mathcal C}c^{(2,2)}(K)=\{(s,t)\in S^1\times S^1\setminus\triangle \,|\, 
f(t_1), f^{\prime}(t_1), f(t_3), \> {\rm and} \> f^{\prime}(t_3) \> 
\mbox{are concircular}\}. 
$$
Generically ${\mathcal C}c^{(2,2)}(K)$ is of dimension 0 
since $(x,y)=(f(s), f(t))$ belongs to ${\mathcal C}c^{(2,2)}(K)$ 
if and only if $f^{\prime}(t)$ coincides with $v_x(y)$, 
the unit tangent vector at $x$ conformally translated to $y$. 
We remark that ${\mathcal C}c^{(2,2)}(K)$ has dimension 1 
when $K$ is a $(p,q)$-torus knot 
which is an orbit of an $S^1$-action defined by 
$$
S^1\ni e^{2\pi it}\mapsto 
\left(\CC^2\supset S^3\ni (z,w)\mapsto (e^{2\pi ipt}z, e^{2\pi iqt}w)
\in S^3\subset \CC^2\right). 
$$

\begin{defi} \rm 
The {\it twice tangent sphere map} $\psi_K^{(2,2)}$ for a knot $K$ 
which assigns the twice tangent sphere $\Sigma_K(f(s),f(t))$ 
to $(s,t)\in {\rm Conf}_2(S^1)\setminus {\mathcal C}c^{(2,2)}(K)$ 
is given by 
$$
\psi_K^{(2,2)}:{\rm Conf}_2(S^1)\setminus {\mathcal C}c^{(2,2)}(K)\ni (s,t)
\mapsto \hskip 3.5cm {}
$$
$$
{}\hskip 1.5cm
\frac{\nu _f^{(2,2)} (s,t)}{\sqrt{L(\nu _f^{(2,2)} (s,t))}}
=\frac{f(s)\wedge f^{\prime}(s) \wedge f(t)\wedge f^{\prime}(t)}
{\sqrt{L(f(s)\wedge f^{\prime}(s) \wedge f(t)\wedge f^{\prime}(t))}}
\in\Lambda\cong {\mathcal S}. 
$$
\end{defi}

\begin{conj} If $K$ is non-trivial then the twice tangent sphere map $\psi_K^{(2,2)}$ is not injective, 
namely, there is a sphere $\Sigma$ that is tangent to $K$ at three points or more.  
\end{conj}

When $t$ approaches $s$, the twice tangent sphere approaches 
an {\it osculating sphere} that is  
the 2-sphere which is the most tangent to $K$ at $f(s)$. It contains 
the osculating circle and is tangent to $K$ in the fourth order. 
Taylor's expansion formula implies 
$$
\nu _f^{(2,2)} (s,t)=\frac1{12}(t-s)^4
f(s)\wedge f^{\prime}(s)\wedge f^{\prime\prime}(s)
\wedge f^{\prime\prime\prime}(s)
+\> \mbox{higher order terms}. 
$$
Put $\nu _f^{(4)}(s)=f(s)\wedge f^{\prime}(s)\wedge f^{\prime\prime}(s)
\wedge f^{\prime\prime\prime}(s)$. 
Suppose $\nu _f^{(4)}(s)\ne 0$. 
Lemma \ref{lorentzform} (3) shows that $\nu _f^{(4)}(s)$ is 
not a time-like vector. If it is a light-like vector then 
\ref{lorentzform} (2) implies that 
$\nu _f ^{(4)}(s)=kf(s)$ for some $k\in\RR$. 
Since $L(f(t),f(t))= 0$ for any $t\in S^1$ 
$$
L(f(t),f^{\prime\prime}(t))=-L(f^{\prime}(t),f^{\prime}(t))
=\sum_{i=1}^{4}(f_i^{\prime}(t))^2\ne 0 
\hskip 0.4cm {\rm for} \> {\rm any} \hskip 0.4cm t\in S^1. 
$$
On the other hand as $L(\nu _f, f^{\prime\prime}(s))=0$, 
we have $k=0$, which is a contradiction\footnote{We assumed that $|f^{\prime}(t)|$ never vanishes. }.
Therefore when $\nu _f^{(4)}(s)\ne 0$ the osculating sphere 
at $f(s)$ is given by 
$\nu _f^{(4)}(s)/\sqrt{\nu _f^{(4)}(s)}$. 


\subsection{The real part as the canonical symplectic form on $T^{\ast}S^3$. } 

In this subsection we show that the real part of the infinitesimal cross-ratio 
can be interpreted as the pull-back of the canonical symplectic form on 
the cotangent bundle $T^{\ast}S^3$ of the 3-sphere. 
This interpretation allows us to deduce the original definition 
(\ref{defofe}) from the cosine formula in terms of the infinitesimal 
cross-ratio, (\ref{cosineformulawithinfinitesimalcrossratio}). 

\subsubsection{The pull-back of the symplectic form to 
$S^n\times S^n\setminus \triangle$. } 

Recall that a cotangent bundle $\pi :T^{\ast}M\to M$ 
of a manifold $M^n$ admits the canonical symplectic form 
$\omega _0$ given by 
$$
\omega _0=\sum dp_i\wedge dq_i,
$$
where $(p_1, \cdots, p_n)$ is a local coordinate of $M$ and 
$(q_1, \cdots, q_n)$ is the local coordinate of $T^{\ast}M$ 
defined by 
\begin{equation}
T_x^{\ast}M\ni v^{\ast}=\sum q_idp_i. 
\label{f4.3.1}
\end{equation}
We remark that $\omega _0$ is an exact 2-form with 
$\omega _0=d(-\sum q_idp_i). $

\begin{lemm} 

\noindent {\rm (1)} 
There is a canonical bijection 
$$
\psi_n :S^n\times S^n\setminus \triangle \to T^{\ast}S^n. 
$$

\noindent {\rm (2)} Let $(x_1, \cdots, x_{n+1})$ be a system of local coordinates of 
$$
U^+_{n+1}=\{(x_1, \cdots, x_{n+1})\in S^n\subset\RR ^{n+1}\, | \, x_{n+1}>0\}. 
$$
It determines the associated system of local coordinates of 
$\pi^{-1}(U^+_{n+1})\subset T^{\ast}S^n$ by (\ref{f4.3.1}). 
The canonical bijection 
$\psi_n :S^n\times S^n\setminus \triangle \to T^{\ast}S^n$ 
is expressed with respect to this local coordinates as 
$$
\psi _n(x(x_1, \cdots, x_{n+1}), y(y_1, \cdots, y_{n+1}))
=\left(x_1, \cdots, x_n, 
\frac{y_1-\frac{y_{n+1}}{x_{n+1}}x_1}{1-x\cdot y}, \cdots, 
\frac{y_n-\frac{y_{n+1}}{x_{n+1}}x_n}{1-x\cdot y}
\right)
$$
on $U^+_{n+1}\times S^n\setminus \triangle$
, where $x\cdot y$ denotes the inner product. 
\label{coordinate}
\end{lemm}

\begin{demo} (1) 
Let $\Pi _x$ be the $n$-dimensional hyperplane in $\RR^{n+1}$ 
passing through the origin 
that is perpendicular to $x\in S^n$, and 
$p_x:S^n\setminus\{x\}\to\Pi _x$ be the stereographic projection. 
We identify $\Pi_x$ with $T_xS^n$. 
Take  an orthonormal bas $\{v_1, \cdots , v_n\}$ of $\Pi _x\cong T_xS^n$. 
Suppose $y\in S^n\setminus\{x\}$ is expressed as 
$$
y=\tilde y_1v_1+\cdots +\tilde y_nv_n+\tilde y_{n+1}x
$$
with respect to the orthonormal basis $\{v_1, \cdots , v_n, x\}$ 
of $\RR^{n+1}$. 
Then 
$$
p_x(y)=\frac{\tilde y_1}{1-\tilde y_{n+1}}v_1+\cdots +
\frac{\tilde y_n}{1-\tilde y_{n+1}}v_n. 
$$
Put 
\begin{equation}
\psi_x(y)=\frac{\tilde y_1}{1-\tilde y_{n+1}}v_1^{\ast}+\cdots +
\frac{\tilde y_n}{1-\tilde y_{n+1}}v_n^{\ast}\in T_x^{\ast}S^n, 
\label{2.4.2}
\end{equation}
where $\{v_i^{\ast}\}$ is the dual basis of $T_x^{\ast}S^n$. 
Then the map $\psi_x:S^n\setminus\{x\}\to T_x^{\ast}S^n$ 
does not depend on the choice of 
the orthonormal basis $\{v_i\}$ of $\Pi_x\cong T_xS^n$. 
Thus the map 
$$
\psi_n : S^n\times S^n\setminus \triangle \ni (x,y)\mapsto (x,\psi_x(y))
\in T^{\ast}S^n
$$
makes the canonical bijection. 

\end{demo}

\begin{lemm} The pull-back $\omega ={\psi_n}^{\ast}\omega _0$ 
of the canonical symplectic form 
$\omega _0$ of $T^{\ast}S^n$ by 
$\psi_n :S^n\times S^n\setminus \triangle \to T^{\ast}S^n$ 
is given by 
$$
\omega =d\left(-\frac{\sum_{i=1}^{n+1}y_idx_i}{1-x\cdot y}\right)
=d\left(\frac{\sum_{i=1}^{n+1}x_idy_i}{1-x\cdot y}\right)
\hskip 2cm {}
$$
$$
=\frac{\sum_{i=1}^{n+1}dx_i\wedge dy_i}{1-x\cdot y}
+\frac{(\sum_{i=1}^{n+1}y_i dx_i)\wedge 
(\sum_{i=1}^{n+1}x_i dy_i)}{(1-x\cdot y)^2}
$$
\label{4.3.2}
\end{lemm}

\begin{demo}
As $\omega _0=-d\sum q_idp_i$ 
$$
{\psi_n}^{\ast}\omega _0=-d\,\sum_{i=1}^{n}
\frac{y_i-\frac{y_{n+1}}{x_{n+1}}x_i}{1-x\cdot y}dx_i\, 
$$
on $U^+_{n+1}\times S^n\setminus \triangle$ by Lemma \ref{coordinate} (2). 
Since $-\sum_{i=1}^{n}x_idx_i=x_{n+1}dx_{n+1}$ on $S^n$, 
it implies the formula. 
\end{demo}

\begin{lemm} Let $p:S^n\setminus\{(0, \cdots , 0, 1)\}\to\RR^n$ 
be the stereographic projection. 
Put 
$$
P_n=p^{-1}\times p^{-1}:\RR^n\times\RR^n\setminus\triangle\to 
S^n\times S^n\setminus\triangle . 
$$
Then the pull-back 
$\omega _{\RR ^n}={P_n}^{\ast}\omega={P_n}^{\ast}{\psi_n}^{\ast}\omega_0$ 
is given by 
\begin{eqnarray*}
\omega _{\RR ^n}&=&\displaystyle{2d\left(\frac{\sum (x_i-y_i)dy_i}{|x-y|^2}\right)
=2d\left(\frac{\sum (x_i-y_i)dx_i}{|x-y|^2}\right)}\\[2mm]
&=&\displaystyle{2\left(\frac{\sum dx_i\wedge dy_i}{|x-y|^2}
-2\frac{(\sum (x_i-y_i) dx_i)\wedge 
(\sum (x_j-y_j)dy_j)}{|x-y|^4}\right).}
\end{eqnarray*}
\label{euclideansmpl}

\end{lemm}

\begin{demo} Suppose 
$$
P_n(x(x_1, \cdots, x_{n}), y(y_1, \cdots, y_{n}))
=(X(X_1, \cdots, X_{n+1}), Y(Y_1, \cdots, Y_{n+1})). 
$$
Then 
$$
{P_n}^{\ast}\left(\frac{\sum_{i=1}^{n+1}X_idY_i}{1-X\cdot Y}\right)
=2\frac{\sum_{i=1}^{n} (x_i-y_i)dy_i}{|x-y|^2}+d\log (|y|^2+1). 
$$
Therefore 
$$
{P_n}^{\ast}\omega 
={P_n}^{\ast}d\left(\frac{\sum_{i=1}^{n+1}X_idY_i}{1-X\cdot Y}\right)
=2d\left(\frac{\sum_{i=1}^{n} (x_i-y_i)dy_i}{|x-y|^2}\right). 
$$
\end{demo}

\begin{prop} The 2-form $\omega _{\RR ^n}$ is invariant 
under the diagonal action of the conformal group ${\mathcal G}$ on 
$\RR^n \times\RR^n\setminus\triangle$ defined by 
$$
g\cdot (x,y)=(g\cdot x, g\cdot y),  
$$
where $g\in{\mathcal M}$ and $(x,y)\in\RR^n \times\RR^n\setminus\triangle$. 
\label{invdiagonalaction}
\end{prop}

\begin{demo} It is obvious that $\omega _{\RR ^n}$ is invariant 
under the diagonal action of multiplication by scalars $(x,y)\mapsto (cx,cy)$ and of addition of vectors, $(x,y)\mapsto (x+a,y+a)$.
The previous Lemma 
implies that $\omega _{\RR ^n}$ is invariant 
under the diagonal action of the orthogonal group $O(n)$. 

Therefore it suffices to show the invariance under the diagonal action 
$I\times I$ of the inversion $I$ with respect to the $(n-1)$-sphere 
with radius 1 whose center is the origin. 

Suppose 
$$
(I\times I)(x(x_1, \cdots, x_{n}), y(y_1, \cdots, y_{n}))
=(X(X_1, \cdots, X_{n}), Y(Y_1, \cdots, Y_{n})).
$$
Then 
$$
(I\times I)^{\ast}\left(\frac{\sum (X_i-Y_i)dY_i}{|X-Y|^2}\right)
=\left(\frac{\sum (x_i-y_i)dy_i}{|x-y|^2}\right)
+\frac12d\log(|y|^2). 
$$
Therefore 
$$
(I\times I)^{\ast}\omega _{\RR ^n}
=(I\times I)^{\ast}2d\left(\frac{\sum (X_i-Y_i)dY_i}{|X-Y|^2}\right)
=2d\left(\frac{\sum (x_i-y_i)dy_i}{|x-y|^2}\right)
=\omega _{\RR ^n}. 
$$
We give another proof of the invariance. 

Let $\Sigma_1$ be the $n$-sphere in $\RR^{n+1}$ with center $0$ and radius 1, 
and $\Sigma_{\sqrt 2}$ be the $n$-sphere in $\RR^{n+1}$ with center 
$(0, \cdots , 0, 1)$ and radius $\sqrt 2$. 
Let $I_{\Sigma}$ be 
the inversions of $\RR^{n+1}\cup\{\infty\}$ 
with respect to an $n$-sphere $\Sigma$ in $\RR^{n+1}$. 

Then 
$$
p^{-1}(I(x)) = I_{\Sigma_{\sqrt 2}}(I_{\Sigma_1}(x))
=I_{I_{\Sigma_{\sqrt 2}}(\Sigma_1)}(I_{\Sigma_{\sqrt 2}}(x))
=I_{\RR^n}(p^{-1}(x)) 
$$
for $x\in\RR^n$. 
Therefore 
$$
P_n\circ (I\times I)=(I_{\RR^n}\times I_{\RR^n})\circ P_n
:\RR^n\times \RR^n\setminus\triangle \to S^n\times S^n\setminus\triangle . 
$$
As 
$$
\psi_n\circ (I_{\RR^n}\times I_{\RR^n})=\psi_n 
: S^n\times S^n\setminus\triangle \to T^{\ast}S^n, 
$$
we get the conclusion: 
$$
(I\times I)^{\ast}\omega _{\RR ^n}
=(I\times I)^{\ast}{P_n}^{\ast}{\psi_n}^{\ast}\omega _0
=(\psi_n\circ P_n\circ (I\times I))^{\ast}\omega _0
=(\psi_n\circ P_n)^{\ast}\omega _0
=\omega _{\RR ^n}. 
$$
\end{demo}

\begin{lemm} Let 
$$
\lambda=\frac{dw\wedge dz}{(w-z)^2} \, , \hskip 1cm (w,z)
\in\CC\times \CC\setminus\triangle
$$
be the 2-form on 
$\CC\times \CC\setminus\triangle\cong\RR^2\times\RR^2\setminus\triangle$ 
obtained as the infinitesimal cross-ratio of $dw$ and $dz$. Then 

\medskip
{\rm (1)} Both $\R \lambda$ and $\I \lambda$ are exact 2-forms. 

\medskip
{\rm (2)} $\R \lambda$ (or $\I \lambda$) is invariant 
(or invariant up to sign, respectively) 
under the diagonal action of the conformal group\footnote{by conformal group we mean the group of transformations generated by the transformations $z\mapsto \frac{az+b}{cz+d}$ and $z\mapsto {\overline z}$ acting on the riemann sphere $\CC \cup \infty$ (or the restriction of elements of this group  to $\RR^2 \setminus \{ \mbox{the point where it is not defined, if necessary} \}$) } on 
$\RR^2\times\RR^2\setminus\triangle$. 
\end{lemm}

\begin{demo} As for the real part, since the next Lemma implies 
$\R \lambda=-\frac12\omega_{\RR^2}$, 
Lemma \ref{euclideansmpl} and Proposition \ref{invdiagonalaction} show that $\R \lambda$ 
satisfies the desired properties. 

As for the imaginary part, direct calculation shows 
\begin{eqnarray*}
\I \lambda
&=&\displaystyle{-2\frac{(x_1-y_1)(x_2-y_2)(dx_1\wedge dy_1-dx_2\wedge dy_2)}
{\{(x_1-y_1)^2+(x_2-y_2)^2\}^2}}\\[2mm]
&&\displaystyle{+\frac{\{(x_1-y_1)^2-(x_2-y_2)^2\}(dx_1\wedge dy_2+dx_2\wedge dy_1)}
{\{(x_1-y_1)^2+(x_2-y_2)^2\}^2}. }
\end{eqnarray*}

(1) If we put 
$$
\rho =\frac{(x_1-y_1)dx_2-(x_2-y_2)dx_1}
{(x_1-y_1)^2+(x_2-y_2)^2}, 
$$
then $\I \lambda=-d\rho$. 

(2) Let $I$ be the inversion with respect to the circle 
with radius 1 whose center is the origin. Then 
$$
(I\times I)^{\ast}\rho =-\rho+d\arctan (\frac{x_2}{x_1}), 
$$
which implies the invariance as in the proof of Proposition \ref{invdiagonalaction}. 

We remark that 
Lemma can be proved by a direct calculation 
with either complex or real coordinates. 
\end{demo}

\begin{lemm} {\rm { (Folklore)}} The real part of the infinitesimal cross-ratio 2-form 
is equal to minus one half of 
the pull-back of the canonical symplectic form of 
the cotangent bundle $T^{\ast}S^2$: 
$$
\R \left(\frac{dw\wedge dz}{(w-z)^2}\right)
=-\frac12\omega_{\RR^2}=-\frac12{P_2}^{\ast}{\psi_2}^{\ast}\omega_0. 
$$
\end{lemm}

\begin{demo} The left hand side is equal to 
\begin{eqnarray*}
&&\displaystyle{\frac{\{(x_1-y_1)^2-(x_2-y_2)^2\}(dx_1\wedge dy_1-dx_2\wedge dy_2)}
{\{(x_1-y_1)^2+(x_2-y_2)^2\}^2}} \\[2mm]
&&\displaystyle{+2\frac{(x_1-y_1)(x_2-y_2)(dx_1\wedge dy_2+dx_2\wedge dy_1)}
{\{(x_1-y_1)^2+(x_2-y_2)^2\}^2}} \\[2mm]
&=&\displaystyle{-\frac{dx_1\wedge dy_1+dx_2\wedge dy_2}{(x_1-y_1)^2+(x_2-y_2)^2}} \\[2mm]
&&\displaystyle{+2\frac{\{(x_1-y_1)dx_1+(x_2-y_2)dx_2\}\wedge \{(x_1-y_1)dy_1+(x_2-y_2)dy_2\}}
{\{(x_1-y_1)^2+(x_2-y_2)^2\}^2}}
\end{eqnarray*}

which is equal to the right hand side. 
\end{demo}

\subsubsection{From cosine formula to the original definition of $E^{(2)}$. 
\label{s4.3.2}}

Now let us recall the definition of the infinitesimal cross-ratio. 
Let $(x_0, y_0)=(f(s_0), f(t_0))\in K\times K\setminus\triangle$. 
Then there is an oriented twice tangent sphere $\Sigma_K(x_0, y_0)$, 
which is given by $\psi_K^{(2,2)}(s_0, t_0)$ when 
$(s_0, t_0)\in {\rm Conf}_2(S^1)\setminus {\mathcal C}c^{(2,2)}(K)$, 
and otherwise by any oriented 2-sphere that contains the circle $C(x_0, x_0, y_0)$. 
Let $T_0$ be a conformal transformation 
of $\RR^3\cup\{\infty\}$ that maps $\Sigma_K(x_0, y_0)$ to 
$\RR^2\times\{0\}\subset\RR^3$ preserving the orientations. 
Let 
$$
\frac{dw\wedge dz}{(w-z)^2}
$$ 
be a complex valued 2-form on 
$$
\RR^3\times \RR^3\setminus\triangle \cong 
\{((w,x_3),(z,y_3))\in (\CC\times\RR )\times (\CC\times\RR )
\setminus\triangle\}. 
$$
Then, as a 2-form, the infinitesimal cross-ratio $\Omega_{CR}$ satisfies 
$$
\Omega_{CR}(x_0, y_0)=(T_0\times T_0)^{\ast}
\left(\frac{dw\wedge dz}{(w-z)^2}\right)(x_0, y_0). 
$$

\begin{lemm} Let $K=f(S^1)$ be a knot. 
Then the real part of the 
infinitesimal cross-ratio $\Omega_{CR}$ defined in definition \ref{defofinfcr}
is equal to 
minus one half of the pull-back of the canonical symplectic form of 
the cotangent bundle $T^{\ast}S^3$: 
$$
\cos\theta\frac{dxdy}{|x-y|^2}=\R \Omega_{CR}
=-\frac12\omega_{\RR^3}|_{K\times K\setminus\triangle}
=-\frac12(\iota _K\times\iota _K)^{\ast}{P_3}^{\ast}{\psi_3}^{\ast}\omega_0, 
$$
where $\iota _K:K\to\RR^3$ is the inclusion map. 
\label{realcr=sympl}
\end{lemm}

\begin{demo} 
We show that their pull-backs $(f\times f)^{\ast}\R \Omega_{CR}$ and 
$-\frac12(f\times f)^{\ast}\omega_{\RR^3}$ 
coincide at any given $(s_0, t_0)\in S^1\times S^1\setminus \triangle$. 

Proposition \ref{invdiagonalaction} implies 
$$
-\frac12(f\times f)^{\ast}\omega_{\RR^3}
=-\frac12((T_0\circ f)\times (T_0\circ f))^{\ast}\omega_{\RR^3}. 
$$
On the other hand, since 
$$
\left((T_0\circ f)^{\ast}dx_3\right)(s_0, t_0)
=\left(\frac{d}{ds}(T_0\circ f)_3\right)(s_0), 
$$
$$
\left((T_0\circ f)^{\ast}dy_3\right)(s_0, t_0)
=\left(\frac{d}{dt}(T_0\circ f)_3\right)(t_0)  {}
$$
we have 

\begin{eqnarray*}
&&\displaystyle{-\frac12((T_0\circ f)\times (T_0\circ f))^{\ast}\omega_{\RR^3}(s_0, t_0)}\\[2mm]
&=&\displaystyle{((T_0\circ f)\times (T_0\circ f))^{\ast}
\left(-\frac{\sum _{i=1}^3dx_i\wedge dy_i}{|x-y|^2}
\right.}\\[2mm]
&&\hskip 3.5cm 
\displaystyle{\left. +2\frac{(\sum _{i=1}^3(x_i-y_i)\wedge dx_i)\wedge 
(\sum _{i=1}^3(x_i-y_i)\wedge dy_i)}{|x-y|^4}\right)(s_0, t_0)}\\[2mm]
&=&\displaystyle{((T_0\circ f)\times (T_0\circ f))^{\ast}
\left(-\frac{\sum _{i=1}^2dx_i\wedge dy_i}{|x-y|^2}
\right.}\\[2mm]
&&\hskip 3.5cm 
\displaystyle{\left.+2\frac{(\sum _{i=1}^2(x_i-y_i)\wedge dx_i)\wedge 
(\sum _{i=1}^2(x_i-y_i)\wedge dy_i)}{|x-y|^4}\right)(s_0, t_0)}\\[2mm]
&=&\displaystyle{(f\times f)^{\ast}(T_0\times T_0)^{\ast}
\Re\left(\frac{dw\wedge dz}{(w-z)^2}\right)(s_0, t_0)}
=(f\times f)^{\ast}\R \Omega_{CR}(s_0, t_0). 
\end{eqnarray*}
\end{demo}

Put for $0<\epsilon \ll 1$ 
$$
N(\epsilon )=\{(s,t)\in S^1\times S^1\, | \,\epsilon\le |s-t|\le 1-\epsilon\}, 
$$
where $S^1=[0,1]/\sim$, and assume $N(\epsilon )$ has a natural orientation 
as a subspace of the torus. 
Then its boundary $\partial N(\epsilon )$ consists of two closed curves 
$L_+$ and $L_-$ with opposite orientations, where 
$$
L_{\pm}=\{(s,t)\in S^1\times S^1\, |\, t-s=\pm \epsilon \>\> 
({\rm mod} \>\> \ZZ)\}. 
$$
Then 
\begin{eqnarray*}
E^{(2)}(K)&=&\displaystyle{\iint_{S^1\times S^1\setminus\triangle}
(f\times f)^{\ast}\left(\frac{dx\, dy}{|x-y|^2}-\R \Omega_{CR}\right)+4}\\[2mm]
&=&\displaystyle{\lim_{\epsilon \to 0}\iint_{N(\epsilon )}
\left(
\frac{ds\, dt}{|f(s)-f(t)|^2}-(f\times f)^{\ast}\R \Omega_{CR}
\right)+4. }
\end{eqnarray*}

On the other hand, as 
$$
-(f\times f)^{\ast}\R \Omega_{CR}=\frac12(f\times f)^{\ast}\omega _{\RR^3}
=d(f\times f)^{\ast}\left(\frac{\sum (x_i-y_i)dy_i}{|x-y|^2}\right)
$$
by Lemmas \ref{realcr=sympl} and \ref{euclideansmpl}, there holds 
$$
E^{(2)}(K)=\lim_{\epsilon \to 0}\left\{
E^{(2)}_{\epsilon}(K) + \int_{L_+\cup L_-}
(f\times f)^{\ast}\left(\frac{\sum (x_i-y_i)dy_i}{|x-y|^2}\right)
\right\}+4
$$
by Stokes' theorem. 
As 
$$
(f(s)-f(t), f^{\prime}(t))=(s-t)+O(|s-t|^3), 
$$
we have 
$$
\int_{L_+\cup L_-}
(f\times f)^{\ast}\left(\frac{\sum (x_i-y_i)dy_i}{|x-y|^2}\right)
=\int_{L_+\cup L_-}
\frac{(f(s)-f(t), f^{\prime}(t))}{|f(s)-f(t)|^2}dt
$$
$$
=\int_0^1\frac{-\epsilon +O(\epsilon ^3)}{\epsilon ^2+O(\epsilon ^4)}dt
+\int_1^0\frac{\epsilon +O(\epsilon ^3)}{\epsilon ^2+O(\epsilon ^4)}dt
=-\frac2{\epsilon}+O(\epsilon). 
$$
Therefore we restored the original definition (\ref{defofe}): 
$$
E^{(2)}(K)=\lim_{\epsilon \to 0}\left(
E^{(2)}_{\epsilon}(K)-\frac2{\epsilon}
\right)+4.
$$


\medskip
\section{The (absolute) imaginary cross-ratio energy, $E_{|\sin\theta|}$.}

In the previous section we introduced the infinitesimal 
cross-ratio $\Omega_{CR}$, which is the unique conformal invariant 
of a pair of infinitesimal curve segments $dx$ and $dy$,  
and showed that $E^{(2)}$ can be expressed 
in terms of the absolute value and the real part of if. 
On the other hand, one can use the imaginary part of $\Omega_{CR}$. 

\subsection{The projection of the inverted open knot.} \label{proj_invert}

\begin{defi} \rm  The {\it absolute imaginary cross-ratio energy} 
or the {\it conformal $|\sin \theta |$ energy}, 
$E_{|\sin\theta|}$, is defined by 
$$
E_{|\sin\theta|}(K) =\iint _{K\times K\setminus\triangle}
|\I \Omega_{CR}|
 =\iint _{K\times K\setminus\triangle}\frac{|\sin\theta |dxdy}{|x-y|^2}
$$
\end{defi} 
$E_{|\sin\theta|}$ is well-defined 
since the argument of cross-ratio, $\theta$, 
is of the same order as $|x-y|^2$ 
near the diagonal set by the following Lemma. 

The interest of the $E_{|\sin \theta |}$ energy has already been noticed by Kusner and Sullivan \cite{Ku-Su}, section 4. 

\begin{rema} 
The lemma \ref{sign_cross} shows that the sign of $sin\theta$ is constant, the absolute value allows us to forget it.
\end{rema}  

\begin{lemm} If a knot is of class $C^4$, 
the conformal angle $\theta _K(x,y)$ satisfies 
$$
\theta _K(x,y)=\frac{\sqrt{k^2\tau ^2+{k^{\prime}}^2}}6|x-y|^2+O(|x-y|^3),
$$
near the diagonal set, where $k=|f^{\prime\prime}|$ is the curvature 
and $\tau$ is the torsion of the knot at $x$. 
\label{orderofarg}
\end{lemm}

\begin{demo} 
Let us first remark that 
the unit tangent vector at $x$ 
conformally translated to $y$, $v_x(y)$, is expressed as 
$$
v_x(y)=2\left(f^{\prime}(s),
\frac{f(t)-f(s)}{|f(t)-f(s)|}\right)\cdot\frac{f(t)-f(s)}{|f(t)-f(s)|}
-f^{\prime}(s). 
$$
since $v_x(y)$ is symmetric to $f^{\prime}(s)$ 
with respect to $f(t)-f(s)$. 

Assume that a knot $f$ is of class $C^4$ and parametrized by the arc-length. 
Put $s=0$ and take the Frenet frame at $x$. 
Then by Bouquet's formula $f(t)$ is expressed as 
$$
f_1(t) =t-\frac{k^2}6t^3+\cdots, \quad 
f_2(t) =\frac{k}2t^2+\frac{k^{\prime}}6t^3+\cdots, \quad 
f_3(t) =\frac{k\tau}6t^3+\cdots. 
$$
Then as $f^{\prime}(0)=(0,0,1)$ 
the unit tangent vector at $x$ 
conformally translated to $y$ is given by 
\begin{eqnarray*}
v_x(y) &=&\displaystyle{\frac{(f_1{}^2(t)-f_2{}^2(t)-f_3{}^2(t), 2f_1(t)f_2(t), 
2f_1(t)f_3(t))}{|f(t)|^2}}\\[2mm]
&=&\displaystyle{\frac{\left(t^2-\frac{7k^2}{12}t^4+O(t^5), 
kt^3+\frac{k^{\prime}}{3}t^4+O(t^5), 
\frac{k\tau}3t^4+O(t^5)\right)}
{t^2\left(1-\frac{k^2}{12}t^2+O(t^3)\right)}.}
\end{eqnarray*}

On the other hand 
$$
\frac{f^{\prime}(t)}{|f^{\prime}(t)|}=\frac{\left(
1-\frac{k^2}{2}t^2+O(t^3), kt+\frac{k^{\prime}}{2}t^2+O(t^3), 
\frac{k\tau}2t^2+O(t^3)
\right)}{1+O(t^3)}. 
$$
Hence
\begin{eqnarray*}
|\sin\theta| &=&\displaystyle{\left|\frac{f^{\prime}(t)}{|f^{\prime}(t)|}\times v_x(y)\right|
 =\frac{\left|\left(O(t^5), \frac{k\tau}{6}t^4+O(t^5), 
-\frac{k^{\prime}}{6}t^4+O(t^5)\right)\right|}
{t^2\left(1-\frac{k^2}{12}t^2+O(t^3)\right)(1+O(t^3))}}\\[2mm]
 &=&\displaystyle{\frac{\sqrt{k^2\tau ^2+{k^{\prime}}^2}}6t^2+O(t^3). }
\end{eqnarray*}
\end{demo}

Let us give a geometric interpretation of $E_{|\sin\theta|}$ using 
the inverted open knots. 

Put 
$$
V_{|\sin\theta|}(K;x)=\int_{K}\frac{|\sin\theta |dy}{|x-y|^2}. 
$$

Let $I_x$, $\tilde{K}_x=I_x(K)$, and $\tilde{y}=I_x(y)$ 
be as in subsection \ref{subsectioncosineformula}. 
Let $\Pi _x$ be the normal plane to $K$ at $x$, 
$n=-f^{\prime}(s)$ be the normal vector to $\Pi _x$, 
and $\pi _x:\RR^3\to\Pi _x$ be the orthogonal projection. 
%
Since ${I_x}_{\ast}(v_x(y))= n$ 
the absolute value of the argument of cross-ratio 
$\theta_K(x,y)$ 
is equal to the angle between $n$ and 
the tangent vector to $\tilde{K}_x$ at $\tilde{y}$. 
Hence 
$$
V_{|\sin\theta|}(K;x)=
\int_{K}\frac{|\sin\theta |dy}{|x-y|^2}
=\int_{\tilde{K}_x}|\sin\theta ||d\tilde{y}|
$$
is the length of the projection $\pi _x(\tilde{K}_x)\subset \Pi _x$. 

\clearpage
\begin{figure}[htb]
\centerline{\input{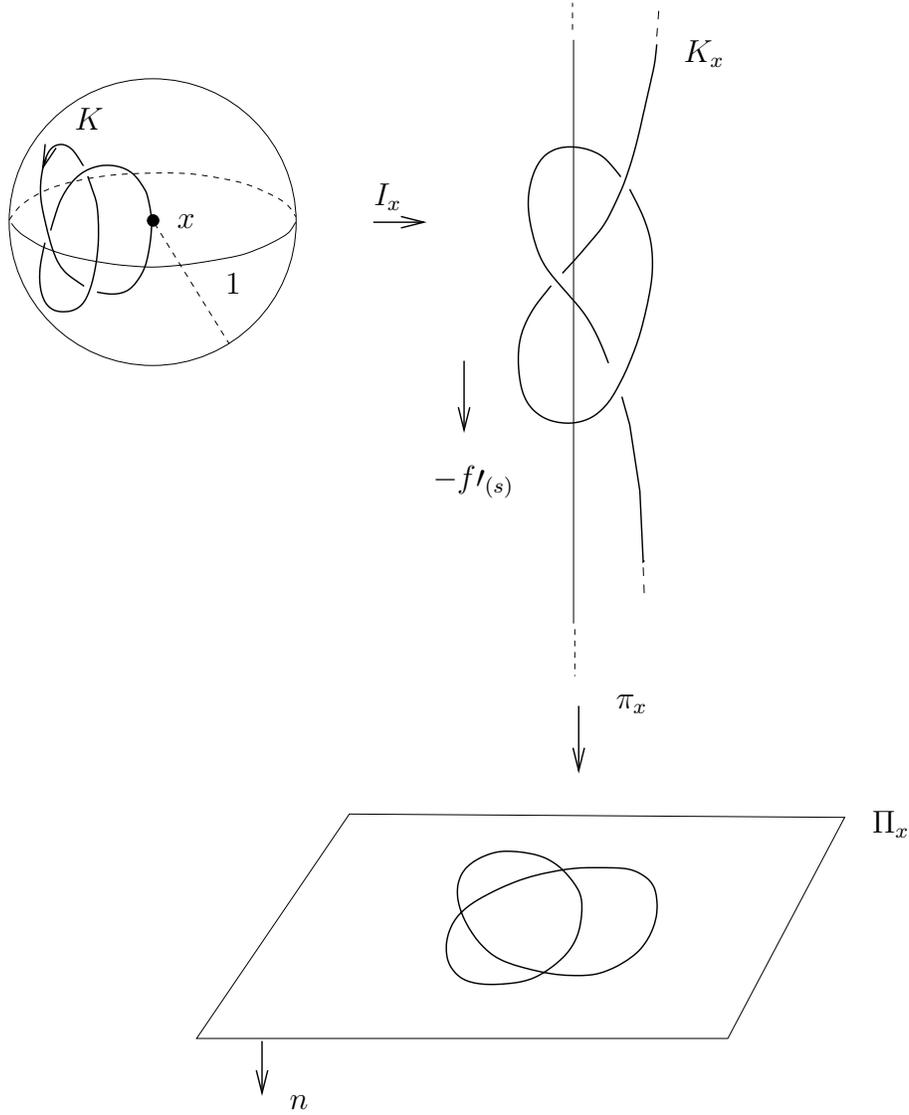}}
\caption{The projection of the inverted open knot.}    
\end{figure}

\noindent
Thus the imaginary cross-ratio energy is equal to 
the {\it total length of the projection of the inverted open knots}: 
$$
E_{|\sin\theta|}(K)=\int_K{\rm Length}(\pi_x(I_x(K)))dx.
$$
\subsection{The properties of $E_{|\sin\theta|}$.}
We show first that $E_{|\sin\theta|}$ can detect the unknot. 

\begin{exam}\rm 
Let $K_{\circ}$ be the standard planar circle. 
Then $\theta_{K_{\circ}}(x,y)=0$ for any $(x,y)$. 
Therefore $E_{|\sin\theta|}(K_{\circ})=0$. 
\end{exam}

We remark that Claim \ref{arg0} implies that $E_{|\sin\theta|}(K)=0$ 
if and only if $K=K_{\circ}$. 

\begin{theo} {\rm (\cite{Lin})}

The imaginary cross-ratio energy  crossing number 
from above, namely, 
$$E_{|\sin\theta|}(K)\ge 4 \pi ac(K)$$
for any knot $K$. 
\end{theo}
\begin{coro}
If $K$ is a non-trivial knot then $E_{|\sin\theta|}(K)\ge 12 \pi$. 
\end{coro}
\begin{demo} 
Assume that a knot is parametrized by the arc-length. 
Put 
$$
u=f^{\prime}(s), v=f^{\prime}(t), w=\frac{f(t)-f(s)}{|f(t)-f(s)|}, 
\quad \mbox{and} \quad  \tilde u=v_{f(s)}(f(t)). 
$$
Recall that $\tilde u=2(u,w)w-u$, and that the conformal angle $\theta$ 
is the angle between $\tilde u$ and $v$. 

Suppose $\theta \not= 0, \pi$. Then the tangent plane $T_{f(t)} \Sigma$ at $f(t)$ to the twice tangent spere $\Sigma = \Sigma_K (f(s), f(t))$ is spanned by $\tilde u$ and $v$. Let $\psi=\psi(s,t)$ be the angle between the chord joining $f(t)$ and $f(s)$ and the tangent plane $T_{f(t)} \Sigma$.

Then the numerator of the integrand of the average crossing number (which is the absolute value of the integrand in Gauss formula) is:
$$
|(u\times v , w)|= |(\tilde u\times v ,w)| = |\sin \psi | |\sin \theta |
$$ 

Therefore 
$$
ac(K)= \frac{1}{4 \pi} 
\iint_{S^1 \times S^1} 
\frac{|\sin \psi (s,t)||\sin \theta (f(s),f(t))|}{|f(s)-f(t)|^2}
dsdt \leq \frac{1}{4\pi} E_{|\sin \theta |}(K).
$$

When $\theta = 0$ or $\pi$ then, $|(u\times v , w)|= \sin \theta =0$. 

\end{demo}

Next we show that $E_{|\sin\theta|}$ is an energy functional for knots in the weak sense. 

\begin{theo} 
\label{sinthetacurvature}
For any real numbers $b, \delta$ 
$(0<\delta\le\frac12)$, and $\kappa_0$, there exists a positive 
constant $C=C(b, \delta, \kappa_0)$ such that if a knot $K$ with 
length $l(K)$ whose curvature is not greater than $\kappa_0$ 
has a pair of points $x,y$ on it which satisfy that 
$\delta _K(x,y)=\delta l(K)$ and that $|x-y|\le Cl(K)$, 
then $E_{|\sin\theta|}(K)\ge b$. 
\end{theo}

\begin{demo}
Fix $\delta$ and $\kappa_0$ $(\kappa_0\ge 2\pi)$, and let $0<d\le\delta$. 
Put 
$$
{\mathcal K}(d)={\mathcal K}_{\delta, \kappa_0}(d)=
\left\{
K: {\rm a}\,  {\rm knot}
\,
\left|
\begin{array}{l} 
\mbox{The length of $K$ is $1$.} \\
\mbox{The  curvature of $K$ is not greater than $\kappa_0$. }\\
\mbox{$\exists x,y\in K$ such that} \\
\mbox{i) the shorter arc-length between $x$ and $y$ is $\delta$,}\\
\mbox{ii) $|x-y|\le d$.} 
\end{array}
\right\}
\right.
$$
We show 
$$
\lim_{d\to+0}\left(\inf_{K\in{\mathcal K}(d)}E_{|\sin\theta|}(K)\right)=\infty. 
$$
\begin{lemm}
Let 
$$
d_0=\min\left\{\frac1{100}, \frac{\delta ^2}5, 
\left(\frac{-1+\sqrt{1+\frac1{50\kappa_0}}}{2}
\right)^2 \right\}.
$$
Let $d\le d_0$. 
Suppose $K=f(S^1)\in{\mathcal K}(d)$ is parametrized by arc-length and 
satisfies that $|f(0)-f(\delta)|=d$. 
If $5d\le t\le\sqrt{d}$ then at least one of 
$$
V_{|\sin\theta|}(K;f(\pm t)), V_{|\sin\theta|}(K;f(\delta\pm t))
$$
is greater than 
$$
\frac1{100}\cdot\frac1{d+t}. 
$$
\label{1/100forsin}
\end{lemm}
We show that Lemma implies Theorem. 
Since $10d\le\sqrt{d}<\frac{\delta}2$ 
$$
f([-\sqrt{d}, -5d])\sqcup
f([5d, \sqrt{d}])\sqcup
f([\delta +5d, \delta +\sqrt{d}])\sqcup
f([\delta -\sqrt{d}, \delta  -5d])
$$
is a disjoint union of curve segments of $K$ where 
$S^1$ is regarded as $\RR$ modulo $\ZZ$. 
Since $V_{|\sin\theta|}(K;y)\ge 0$ for any $y\in K$, 
\begin{eqnarray*}
E_{|\sin\theta|}(K)  
&\ge&\displaystyle{\int_{5d}^{\sqrt{d}}\{V_{|\sin\theta|}(K;f(-t))
+V_{|\sin\theta|}(K;f(t))} \\[2mm]
&&\displaystyle{+V_{|\sin\theta|}(K;f(\delta
-t)) +V_{|\sin\theta|}(K;f(\delta +t))\}dt} \\[2mm]
&\ge&\displaystyle{\int_{5d}^{\sqrt{d}}\frac1{100}\cdot\frac1{d+t}dt
=\frac1{100}\left\{\log\left(1+\frac1{\sqrt{d}}\right)-\log 6\right\},}
\end{eqnarray*}
which explodes as $d$ goes down to 0, which completes the proof of theorem \ref{sinthetacurvature} . 
\end{demo}

\medskip
\noindent
{\bf Proof of Lemma \ref{1/100forsin}:} Assume there is a $t$ with $5d\le t\le\sqrt{d}$ 
such that 
$$
V_{|\sin\theta|}(K;f(\pm t)), V_{|\sin\theta|}(K;f(\delta\pm t)) \ge
\frac1{100}\cdot\frac1{d+t}. 
$$
Put $x_{\pm}=f(\pm t)$. 
Let $R_t=1/|f^{\prime\prime}(t)|$ be the radius of curvature at $x_+$, 
and $C(x_+, x_+, x_+)$, $\tilde K_{x_+}=I_{x_+}(K)$, 
$\Pi_{x_+}$, and $\pi_{x_+}$ be as before. 
Define $x^{\prime}_+$ by 
$$
C(x_+, x_+, x_+)\cap\Pi_{x_+}=\{x_+, x^{\prime}_+\}.
$$
Let $\hat x_+=I_{x_+}(x^{\prime}_+)$. Then 
$$
\lim_{K\setminus\{x_+\}\ni y\to x_+}\pi_{x_+}(y)=\hat x_+. 
$$

\begin{figure}[ht]
\centerline{\input{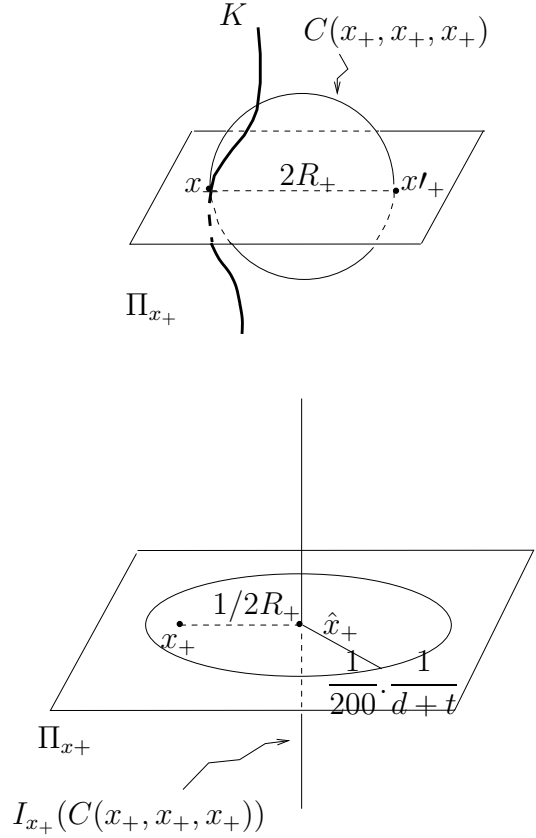}}
\caption{The osculating circle, $x_+, x^{\prime}_+$, and $\hat x_+$.}
\label{osccirc}
\end{figure}

\begin{figure}[ht]
\centerline{\input{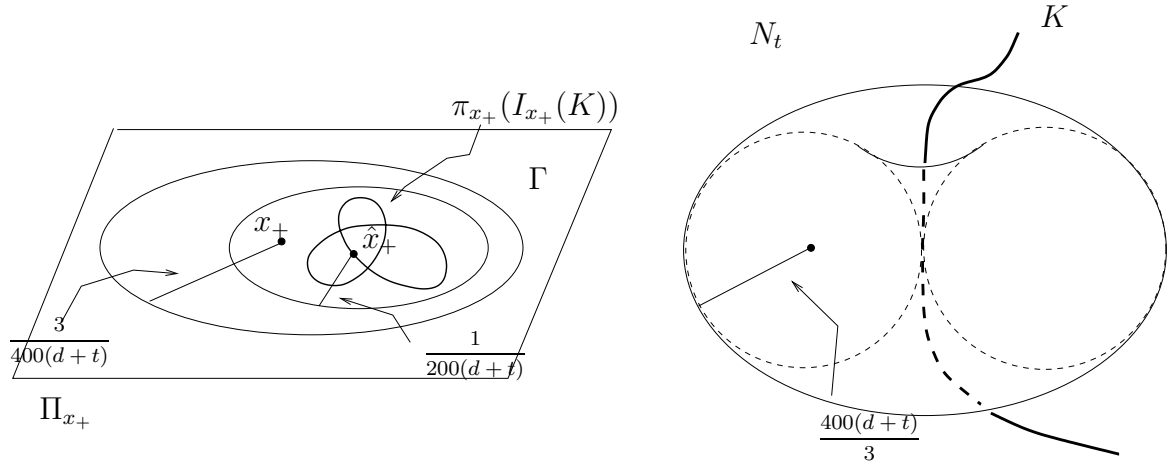}}
\caption{The circle $\Gamma$ and the degenerate solid torus.}
\label{degst}
\end{figure}

\begin{figure}[ht]
\centerline{\input{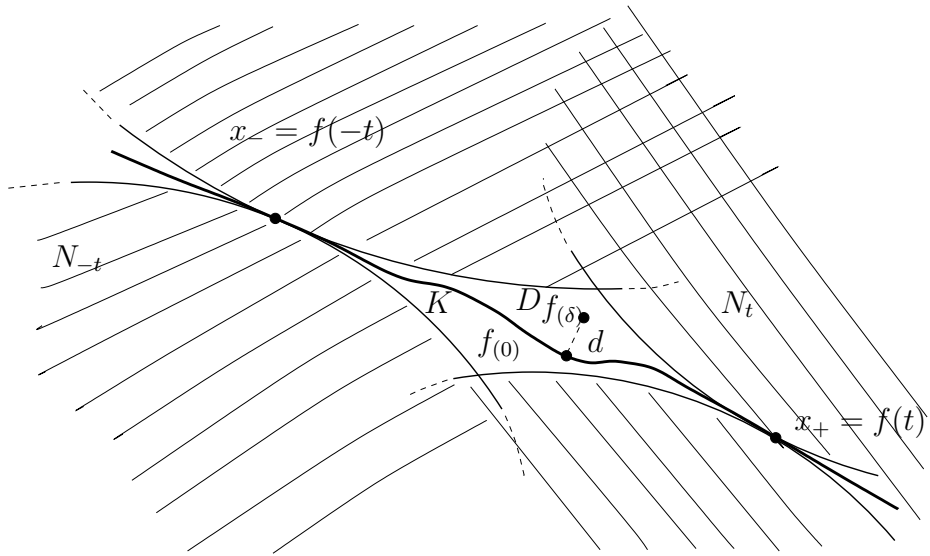}}
\caption{The domain $D$ illustrated as a 2-dimensional figure.}
\label{domaind}
\end{figure}

\noindent
Since 
$$
V_{|\sin\theta|}(K;x_+)={\rm Length} (\pi_{x_+}(\tilde K_{x_+}))
\le\frac1{100}\cdot\frac1{d+t}
$$
$\pi_{x_+}(\tilde K_{x_+})$ lies inside the circle on $\Pi_{x_+}$ 
with center $x_+$ and radius $1/(200(d+t))$. 
Since 
$$
\sqrt d\le \frac{-1+\sqrt{1+\frac1{50\kappa_0}}}{2}
$$
there holds
$$
|x_+-\hat x_+|=\frac1{2R_t}
\le \frac12\kappa_0
\le \frac12\cdot\frac1{200}\cdot\frac1{d+\sqrt{d}}
\le \frac1{400}\cdot\frac1{d+t}. 
$$
Therefore $\pi_{x_+}(\tilde K_{x_+})$ lies inside the circle $\Gamma$ 
on $\Pi_{x_+}$ with center $x_+$ and radius $3/(400(d+t))$. 
This means that $\tilde K_{x_+}$ lies inside the cylinder 
$D^2\times\RR$, where $\partial D^2=\Gamma$ and the direction 
of $\RR$ is $f^{\prime}(t)$. 
If we apply the inversion $I_{x_+}$ again, this implies that $K$ lies 
outside the ``{\it degenerate}
\footnote{By a ``{\it degenerate solid torus}'' we mean a torus of 
revolution of a circle around a tangent line. }
 {\it open solid torus}'' $N_t$ whose meridian disc has radius $(400(d+t))/3$, 
as illustrated in Figure \ref{degst}. 
\medskip
\noindent
The same argument applies to $-t$ to conclude that $K$ lies 
outside the ``degenerate open solid torus'' $N_{-t}$. 
Since 
$$
|f(t)-f(-t)|\le 2t\ll\frac{400}3(d+t)
$$
$f([-t,t])$ is contained in a region $D$ 
bounded by $N_t$ and $N_{-t}$ as illustrated in Figure \ref{domaind} 
\footnote{$D$ is the bounded component of 
$\RR^3\setminus (N_t\cup N_{-t})$. }
\medskip
\noindent

\begin{subl} 
$f(\delta)$ is contained in $D$. 
\label{sublemma}
\end{subl}

\vskip 0.5cm
We show that Sublemma implies Lemma \ref{1/100forsin}. 
If $f(\delta)$ is contained in $D$, then the knot $K$ must have 
double points at $x_+$ and $x_-$, which is a contradiction. 
(End of the proof of Lemma \ref{1/100forsin}.) 
\hfill$\square$

\medskip
\noindent
{\bf Proof of Sublemma \ref{sublemma}:} 
Since $|f(0)-f(\delta )|=d$ which is much smaller than 
the radius $(400(d+t))/3$ of the meridian discs 
of the degenerate solid tori $N_t$ and $N_{-t}$, 
it suffices to show that 
$$
|f(0)-f(\pm t)|>d
$$
to prove Sublemma. 
First we remark that
$$
t\le \sqrt d\le \frac{-1+\sqrt{1+\frac1{50\kappa_0}}}{2} <\frac1{200\kappa_0}. 
$$
For $0<s\le t$ 
$$
\frac{d}{ds}(f^{\prime}(0), f^{\prime}(s))=(f^{\prime}(0), f^{\prime\prime}(s))
\le\kappa_0, 
$$
hence 
$$
(f^{\prime}(0), f^{\prime}(s))\le 1-\kappa_0 s. 
$$
As $s\le\frac 1{\kappa_0}$ we have 
$$
|f(0)-g(t)\ge\int_0^t(1-\kappa_0 s)ds
=t(1-\frac{\kappa_0}2t)
>d. 
$$
\hfill$\square$
\section{Measure of non-trivial spheres.}

\subsection{Spheres of dimension 0}
We will start with spheres of dimension 0 in  $S^1$, and study their
 positions with respect to a ``torus" $T$ made of $4$ distinct points.

Notice that, if the conformal image of four points is
$(-R,-1,1,R)$, the cross ratio is $R$. We will use this remark later, to
 give an interpretation of the modulus of a ring or of the zone between two
 non-intersecting spheres.

An oriented sphere  $\sigma$ disjoint from  $T$ bounds an interval I. We
 will say that $\sigma$ is non-trivial if I contains two points of $T$.
Informally we may say that the small enough spheres will all be trivial.

\begin{figure}[ht]
\begin{center}
\centerline{\input{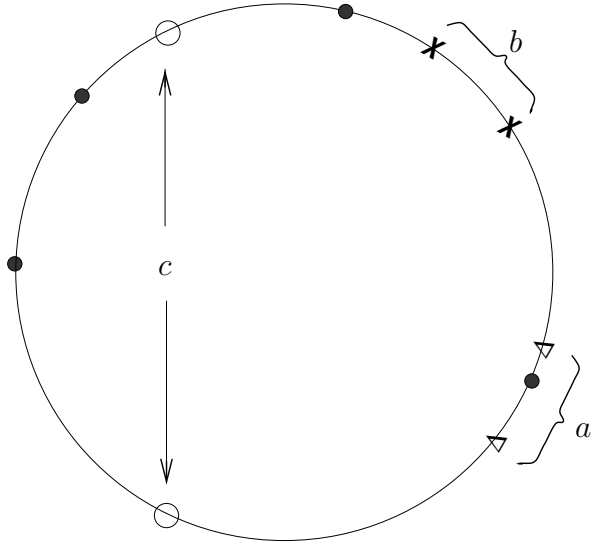}}
\caption{A non-trivial (c), and two trivial (a and b) 0-spheres.}
\end{center}
\end{figure}

\begin{prop}

 The ``torus" $T$ which minimizes the measure of the set of the non-trivial
 spheres is
the ``torus" made of the four vertices of a square (or its image by the
 conformal group of the circle).
\end{prop}

The domain $Z$ of $\cal{S}$ formed by the non-trivial spheres is bounded by 
 segments of light rays formed by the spheres containing one of the four
 points of $T$.

\begin{figure}[ht]
\centerline{\input{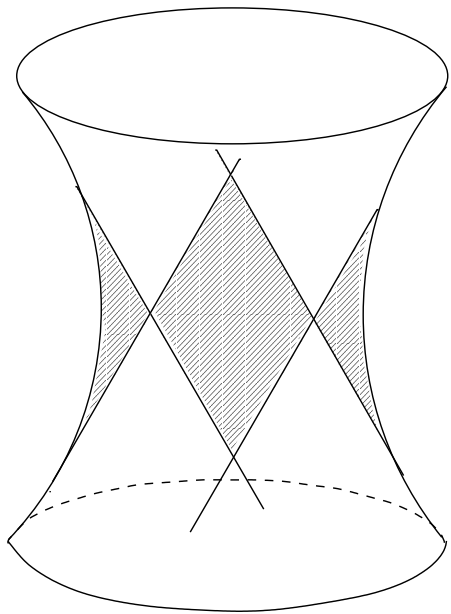}}
\caption{The set of the non-trivial 0-spheres}
\end{figure}

As the only conformal invariant of a set of four points is their 
cross-ratio, The  measure $m(Z)$ is a function of this cross-ratio.
 
\begin{demo}
The proof of the proposition is a computation. 
Using the stereographic projection of $S^1$ on $\RR$ the measure on
 ${\mathcal S} =\{ pairs \, \, of \, \, points\, \,  of\, \,  S^1 \}$ is
 $\frac{2}{(y-x)^2}|dx\wedge dy|$. Without loss of generality, we can
 suppose that the four points of the ``torus" $T$ are $\{\infty \, , \, 0 \,
 , \, 1\, , \, z\}$. We will make the computation of the measure of ``half"
 of the points of $Z$, that is $\{\infty<x<0; 1<y<z\}$, supposing $z>1$. The
 other cases are analogous.
One has $m(Z_z )=m(\{ \infty<x<0; 1<y<z \})+m(\{ 0<x<1; z<y<\infty \})$.
 One has~:
$$m(\{ \infty<x<0; 1<y<z \} )=\int_{1<y<z} \int_{-\infty <x<0}
 \frac{2}{(y-x)^2}|dx\wedge dy|=2\log(z).$$
In the same way we compute:
$$m(\{ 0<x<1; z<y<\infty \})=2\log(z)-2\log(z-1)$$
The minimum of $m(Z_z )$ is achieved for $z=2$, $m(Z_2 )= 4 \log(2)$.
This correspond to the ``square" ``torus" $T=\{ e^{ik\pi /2} \} $.

\end{demo}

\subsection{Non-trivial spheres in $\mathcal S$, tangent spheres and twice
 tangent spheres.}

Observe first that a circle is intersected by a sphere in at most two points 
if the sphere does not contain the circle. 
Generically, a sphere can intersect a knot $K$ in 0, 2 or any even 
number of points. 

We remark that the volume of the set of spheres that intersect 
a given knot $K$ is infinite for any knot $K$. 
In order to have finite valued functional, we should get rid of the 
contribution of the set of ``small spheres'', which intersect $K$ 
in at most 2 points. 

\begin{defi} \rm 
An oriented 2-sphere $\Sigma$ is called a {\it non-trivial sphere} 
for a knot $K$ if $\Sigma$ intersects $K$ in at least 4 points, 
where the number of the intersection points in $K\cap\Sigma$ 
is counted with multiplicity. 
An oriented 2-sphere $\Sigma$ with $\Sigma =\partial D^3$ is called 
a {\it non-trivial sphere in the strict sense} 
for a knot $K$ if ${\rm Int} D^3\cap K$ 
contains at least 2 connected components of $K\setminus \Sigma$. 

Let $\overline{{NT}}(K) \subset \Lambda$ (or respectively, ${NT}(K)\subset \Lambda$) 
denote the set of the non-trivial spheres 
(or respectively, the non-trivial spheres in the strict sense). 
\end{defi}

If we fix a metric on the ambient sphere $S^3$, then for a
 constant $\alpha$ depending on the curve $K$, spheres of 
geodesical radius smaller
 than $\alpha$ can intersect $K$ in at most two points. 
This implies that the set $\overline{{NT}}(K)$ of non-trivial spheres 
is bounded in ${\mathcal S}$. It has therefore a finite volume. 

\begin{defi} \rm Let ${\it mnts}_1(K)$ denote the volume of the set 
$\overline{{NT}}(K)$ of non-trivial spheres for a knot $K$. 
\end{defi}

The subset  $\overline{{NT}}(K)\subset {\mathcal S}$ is, when $K$ is not a circle,  a 4-dimensional 
region whose boundary consists of 
spheres which are tangent to the knot $K$. 
This boundary is not smooth. 
In particular, it has ``corners" which are the transverse intersections of two 
folds of the set of the tangent spheres. 
Therefore the ``corner" of $\partial \overline{{NT}}(K)$ 
consists of spheres which are tangent to the knot 
in two distinct points, that is, the twice tangent spheres. 
The closure of the boundary $\partial \overline{{NT}}(K)$ and the closure of its corner 
may have other singularities. 
\begin{defi} \rm 
Let $TWT(K)$ denote the set of the twice tangent spheres 
and $\mbox{{\it atwt}}(K)$ its area. 
\end{defi}
\begin{rema}
(1) Then the set $TWT(K)$ of the twice tangent spheres 
forms a surface of space type. 
Namely, the restriction of the Lorentz form to the tangent space 
at each regular point of this surface is positive definite. 
Therefore its area $\mbox{{\it atwt}}(K)$ is a positive number, 
which is a conformal invariant of the knot $K$. 

\medskip 
(2) Unlike $E^{(2)}$ or $E_{|\sin\theta|}$, $\mbox{{\it atwt}}(K)$ 
can not be expressed as an integral over $T^2$ of a function of 
the infinitesimal cross-ratio $\Omega_{CR}$ which does not 
contain its derivatives. 
This is because the local contribution of a neighborhood of a pair 
of points on $K$ to $\mbox{{\it atwt}}(K)$ depends on the curvature, 
whereas $\Omega_{CR}$ itself is independent of $f''$. 

\medskip
\begin{ques} Can $\mbox{{\it atwt}}(K)$ be expressed in terms of 
 of the infinitesimal cross-ratio and its derivatives:
$$
\frac{e^{i\theta}}{|f(s)-f(t)|^2}, \, 
\frac{\partial}{\partial s}\left(\frac{e^{i\theta}}{|f(s)-f(t)|^2}\right), 
etc. ?
$$
\end{ques}

\medskip 
(3) The area $\mbox{{\it atwt}}(K)$ is not an energy functional for knots. 

it is enough to check that the contribution to $\mbox{{\it atwt}}(K_n)$ of two ortogonal segments of length $1$, $I_n \subset K_n$ of middle point $p_n$, and $J_n\subset K_n$ of middle point $q_n$, such that $d(I_n ,J_n )\, =\, d(p_n ,q_n )\, =\,  1/n$ does not blow up when $n$ goes to $\infty$. This contribution is bounded by the area of the spheres tangent to two orthogonal lines distant of $1$. 
\end{rema}


\subsection{The volume of the set of the non-trivial spheres. }

The measure ${\it mnts}_1(K)$ of the set of the non-trivial spheres 
is not expressed as a ``bilocal integral'' 
like $E^{(2)}$ or $E_{|\sin\theta|}$ 
since it does not takes into account the number of intersection points of 
$K$ with spheres. 
In order to get such an integral, we have to count each sphere $\Sigma$ 
with multiplicity, which is a number of pairs of points in $\Sigma \cap K$ 
where the algebraic intersection number has the same sign. 

\begin{defi} \rm 
The {\bf measure of acyclicity} 
${\it mnts}(K)$ 
is defined by 
$$
{\it mnts}(K) = \int_{NT(K)} C^n _2 \,
d\Sigma , 
\>\> {\rm where} \>\>
n=n_K(\Sigma )= \frac{\sharp(\Sigma \cap K)}{2},
$$
where the number of the intersection points $\sharp(\Sigma \cap K)$ 
is counted with multiplicity, and $d\Sigma$ denotes the 
${\mathcal G}$-invariant measure of $\Lambda$. 
\end{defi}

Another related functional of interest is 
the {\it measure of the image of the 4-tuple map for a knot 
with multiplicity}, which is given by 
$$
{\it mnts}^{(4)}(K)=
\int_{{\rm Conf}_4(S^1)\setminus {\mathcal C}c(K)}
{\psi _K}^{\ast}\omega _{\Lambda} 
=\int_{NT(K)} C^{2n} _4 \,
d\Sigma , 
$$
where $\omega _{\Lambda}$ is ${\mathcal G}_+$-invariant volume 4-form 
of $\Lambda$. 
This definition of ${\it mnts}^{(4)}$ 
is based on the suggestion by J. Cantanella to the second author 
when he gave a talk on ${\it mnts}(K)$. 
It can not be expressed as a ``bilocal integral''. 

\begin{clai} 
These three functionals, ${\it mnts}_1$, ${\it mnts}$, and ${\it mnts}^{(4)}$ 
are conformally invariant. 
\end{clai}

\begin{demo} 
Let $g\in{\mathcal G}$ be a conformal transformation. 
Then $g$ maps a non-trivial sphere $\Sigma$ of $K$ 
to a non-trivial sphere $g\cdot\Sigma$ of $g\cdot K$ 
preserving the number of the intersection points. 
Therefore $g(NT(K))=NT(g\cdot K)$ and 
$n_K(\Sigma )=n_{g\cdot K}(g\cdot\Sigma )$. 
Since the measure $d\Sigma $ is ${\mathcal G}$-invariant, i.e. 
$d\Sigma =d(g\cdot\Sigma)$, we have 
${\it mnts}(g\cdot K)={\it mnts}(K)$ etc. 
\end{demo} 

As we mentioned before, 
the measure ${\it mnts}_1(K)$ of the set $\overline{{NT}}(K)$ 
of the non-trivial spheres 
is finite for any smooth knot $K$ since $NT(K)$ is bounded in $\Lambda$. 
On the other hand, it is not easy to control 
the number of the intersection points of $K$ and a sphere $\Sigma$. 
We will show later at the end of the next subsection that 
the measure of acyclicity is finite for a smooth knot using formula \ref{formulaofmeasure}. 

\begin{ques}
Is ${\it mnts}^{(4)}$ finite for any smooth knot? 
\end{ques}

We conjecture that for any smooth knot $K$ 
there exists a natural number $n$ such that 
the measure of the set of spheres that intersect $K$ 
in at least $n$ points is 0, 
and hence the answer to the above question is affirmative. 

\begin{rema} 
These three functionals, ${\it mnts}_1$, ${\it mnts}$, and ${\it mnts}^{(4)}$ 
are different. 
There holds 
${\it mnts}_1(K)\le {\it mnts}(K)\le {\it mnts}^{(4)}(K). $ 
In fact, 
if an oriented 2-sphere $\Sigma$ has exactly 6 transversal 
intersection points with $K$, 
then $\Sigma$ is counted $1$, $3=C_2^3$, 
and $15=C_4^6$ times respectively in these functionals. 
\end{rema}

\medskip
\subsection{The measure of acyclicity in terms of the infinitesimal cross-ratio. }

Assume that $K$ is oriented. 
Let $p\in K\cap\Sigma$ be a transversal intersection point. 
We say that {\it $\Sigma$ intersects $K$ at $p$ positively} 
if the orientation of $K$ is ``{\sl outward}'' at $p$ 
with respect to the orientation of $\Sigma$, namely, 
if the algebraic intersection number $K\cdot \Sigma$ 
is equal to $+1$ at $p$. 
We call $p$ a {\it positive intersection point}. 
Then the measure of acyclicity satisfies 
$$
{\it mnts}(K) = \int_{NT(K)} C^{n^{\prime}} _2 \,\omega , 
$$
where $n^{\prime}$ is the number of positive intersection points of 
$K$ and $\Sigma$. 
We remark that the above definition does not depend on the choice 
of the orientation of a knot $K$. 

\vskip 3 mm

Let $|{NT}_K^{(4)}(x,x^{\prime};y,y^{\prime})|$ denote 
the volume of the set of the non-trivial spheres for $K$ 
which intersect $xx^{\prime}$ and $yy^{\prime}$ positively, 
where $xx^{\prime}$ denotes the curve segment of $K$ 
between $x$ and $x^{\prime}$ with the natural orientation 
derived from that of $K$. 
Put $|{NT}_K^{(4)}(dx, dy)|=|{NT}_K^{(4)}(x,x+dx;y,y+dy)|$. 
Let 
$$|{nt}_K^{(4)}(x,y)|=\displaystyle{\lim_{|x-x^{\prime}|,|y-y^{\prime}|\to 0}
\frac{|{NT}_K^{(4)}(x,x^{\prime};y,y^{\prime})|}
{|x-x^{\prime}||y-y^{\prime}|}}.$$ 
Then the measure of acyclicity of $K$ is then equal to 
$$
{\it mnts}(K)=\lim_{\max |x_i-x_{i+1}|\to 0}\sum_{i < j} 
|{NT}_K^{(4)}(x_i, x_{i+1}, x_j, x_{j+1})|
= \frac{1}{2} \iint_{K\times K\setminus\triangle} |{nt}_K^{(4)}(x,y)| dxdy. 
$$

\begin{prop}
The measure of acyclicity 
is expressed as 
$$
{\it mnts}(K)
=\frac{\pi}{4}
\iint _{K\times K}(\sin\theta -\theta\cos\theta )\frac{dxdy}{|x-y|^2}, 
$$
where $\theta$ is the argument of cross-ratio. 
\label{formulaofmeasure}
\end{prop}
The proof decomposes into the following two lemmas.

\begin{lemm} There holds the dimension reduction formula: 
$$
|{NT}_K^{(4)}(dx, dy)|=
\frac{\pi}2|{NT}_K^{(3)}(dx, dy)|, 
$$
where ${NT}_K^{(3)}(dx, dy)$ denotes the  the set of the 
non-trivial circles of the sphere $\Sigma (dx,dy)$ tangent to the knot at $x$ and $y$ which intersect both $dx$ and $dy$ positively. We denote by $|{NT}_K^{(3)}(dx, dy)|$ the volume of this set.
\end{lemm}

\begin{demo}
Since $|{NT}_K^{(4)}(dx, dy)|$ is conformally invariant, 
we may assume that $x, x+dx, y$ and $y+dy$ are in 
$\RR ^2=\{(X, Y, 0)\}\subset\RR ^3$ 
by a stereographic projection from the twice tangent sphere. 
Suppose the set of the non-trivial spheres or circles is parametrized by 
\begin{eqnarray*}
{{CTR}}^{(3)}(dx, dy)&=&\left\{(X,Y,r)\in\RR ^3_+\left|
\begin{array}{l}
\exists C \in NT_K^{(3)}(dx,dy)\\
\mbox{with center $(X,Y)$ and radius $r$}
\end{array}
\right\}\right. \\[2mm]
{{CTR}}^{(4)}(dx, dy)&=&\left\{(X,Y,Z,r)\in\RR ^4_+\left|
\begin{array}{l}
\exists \Sigma \in {NT}_K^{(4)}(dx, dy) \\
\mbox{with center $(X,Y,Z)$ and radius $r$}
\end{array}
\right\}\right.\\[2mm]
&\cong&\{(X,Y,Z,\sqrt{r^2+Z^2})\,|\,(X,Y,r)\in{{CTR}}^{(3)}(dx, dy), 
Z \in\RR \}. 
\end{eqnarray*}
Then 
\begin{eqnarray*}
|{NT}_K^{(4)}(dx, dy)|
\, &=&\displaystyle{\int_{{{CTR}}^{(4)}(dx, dy)}\frac1{r^4}dXdYdZdr }\\[2mm]
&=&\displaystyle{\int_{{{CTR}}^{(3)}(dx, dy)}
\left(\int_{-\infty}^{\infty}\frac1{(r^2+Z^2)^2}dZ\right)dXdYdr}\\[2mm]
&=&\displaystyle{\frac{\pi}2\int_{{{CTR}}^{(3)}(dx, dy)}\frac1{r^3}dXdYdr
=\frac{\pi}2|{NT}_K^{(3)}(dx, dy)|. }
\end{eqnarray*}
\end{demo}

In what follows we calculate $|{NT}_K^{(3)}(dx, dy)|$ 
in the case when $\theta$ satisfies $0\le\theta\le\frac{\pi}2$. 
Similar calculation works for $\frac{\pi}2\le\theta\le\pi$. 

We may assume that 
$$
x=\infty, \, x+dx=(0,0), \, y=(1,0), \> {\rm and} \>\> 
y+dy= (1+a,\  b) \quad (b\ge 0) 
$$

by a suitable orientation preserving conformal transformation. 
Then, as the cross-ratio is preserved, 
$$
a=\cos\theta\frac{|dx|\, |dy|}{|x-y|^2}, \, 
b=\sin\theta\frac{|dx|\, |dy|}{|x-y|^2}. 
$$

Then ${{CTR}}^{(3)}(dx, dy)$ is given by 
$$
{{CTR}}^{(3)}(dx, dy)=\left\{(X,Y,r)\in \RR ^3_+\left| 
\begin{array}{l}
X^2+Y^2<r^2\\
(X-1)^2+Y^2>r^2\\
(X-(1+a))^2+(Y-b)^2<r^2
\end{array}
\right.
\right\}. 
$$
Let ${\mathcal N}(r_0)$ be the intersection of ${{CTR}}^{(3)}(dx, dy)$ and 
the level plane defined by $\{r=r_0\}$. 
Then 
$$
|{NT}_K^{(3)}(dx, dy)|
=\int_0^{\infty}\frac1{r^3}\, {\rm area}({\mathcal N}(r))dr. 
$$
Let $C_i(r)$ be circles defined by 
$$
\begin{array}{l}
C_1(r)\, =\, \{(X, Y) \,|\, X^2+Y^2=r^2\},\\
C_2(r)\, =\, \{(X, Y) \,|\, (X-1)^2+Y^2=r^2\},\\
C_3(r)\, =\, \{(X, Y) \,|\, (X-(1+a))^2+(Y-b)^2=r^2\}.
\end{array} 
$$

\begin{figure}[ht]
\centerline{\input{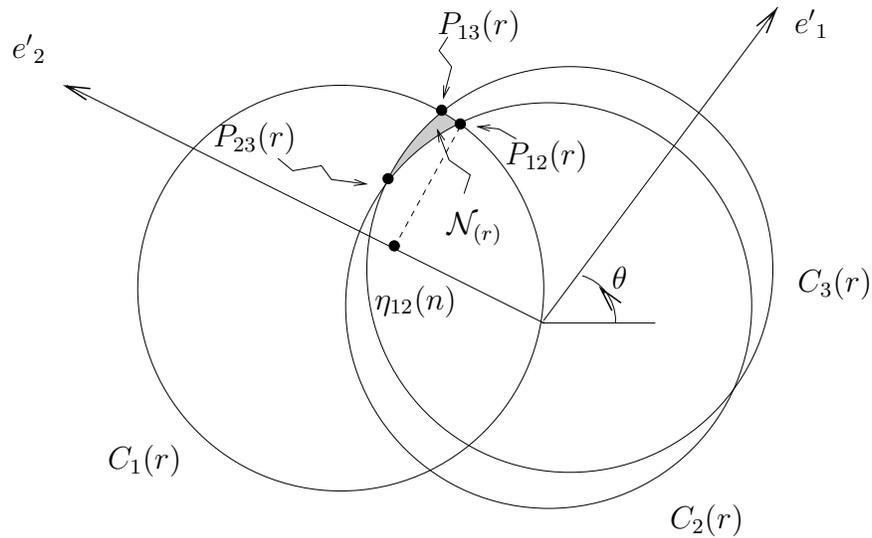}}
\caption{${\mathcal N}(r)$ and $P_{ij}(r)$} 
\label{regionMetP}   
\end{figure}

Let $P_{ij}(r)=(X_{ij}(r), \, Y_{ij}(r))$ $(i\ne j)$ be 
one of the two intersection points of $C_i(r)\cap C_j(r)$ 
(we chose the one with bigger $y$-coordinate as is shown in Figure \ref{regionMetP}).

\begin{eqnarray*}
P_{12}(r)& =&\displaystyle{ \left(\frac12, \sqrt{r^2-\frac14}\right), }\\[2mm]
P_{13}(r)&=&\displaystyle{ \left(\frac{1+a}2-b\sqrt{\frac{r^2}{(1+a)^2+b^2}-\frac14}, \, 
\frac b2+(1+a)\sqrt{\frac{r^2}{(1+a)^2+b^2}-\frac14}
\right), }\\[2mm]
P_{23}(r)& =&\displaystyle{ \left(1+\frac{a}2-\frac b{\sqrt{a^2+b^2}}
\sqrt{r^2-\frac{a^2+b^2}{4}}, \, 
\frac b2+\frac a{\sqrt{a^2+b^2}}
\sqrt{r^2-\frac{a^2+b^2}{4}}
\right). }
\end{eqnarray*}
Let $(\xi _{ij}(r), \eta _{ij}(r))$ be the coordinate of 
$P_{ij}(r)$ with respect to the frame 
$((1,0), e^{\prime}_1, e^{\prime}_2)$ 
where $e^{\prime}_1=(\cos\theta, \sin\theta)$ and 
$e^{\prime}_2=(-\sin\theta, \cos\theta)$. 
\begin{lemm} There holds 
\begin{eqnarray*}
|{NT}_K^{(3)}(dx, dy)|
& =& \displaystyle{\int
_{\frac1{2\sin\theta}}
^{\infty}
\frac 1{r^3}\cdot\sqrt{a^2+b^2}\, (\eta_{23}(r)-\eta_{12}(r))dr
+o(\sqrt{a^2+b^2}) }\\[2mm]
&=&\displaystyle{\left(\int_{\frac1{2\sin\theta}}^{\infty}
\frac{\eta_{23}(r)-\eta_{12}(r)}{r^3}dr\right)
\frac{|dx|\, |dy|}{|x-y|^2}
+o\left(\frac{|dx|\, |dy|}{|x-y|^2}\right). }
\end{eqnarray*}
\label{approximationofnt3}
\end{lemm}

\begin{demo}
 
We first remark that ${\mathcal N}(r)$ is not empty if and only if 
$X_{12}(r)\ge X_{23}(r)$, which is equivalent to 
$1-r\sin\theta\le\frac12$ when $a,b \ll 1$.

The second equality is a consequence of 
$$
\sqrt{a^2+b^2}=\frac{|dx|\, |dy|}{|x-y|^2}. 
$$

Let ${\mathcal P}(r)\subset\RR^2$ be a subset which is swept 
by the arc in $C_2(r)$ between $P_{23}(r)$ and $P_{12}(r)$ 
when it moves parallel by $\sqrt{a^2+b^2}$ in the direction 
$e^{\prime}_1$. 
Let $M_1(r)$ be a curved triangle bounded by $C_1(r)$, $C_3(r)$ 
and the half-line startin from $P_{12}(r)$ in the direction 
$e^{\prime}_1$, and let $M_2(r)$ be a subset bounded by $C_3(r)$ 
and the half-line startin from $P_{23}(r)$ in the direction 
$e^{\prime}_1$ as in Figure \ref{regionM}.

\begin{figure}[ht]
\centerline{\input{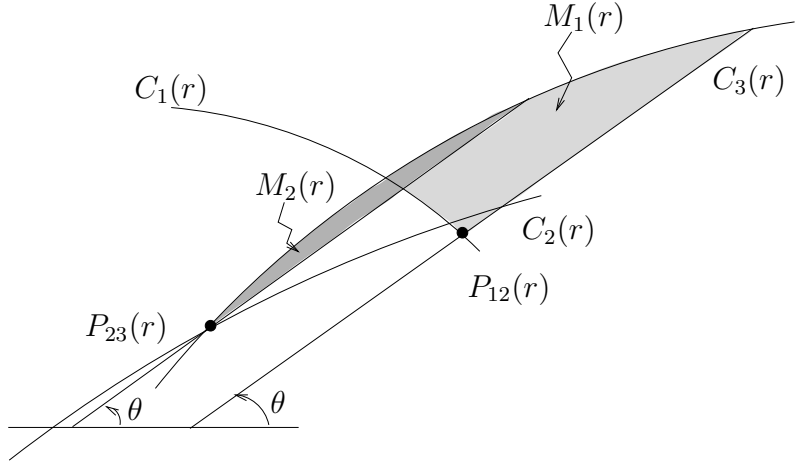}}
\caption{The regions $M_1 (r)$ and $M_2 (r)$} 
\label{regionM}   
\end{figure}
\medskip
\noindent
Then ${\mathcal N}(r)=({\mathcal P}(r)\cup M_2(r))\setminus M_1(r)$. 
We have only to show that 
\begin{equation}
\int_{\frac1{2\sin\theta}}^{\infty}
\frac 1{r^3}\cdot {\rm area}(M_i(r))\, dr
=o(\sqrt{a^2+b^2}) \quad (i=1,2). 
\label{5.4.1}
\end{equation}
Since ${\rm area}(M_1(r))\le a^2+b^2$, (\ref{5.4.1}) holds for $i=1$. 
On the other hand 

\begin{eqnarray*}
&&\displaystyle{{\rm area}(M_2(r))=\frac12\sqrt{a^2+b^2}\, 
(\eta_{13}(r)-\eta_{12}(r))
+o(\sqrt{a^2+b^2}) }\\[2mm]
&&\displaystyle{=\frac12\sqrt{a^2+b^2}\, \{\sin\theta\cdot(X_{12}(r)-X_{13}(r))
-\cos\theta\cdot(Y_{12}(r)-Y_{13}(r))\}+o(\sqrt{a^2+b^2}).  }
\end{eqnarray*}

Since 
$$
X_{12}(r)-X_{13}(r)=\sqrt{a^2+b^2}\left(-\frac{\cos\theta}2
+\sin\theta\sqrt{\frac{r^2}{(1+a)^2+b^2}-\frac14}\right), 
$$
we have 
\begin{eqnarray*}
&&\displaystyle{\int_{\frac1{2\sin\theta}}^{\infty}
\frac 1{r^3}\cdot\frac12\sqrt{a^2+b^2}\cdot \sin\theta\cdot
(X_{12}(r)-X_{13}(r))\, dr }\\[2mm]
&&\displaystyle{\le\frac{\sin\theta}2(a^2+b^2)
\int_{\frac1{2\sin\theta}}^{\infty}
\left(-\frac{\cos\theta}{2r^3}
+\frac{\sin\theta}{r^2}\right)dr
=O(a^2+b^2). }
\end{eqnarray*}
On the other hand since 
\begin{eqnarray*}
Y_{13}(r)-Y_{12}(r)&=&\displaystyle{\frac b2
+a\sqrt{\frac{r^2}{(1+a)^2+b^2}-\frac14}
+\sqrt{\frac{r^2}{(1+a)^2+b^2}-\frac14}
-\sqrt{r^2-\frac14}}\\[2mm]
&=&\displaystyle{\frac b2+a\sqrt{\frac{r^2}{(1+a)^2+b^2}-\frac14}
-\frac{\frac{2a+a^2+b^2}{(1+a)^2+b^2}r^2}
{\sqrt{r^2-\frac14}+\sqrt{\frac{r^2}{(1+a)^2+b^2}-\frac14}} }
\end{eqnarray*}
we have 
\begin{eqnarray*}
&&\displaystyle{\int_{\frac1{2\sin\theta}}^{\infty}
\frac 1{r^3}\cdot\frac12\sqrt{a^2+b^2}\cdot \cos\theta\cdot
(Y_{13}(r)-Y_{12}(r))\, dr   }\\[2mm]
&&\displaystyle{\le\frac{\cos\theta}2(a^2+b^2)
\int_{\frac1{2\sin\theta}}^{\infty}
\left(\frac{\sin\theta}{2r^3}+\frac{\cos\theta}{r^2}
+\frac{\cos\theta}{r\sqrt{r^2-\frac14}}
\right)dr+o(\sqrt{a^2+b^2}),  }
\end{eqnarray*}
as $\sqrt{r^2-\frac14}\ge r\cos\theta$ when \, $r\ge\frac{1}{2\sin\theta}$, 
the right hand side satisfies 
$$
\qquad \le\frac{\cos\theta}2(a^2+b^2)
\int_{\frac{1}{2\sin\theta}}^{\infty}
\left(\frac{\sin\theta}{2r^3}+\frac{1+\cos\theta}{r^2}
\right)dr+o(\sqrt{a^2+b^2})=O(a^2+b^2). 
$$
This completes the proof of  Lemma \ref{approximationofnt3}.  
\end{demo}

\noindent
{\bf Proof of Proposition \ref{formulaofmeasure}}: Since 
\begin{eqnarray*}
\eta_{23}(r)-\eta_{12}(r)&=&\displaystyle{\sin\theta(X_{12}(r)-X_{23}(r))
-\cos\theta(Y_{12}(r)-Y_{23}(r)) }\\[2mm]
&=&\displaystyle{\sqrt{r^2-\frac{a^2+b^2}{4}}-\frac{\sin\theta}2
-\cos\theta\sqrt{r^2-\frac14} }
\end{eqnarray*}
we have 
$$
\int_{\frac1{2\sin\theta}}^{\infty}
\frac {\eta_{23}(r)-\eta_{12}(r)}{r^3}dr
=\int_{\frac1{2\sin\theta}}^{\infty}
\frac 1{r^3}\left(r-\frac{\sin\theta}2
-\cos\theta\sqrt{r^2-\frac14}
\right)dr+o(\sqrt{a^2+b^2})
$$

\begin{eqnarray*}
&=&\displaystyle{2\int_0^{\sin\theta}
(1-u\sin\theta-\cos\theta\sqrt{1-u^2}\,)du 
+o(\sqrt{a^2+b^2}) \qquad \left(u=\frac1{2r}\right)  }\\[2mm]
&=&\sin\theta -\theta\cos\theta +o(\sqrt{a^2+b^2}), 
\end{eqnarray*}
which, combined with Lemma \ref{approximationofnt3}, 
completes the proof of Proposition \ref{formulaofmeasure}. 
\hfill$\square$

Lemma \ref{orderofarg} implies that the integrand 
$(\sin\theta -\theta\cos\theta)/|x-y|^2$ 
of the measure of acyclicity is 0 at the diagonal set. 
Therefore, the measure of acyclicity is finite for any knot of class $C^4$. 
Lemma \ref{orderofarg} also implies that 
$E^{(2)}$, $E_{|\sin\theta|}$, and ${\it mnts}$ are independent. 

We remark that Claim \ref{arg0} implies that the measure of acyclicity of $K$, 
and hence ${\it mnts}^{(4)}(K)$, is equal to 0 
if and only if $K$ is the standard circle. 


\section{Non-trivial zones. }

\subsection{Non-trivial zones.}

A circle is generically intersected by a sphere in $0$ or two point. It is a
 necessary and sufficient condition for a knot to be the circle. 
A curve which is not a circle should therefore admit spheres intersecting it
 transversally in at least $4$ points. Thickening such a non-trivial
 sphere, we get a region bounded by two disjoint spheres which is crossed by
 at least $4$ strands.

\begin{defi}\rm 
(1) A {\it zone} $Z$ is a region of $S^3$ or $\RR^3$ diffeomorphic to $S^2\times [0,1]$ 
bounded by two disjoint spheres. 

\medskip 
(2) Let $Z$ be a zone with $\partial Z=S_1\cup S_2$ and 
$T:\RR^3$ (or $S^3$) $\to\RR^3$ (or $\RR^3\cup\{\infty\}$ respectively) 
be a conformal transformation which maps the two boundary spheres of $Z$ 
into a concentric position. 
Then the {\it modulus} $\rho (Z)>0$ of the zone $Z$ 
is 
$\rho (Z) = |\log \frac{R_1}{R_2}|$
where $R_1$ and $R_2$ are the radii of the two concentric boundary spheres of $T(Z)$. 

\medskip 
(3) Instead of using a particular conformal transformation $T$, 
one can use the {\it Lorentzian modulus} $\lambda (Z)=|L(S_1, S_2)|$ 
of a zone $Z$ with $\partial Z=S_1\cup S_2$, 
where $S_1$ and $S_2$ are considered
\footnote{$S_1$ and $-S_2$ are endowed with the orientation 
as the boundary of $Z$} 
 as points in $\Lambda$. 
\end{defi} 

\begin{rema} \label{rem_module}

\noindent (1) A zone $Z$ is uniquely determined by its modulus 
up to a motion by a conformal transformation. 

\medskip
\noindent (2) The modulus $\rho (Z)$ is an increasing function of the Lorentzian 
modulus $\lambda (Z)$ with: 
\begin{equation}\label{relation}
\rho =\log (\lambda +\sqrt{\lambda ^2-1}) \ \mbox{and}  \ \lambda =Ch \rho. 
\end{equation}
In fact, one can assume, after a conformal transformation, that $S_1$ is an equator 2-sphere of $S^3$ 
and $S_2$ is a parallel  2-sphere to $S_1$. 
Then the Poncelet pencil in $\Lambda$ generated by $S_1$ and $S_2$, 
which is a hyperbola in the 2-plane in $\RR^5$ generated by $S_1$ and $S_2$, 
can be parametrized by $(Ch t, Sh t)$, where $t=0$ corresponds to $S_1$. 
The restriction of $L$ to this Poncelet pencil is a quadratic form of type 
$(1,1)$. Let $t_0$ be the parameter for $S_2$. 
Then $\rho (Z)= t_0$ whereas $\lambda (Z)=Ch t_0$. 
\end{rema}

\begin{figure}[ht]
\centerline{\input{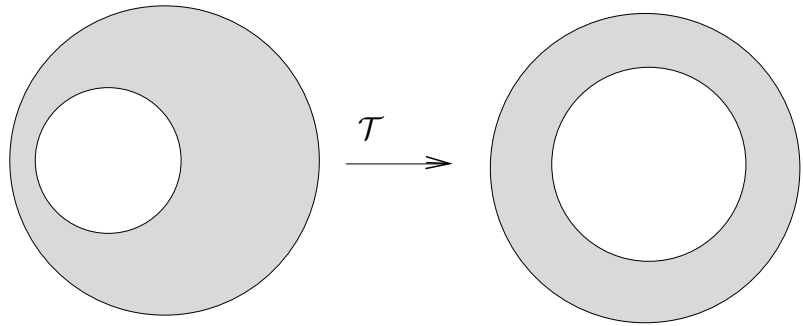}}
\caption{A zone and its privileged representative}
\end{figure}

 It can also be computed from the cross ratio of the four points
 intersections of the two spheres with a circle containing the two base
 points of the pencil generated by the two spheres.

\begin{rema}
 If the zones $Z_1$ and $Z_2$ satisfy $Z_1 \subset Z_2$, one has 
 $\rho(Z_2) \geq \rho(Z_1 ).$ 
\end{rema}

\begin{defi} \rm 
Let $K$ be a knot and $Z$ be a zone with 
$\partial Z=S_1\cup S_2$. 

(1)A {\it strand} is an arc  of $K$ that is 
contained in the zone $Z$ and whose end points are on $\partial Z$. 
 A {\it small strand (of $K$ in $Z$)} is a closure 
of a connected component of the intersection of the interior of $Z$ 
with $K$ whose two end points 
lie on the same boundary sphere of $Z$. 
A {\it degenerate small strand} is a closed connected component 
of the intersection of $K$ with $\partial Z$. 

\medskip
(2) A {\it cross-strand} is a strand of $K$ 
with one end point contained in $S_1$ and another in $S_2$. 
A {\it minimal cross-strand} is a cross-strand 
whose intersection with the interior of $Z$ is connected. 
If a cross-strand has a common end point with a  small strand, 
one can connect them to obtain a longer cross-strand. 

\medskip
(3) A zone will be called {\it non-trivial} for a knot $K$ 
if its intersection with $K$ contains at least $4$ closed cross-strands of the knot $K$ with disjoint interiors. 
\end{defi}


The existence of spheres intersecting a curve in at least four points 
transversally implies the existence of non-trivial zones. 
In particular such non-trivial zones exist for a non-trivial knot $K$. 

\begin{figure}[ht]
\centerline{\input{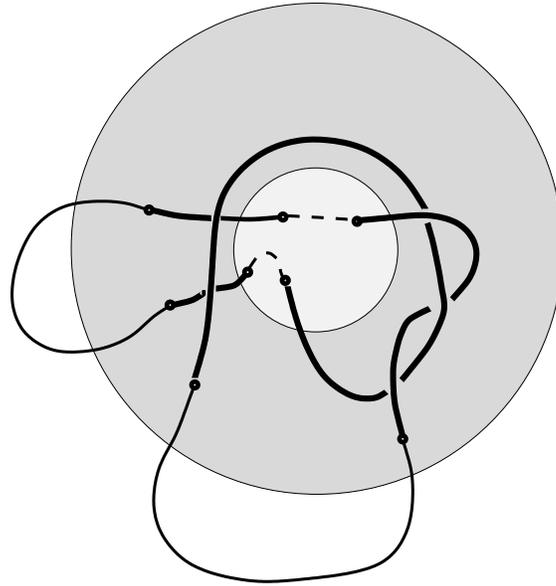}}
\caption{Non-trivial zone}
\end{figure}

\begin{defi} \rm 
The {\it modulus $\rho$ of a knot} $K$ is 
the supremum of the modulus of a non-trivial zone for $K$. 
$$\rho = \sup_{A:\mbox{{\rm  non-trivial zone}}} \rho(A)$$
\end{defi}

\begin{defi} \rm 
A non-trivial zone for a knot $K$ is called {\it maximal} if it 
has the maximum modulus among the set of non-trivial zones for $K$. 
\end{defi}

\begin{rema} 
There exists a maximal non-trivial zone for any 
non-trivial knot $K$ since the set of pairs of boundary spheres 
of non-trivial zones is a non-empty compact subset of 
$\Lambda\times \Lambda$, as a ``small" sphere intersects the knot in two points at most, and in only one point if it is tangent to $K$. 
\end{rema}

Let us now look for properties of maximal  non-trivial zones.

\begin{prop} \label{maxnontriv}
Let $Z$ be a maximal non-trivial zone for $K$ with 
$\partial Z=S_1\cup S_2$. appying a suitable conformal transformation, we can suppose that $S_1$ and $S_2$ are spheres of $\RR^3$ centered at the origin. Then it satisfies the following. 

$K$ is contained in $Z$. 
\end{prop}

\begin{prop} \label{tangencypositions}
Under the same assumptions as above, 
 
\noindent {\rm (1)} $K\cap S_1$ (or $K\cap S_2$) consists of two points
 where $K$ is tangent to $S_1$ (or $ S_2$, respectively). 

\medskip 
\noindent {\rm (2)} The two tangent points of $K\cap S_1$ (or $K\cap S_2$) are antipodal. 
\end{prop}

let us prove first the proposition \ref{maxnontriv}

\begin{demo}
The knot $K$ is a finite union of open minimal cross-strands: $m_1 ,...,m_n$, and arcs or points contained in $\partial Z$: $c_1 ,...,c_n$; $K=m_1 \cup c_1 \cup m_2 \cup c_2 ...\cup m_n \cup c_n $ in cyclic order on $K$. Each arc $c_i$ meets either $S_1$ or $S_2$ but cannot meet both.

Suppose that one of these arcs, say $c_1$, meets $S_1$ and exits $Z$. At least two other points or arcs $c_n$ and $c_2$ meet $S_2$. The arcs $m_3$ and $m_n \, n\geq 4$ exist. The difference $K\setminus (m_n \cup c_n \cup m_1 \cup c_1 \cup m_2 \cup c_2 \cup  m_3)$ is a closed arc containing (at least) a point $A\in S_1$ . Moving $S_1$ in the pencil of spheres tangent to $S_1$ at $A$, we can enlarge the zone $Z$ keeping minimal strands containing $m_1$ and $m_2$.  

\end{demo}

\begin{figure}[ht]
\centerline{\input{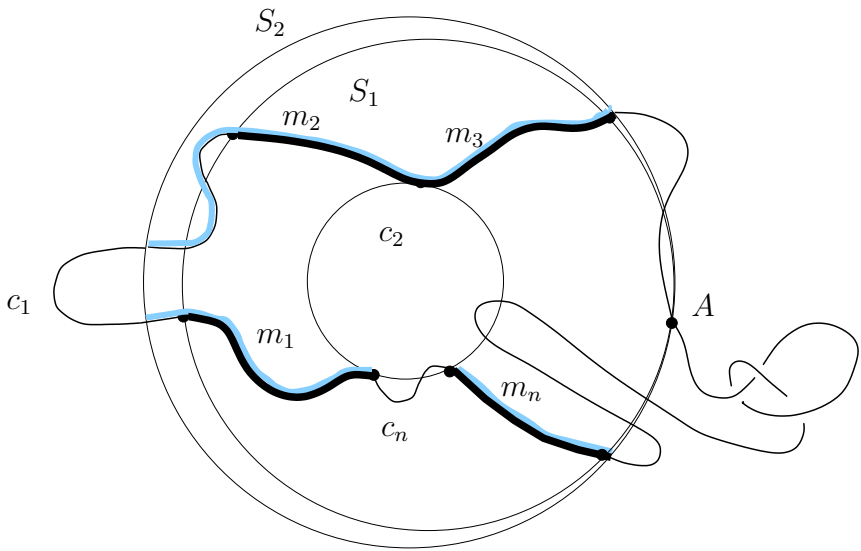}}
\caption{Increasing the modulus of the zone $Z$}
\end{figure}

\begin{coro}
If a non-trivial zone $Z$ with $\partial Z=S_1 \cup S_2$ is maximal, 
then the points $c_i$ or both endpoints of the arcs $c_i$ are alternatively  contained  in  $S_1$ and  in  $S_2$. Therefore at least two are contained in or touch each of spheres $S_1$ and $S_2$ 
\end{coro}

To prove the proposition \ref{tangencypositions} we will show that, if we can find on  the sphere $S_1$ two points $x_1 \in c_i$, $x_3 \in c_j$ belonging to two different arcs $c_i \not= c_j$ which are not antipodal, or in the sphere $S_2$ two points $x_2 \in c_k$, $x_4 \in c_l $ belongings to two different arcs $c_k \not= c_l$ which are not antipodal, then the zone $Z$ is not maximal.
The numeration of the four points can be chosen such that they are in cyclic order on $K$. 
\subsubsection{The modulus of 4 points}
For a moment, let us consider only the four points $x_1 ,x_2 ,x_3 ,x_4$ and forget the rest of the knot.

\begin{defi} \rm 
Let $x_1, x_2, x_3$, and $x_4$ be ordered four points 
in $\RR^3$ (or $S^3$). 

\medskip
(1) A zone $Z$ with $\partial Z=S_1\cup S_2$ is called to be 
{\it separating} if 
$x_1$ and $x_3$ are contained in one of the connected component of 
$\RR^3\setminus \mbox{{\rm Int}} Z$ (or $S^3\setminus \mbox{{\rm Int}} Z$)
and 
$x_2$ and $x_4$ in another component. 

\medskip
(2) The {\it maximal  modulus} $\rho (x_1, x_2, x_3, x_4)$  {\it of the ordered four points} 
$(x_1,x_2,x_3,x_4)$ is the supremum of the moduli 
$\rho (Z)$  
of separating zones $Z$ with $\partial Z=S_1\cup S_2$ (when  $S_1 =S_2 , \, \rho = 0$). 

\medskip
(3) A separating zone that attains the maximal modulus will be called 
the {\it maximal separating zone of the four point}. 

\medskip
(4) We can also consider $L(S_1, S_2)$,  where with $\partial Z=S_1\cup S_2$, and define $\lambda (x_1, x_2, x_3, x_4)$ to be the value of $L(S_1, S_2)$ for a maximal separating zone
(when  $S_1 =S_2 , \, \lambda =1$).

\end{defi}

We remark that if $Z$ is a maximal separating zone of the four points, 
then one of its boundary sphere should contain both $x_1$ and $x_3$, 
and the other $x_2$ and $x_4$. 

\begin{lemm} \label{max_sep_zone}

\noindent {\rm (1)} There are no separating zone, i.e. 
 $\rho =0; \ \mu _L=1$ if and only if $x_i$'s are concircular 
in such a way that $x_1$ and $x_3$ are not adjacent, i.e. 
when the the cross-ratio $(x_2 ,x_3 ;x_1 ,x_4) $ is a  real number between $0$ and $1$.

\medskip
\noindent {\rm (2)} Let $Z$ be a maximal separating zone for $(x_1, x_2, x_3, x_4)$ with $\partial Z=S_1\cup S_2$. Suppose we transform the picture by a conformal map so that the images of $S_1$ and $S_2$ are concentric. 
Then on the new picture, $x_1$ and $x_3$ are antipodal in one of the boundary spheres of $Z$, 
and $x_2$ and $x_4$ are antipodal in the other boundary sphere. 

\medskip
\noindent {\rm (3)} When separating zones exist, the maximal separating zone $Z$ is unique. 

\label{tangencieswith}
\end{lemm}

\begin{figure}[ht]
\centerline{\input{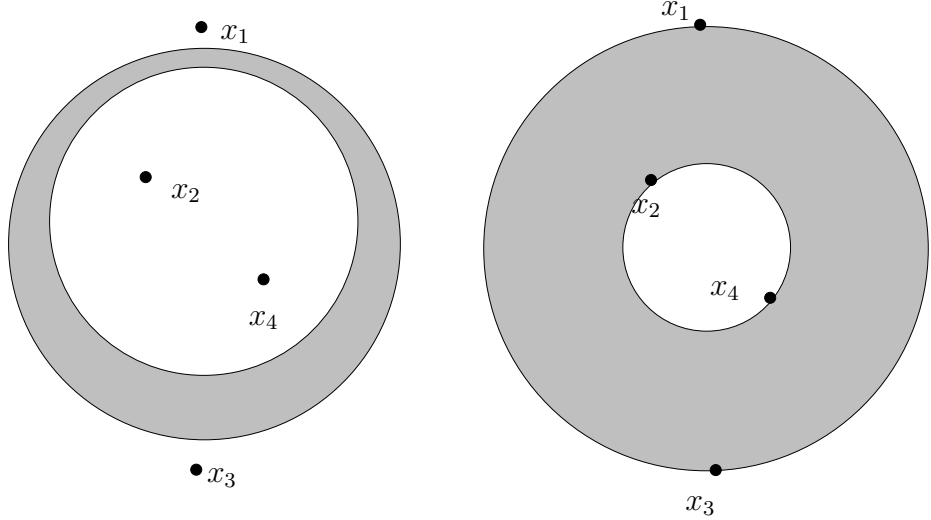}}
\caption{Maximal modulus for four points}
\end{figure}

\begin{rema}
 
\noindent (1) The condition of being antipodal does not seem to be 
conformally invariant.  

The condition can be rephrased as follows: 
the two points $x_1$ and $x_3$ have to belong to a circle 
containing the two limit points of the pencil generated by the two boundary 
spheres of the zone; the two points $x_2$ and $x_4$ satisfy a similar  condition (the  circle containing $x_2$ and $x_4$ belong to the same pencil as the circle containing $x_1$ and $x_3$). The conformal transformation such that the images of $S_1$ and $S_2$ are antipodal has to send the limit poins of the pencil generated by the two boundary 
spheres of the zone to $0$ and $\infty$; circles through the limit points become lines through the origin. 

\noindent (2) The proof of statement (3) will be given in the appendix.

\noindent (3) The statement (2) and (3)of the lemma also holds when the four points are concircular and  $(x_1 ,x_3 )$ are consecutive (and also $(x_2 ,x_4)$).

\end{rema}

The proof is the consequence of the following sublemma
\footnote{We state other proofs in the 
subsection \ref{appendix1}. }: 

\begin{subl}
Let $\mathcal P$ be a Poncelet pencil of spheres; it is the intersection of
 $\Lambda $ with a plane $P$. Let $\mathcal B$ be the two dimensional linear
 family of spheres
${\mathcal B} = B\cap \Lambda$ where $B=P^{\bot}$. 
At a point $\sigma \in \mathcal P$ the 3-plane orthogonal to
$\mathcal P$ is generated by the vectors tangent at $\sigma$ to the pencils
 of spheres with base circle a circle of the form $(\sigma^{\bot} \cap S_{\infty}) \cap (b^{\bot} \cap S_{\infty}\, (b\in
 \mathcal B)$.
\label{sublemma716}
\end{subl}

\begin{figure} \label{3plane}
\centerline{\input{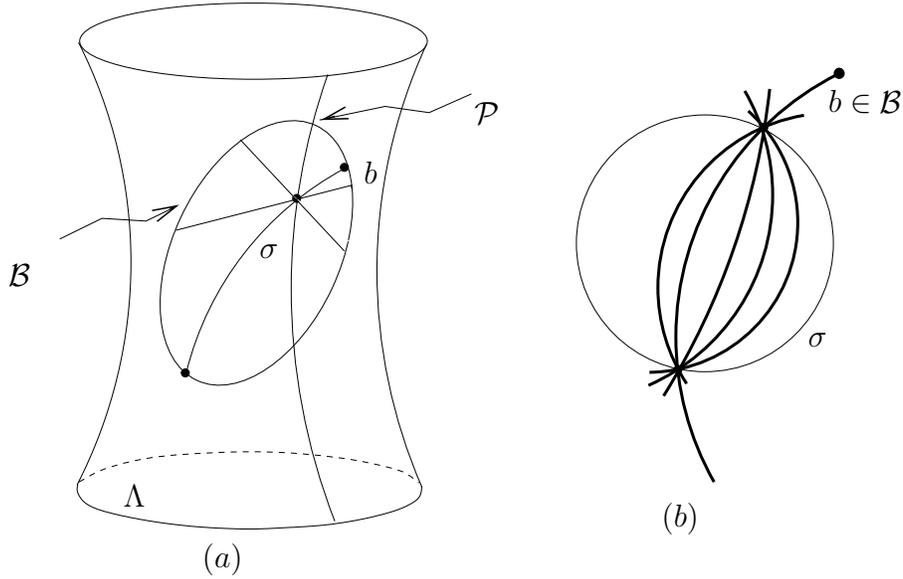}}
\caption{The 3-plane in $T_{\sigma} \Lambda$ orthogonal  at $\sigma$ to $T_{\sigma} {\mathcal P}$}
\end{figure}

{\bf Proof of the sublemma} \ref{sublemma716}{\bf : } 

The intersection ${\mathcal B} = B\cap \Lambda$ is the set of spheres
 containing the two limit points of the Poncelet pencil $\mathcal P$.  Let $F=B \oplus \RR  \cdot \sigma$. The tangent vectors to $F \cap \Lambda$ 
are in $F \cap T_{\sigma} \Lambda = F \cap \sigma^{\bot}= B$; therefore the tangent space at $\sigma$ to $F\cap \Lambda$ is an affine 3-space parallel to $B$. Hence it is orthogonal to  $\mathcal P$. 

The tangent at $\sigma$ to  $\mathcal P$ is time-like. The vector space generated by $\sigma$ and a point $b\in {\mathcal B}$ is space like. It intersects $\Lambda $ in a closed geodesic. The tangent line at $\sigma$ to this geodesic is parallel to $b$. This closed geodesic is a pencil of spheres with base circle $(\sigma^{\bot} \cap S_{\infty}^3) \cap (b^{\bot} \cap S_{\infty}^3)$ (see fig.\ref{3plane}).
\hfill$\square$

\noindent
{\bf Proof of Lemma }\ref{tangencieswith} {\bf :} 
The last sublemma implies the lemma as,  if  two points contained in  $S_1$
  are not opposite, there exists a non geodesic circle $\gamma$ containing
 them. 

Let us call ${\mathcal B}_{\gamma}$ the pencil with base circle
 $\gamma$. Consider the Poncelet pencil $\mathcal P$ generated by $S_1$ and
 $S_2$. The Lorentz modulus $\lambda $ of the region bounded by the disjoint spheres $S_2$ and
 $\sigma \in \mathcal P$ is just
$|L(S_2 ,\sigma)|$.  Let us call this function $\varphi (\sigma )$. The connected component of its level through $S_1 $ is therefore  the  hypersurface of equation $L(S_2 ,\tau )=L(S_2 , S_1)$. Its (affine) tangent plane
 at $S_1$ has the equation: $L(S_2 ,v- S_1) = 0$. The set of tangent vectors 
to this level is $B= P^{\bot}$ (where $P$ is the plane such that 
${\mathcal P} = P\cap \Lambda$).

Before finishing  the proof, we transform the picture by a conformal map sending the limit points of 
$\mathcal P$ to $0$ and $\infty$. Therefore the spheres $S_1$ and $S_2$ are concentric.

Spheres $b\in \mathcal B$ (notation of the sublemma) appear as planes through the origin. Therefore the circle basis of the pencil generated by $S_1$ and $b$ is a geodesic circle on $S_1$. 

It is equivalent to say that a vector $w \in T_{S_1} \Lambda$ is not orthogonal to $\mathcal P$ or to say that the vector $w$ is tangent at $S_1$ to a pencil with base circle which is not a geodesic circle of $S_1$.

Therefore, for $w$ a vector tangent to ${\mathcal B}_{\gamma}$ at $S_1$ we
 have $L(S_2 ,w) \not= 0$. This implies that the derivative of the Lorentzian modulus $\lambda$
 of the region between $S_2$ and a sphere $\tau \in \mathcal B$ is not zero in
 $S_1$, allowing to increase the modulus of the zone  keeping the four  points $x_1$, $x_2$, $x_3$ and $x_4$ in the boundary spheres.

What we did for the Lorentzian modulus $\lambda$ is enough to prove the same result for the modulus $\rho$, as the correspondance between $\lambda$ and $\rho$ is a diffeomorphism, see the formula (\ref{relation}).
\hfill$\square$

\subsubsection{Knots of small moduli}

The goal of this section is to prove the:

\begin{theo}
There exists a constant $a>0$ such that if a knot 
$K$ is a non-trivial knot, then its modulus is larger than $a$.
\label{modulus_knot}
\end{theo}

Let us now start with a curve $K$ of modulus $\rho$.
Let us consider a stereographic projection to $\RR^3$  such that the two
 boundary spheres of a maximal zone $Z$ are concentric, with center the
 origin.
We may imagine a very thin knitted curve contained in a zone of small
 modulus. In order to prove knot is trivial,
 we need to use the condition that the maximal modulus is very small 
to control the knot at all scales.

 \begin{figure}[ht]
\centerline{\input{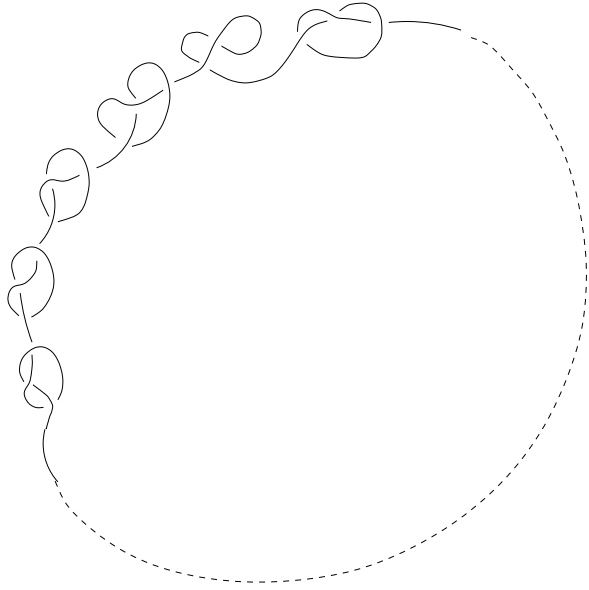}}
\caption{Knitted curve}
\end{figure}

We remark that the modulus $\rho (K)$ of a knot $K$ is 
larger than the maximal modulus $\mu (x_1 ,x_2 ,x_3 ,x_4)$ 
for any four points on the knot chosen in that cyclic order, 
since any separating zone of $(x_1 ,x_2 ,x_3 ,x_4)$ is 
a non-trivial zone for $K$.

Fixing three points $x_1 ,x_2 ,x_3$ defines a function $\mu_{x_1 ,x_2 ,x_3}$
 by:
 
\noindent $\mu_{x_1 ,x_2 ,x_3}(x) =\mu (x_1 ,x_2 ,x_3 ,x)$.
It is a continuous function, the zero level of which is the arc $\gamma$ of
 the circle trough $(x_1 ,x_2 ,x_3 )$ which does not contain the point
 $x_2$. Moreover, given a small value $\epsilon$ of $\mu_{x_1 ,x_2 ,x_3}$,
 if the sphere $S^2$ is endowed with a metric such that the length of the
 circle through $x_1 ,x_2 ,x_3 $ is of the order of $1$, then the distance
 from the level $\mu_{x_1 ,x_2 ,x_3}= \epsilon$ to the circle $\gamma$
 satisfies:
$$\forall x\in \{ \mu_{x_1 ,x_2 ,x_3} = \epsilon \}\, ; \,  a\cdot
 {\rm length}(\gamma)\, <\, d(x\, , \, \gamma ) \, 
<\, A\cdot {\rm length}(\gamma) .$$

Let us now suppose that the curve $K$ has a small modulus, say less than
 $\frac {1}{10000}$, and chose a stereographic projection of $S^3$  on
 $\RR^3$ such that the two boundary spheres $S_1 $ and $S_2$ of a maximal
 zone $Z$ for $K$ are concentric with center the origin (the north pole of
 the projection is one of the limit points of the pencil generated by $S_1$
 and $S_2$, the tangency point of $S^3$ and $\RR^3$,which we chose as origin
 of $\RR^3$, is the other limit point). Finally, compose the stereographic
 projection with an homothety to transform the inner sphere into the unit
 sphere centered at the origin. 

\begin{lemm} \label{tube1knot}
The knot $K$ is contained in a thin tubular neighborhood of a geodesic 
circle of the middle sphere $S_m $ (defined by the intersection of the ray 
containing $\frac{1}{2} ( \sigma_1 + \sigma_2)$ and  $\Lambda$, where 
$\sigma_1$ and $\sigma_2$ are the points in $\Lambda$ corresponding to the 
spheres $S_1$ and $S_2$) of the zone $Z$. 
\end{lemm}

\begin{demo}
We know that the knot $K$ is tangent to $S_1$ in two antipodal points, they 
will be $x_1$ and $x_3$.  We chose $x_2$ on the intersection of the knot 
with the ``equatorial" plane associated to the north and south poles $x_1$ and $x_3$. 
This defines a 
strand $\gamma_1$, of extremities $x_1$ and $x_3$, containing $x_2$, of 
$K$. The strand of $K$ joining $x_3$ to $x_1$ which does not contain $x_2$, 
that we will call $\gamma_2$, have to stay in a neighborhood of the arc 
of extremities $x_1$ and $x_3$, of the circle defined by the three points 
$x_1 ,x_2 ,x_3$ which does not contain $x_2$. Let $x_4$ be a point of the 
intersection of $\gamma_2$ with the equatorial plane previously defined. 
Using now the circle through $x_1 ,x_2 ,x_4$ we deduce that the strand 
$\gamma _1 $ is contained in a neighborhood of the arc of extremities 
$x_1$ and $x_2$ of the circle determined by $x_1 ,x_2 ,x_4$ which does not 
contain $x_4$.  

 \begin{figure}[ht]
\centerline{\input{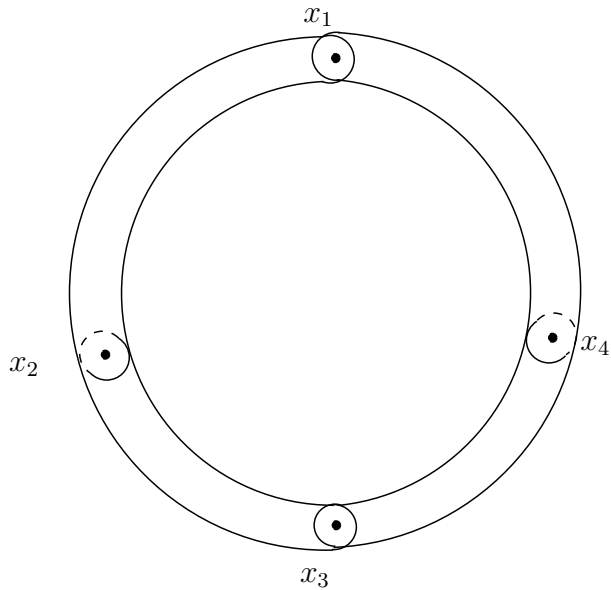}}
\caption{The knot is trapped in a neighborhood of a circle}
\end{figure}

As the two circles containing  $x_1 ,x_2 ,x_3$ and containing 
$x_1 ,x_2 ,x_4$ are close to each other, and close to a great circle $\Gamma_1$ of the
 sphere $S_m$, which completes the proof.
\end{demo}

Consider now the tube of radius $\delta _1$ around $\Gamma_1$ containing $K$, and geodesic discs $D^1_1 ,... D^1_{n_1}$ of
 this tube, normal to $\Gamma_1$  through equidistant point of $\Gamma_1$, say, such that the length
 of those arc is roughly 100 times the diameter $\delta _1$ of the tube. Chose a point $x^2_i$ of $K$ in each disc
 $D^1_i$. Joining consecutive (for the cyclic order) points $x^2_i$ and
 $x^2_{i+1}$ by a small arc of the circle defined by the two consecutive
 points $x^2_i$ and $x^2_{i+1}$ and one of the almost antipodal points to
 $x^2_i$: $x^2_{I}$, that we will call $\gamma^2_i$, we get an unknotted
 polygon $\Gamma_2$ with all the angles almost flat (the tangent of those angle
 is at most, say, $\varphi_2 \, <\, 1.5\cdot Arctg(\frac{1}{100})\, <\,
 \frac{3}{100}$).

 \begin{figure}[ht]
\centerline{\input{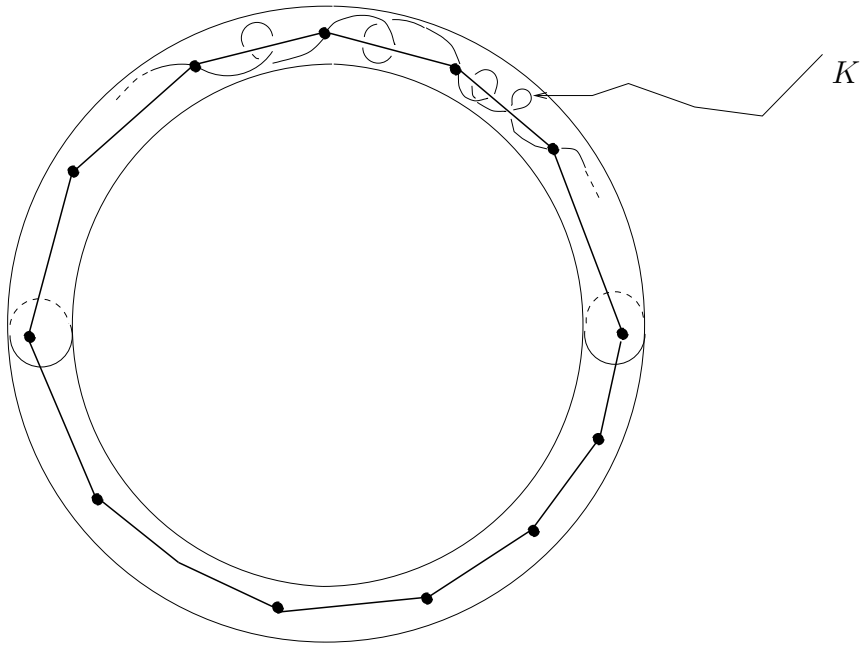}}
\caption{The first polygonal approximation}
\end{figure}

Using the the same idea as for the proof of the lemma \ref{tube1knot},
 we deduce that
 the ``small" arc of the knot joining $x^2_i$ to $x^2_{i+1}$ has to stay in
 a neighborhood of the small arc of the circle $\gamma^2_i$ of ``diameter"
 of the order of $\delta_1 \cdot (\mbox{length of} \, \gamma^2_i)$.

The knot is now confined in a  thin neighborhood
of $\Gamma_2$.  Let us call $\delta_2$ the diameter of this tube, which is also of the order of $\delta_1 \cdot (\mbox{length of} \, \gamma^2_i)$. Consider now in the tube a sequence of consecutive normal geodesic discs $D^2_1...D^2_{n_2}$ of radius  $2\cdot \delta_2$. We chose them to be distant distant of
roughly $100\cdot \delta_2$.
We can define them with no ambiguity, except in
 a neighborhood of the vertices of $\Gamma_2$.  The choice of the discs near the vertices of 
 $\Gamma_2$ has to be made with more care: they should be roughly at a
 distance of $50\delta_2$ form the vertices of $\Gamma_2$. Again we chose a
 point $x^3_i$ in each normal disc $D^2_i$. We get that way a new unknotted
 polygon $\Gamma_3$. The knot $K$ is now contained in a thin tubular neighborhood
 of $\Gamma_3$. The diameters $\delta _i$ of the successive tubes  are related to the of one order of magnitude thinner.

The only obstruction left to proceed with inductive construction of
 unknotted polygons inscribed  in $K$ is 
the angles of the polygons. We need
 to guarantee they stay flat enough. The angles $\theta^2_i$ of $\Gamma_2$
 are bounded by $\varphi_2 $ see next picture.

 \begin{figure}[ht]
\centerline{\input{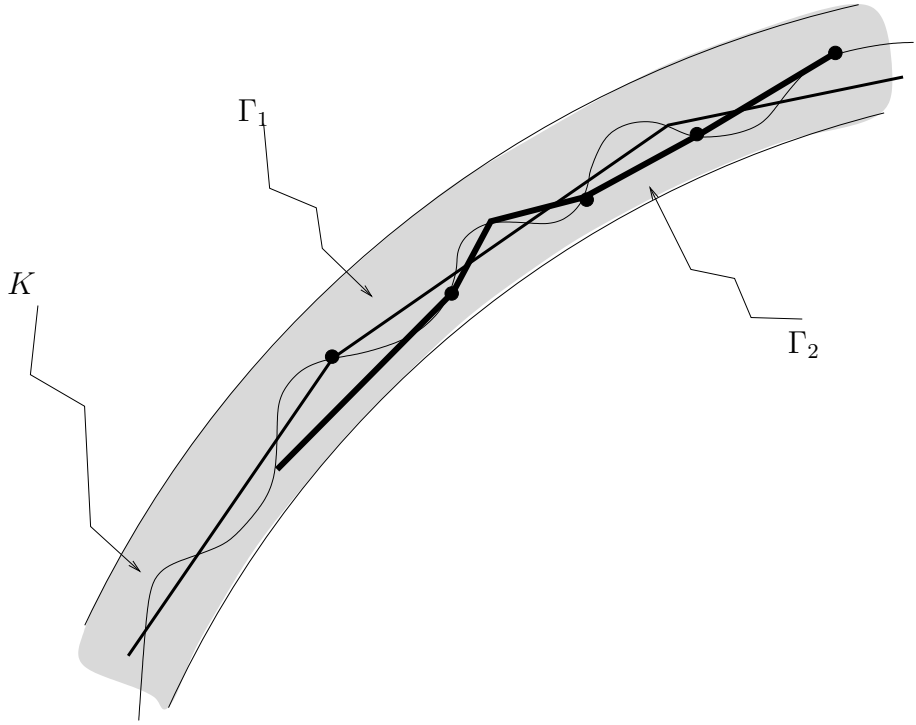}}
\caption{The second polygon $\Gamma_2$ inscribed in $K$}
\end{figure}

That is also true for the angles at all the vertices 
of the successive polygons 
except for the ones which are adjacent to some vertex 
of the previous polygon.
A priori the angle cannot grow more than linearly with the index of the
 polygone. In fact $8\cdot \varphi_2$ bounds all the angles of all the
 $\Gamma_k$. If the angle $\theta$ at a vertex of $\Gamma_{k-1}$ is larger
 than $5\cdot \varphi_2$, then our construction will replace it by two
 angles of $\Gamma_k$  smaller than $\frac{\theta}{2} +3\varphi_2$,
 therefore never reaching $8\cdot \varphi_2$


Getting to a polygon $\Gamma_k$ such that the diameter $\delta_k$ of its
 tubular neighborhood containing the knot $K$ is smaller than,
 $\frac{1}{3}\cdot \Delta (K)$ where $ \Delta (K)$ is the global 
radius of curvature defined by Gonzalez and Maddocks \cite{Go-Ma} 
(see section \ref{globalcurvature}), 
we prove that the
 knot $K$ is trivial. This ends the proof of the theorem \ref{modulus_knot}.

\begin{rema}
Freedman and He \cite{Fr-He} give a  definition of the modulus of a solid torus using the degree of maps from the solid torus $\bf T$ to the circle $\RR /\ZZ$. 
From our viewpoint (topology implies a jump of some geometrical invariant), it has the defect of beeing very small when the solid torus is very long and very thin.

But, in the same article Freedman and He defines a dual notion which is bounded below by a positive power of the average crossing number of knot type of $\bf T$. It would be very interesting to understand possible relations between this modulus they denote by $m^* (\bf T)$ and our modulus of a knot, for example via a ``thickest" tubular neighbourhood of a representative of a knot in its isotopy class.  

\end{rema}

\subsection{Jump of the measure of acyclicity for non-trivial knots.}

Given the curve $K$, the spheres contained in a non-trivial zone $A$
 intersect $K$ in at least four points. The region $R_A$ in $\Lambda$
 corresponding to the spheres contained in the annulus  $A$ is bounded by
 two light cones (half light cones to be precise) with vertex the boundary
 spheres $S_1$ and $S_2$ of $A$. The volume of $R_A$ is a monotonous
 function of the modulus of $A$. Therefore the existence of a strictly
 larger than $0$ lower bound for the moduli of maximal non-trivial zones
 associated to non-trivial knots implies the existence of a strictly
 positive constant $c>0$ such that for any non-trivial knot, the measure
 ${\it mnts}_1 (C)$ of the set of non-trivial spheres is larger or equal to $c$.

The modulus and the measure of acyclicity probably 
do not behave the same way for connected sums of knots.
Let $\rho ([K])$ (or ${\sl mnts}([K])$) denote 
the infimum of the modulus (or the measure of acyclicity respectively) 
of a knot that belongs to an isotopy class $[K]$. 

\begin{ques}
Is it true that:
$$\rho ([K_1 \sharp K_2] ) =
\max\{\rho([K_1] ) \, ,
\, \rho([K_2] )\} \, ?$$
\end{ques}

\begin{ques}
Is it true that:
$${\it mnts} ([K_1 \sharp K_2] ) = (\mbox{or}\, \leq \, \mbox{or} \, \geq )\, 
{\it mnts}([K_1] ) + {\it mnts}
 ([K_2] ) \, ?$$
\end{ques}

\begin{ques}
Do the modulus and the measure of acyclicity take 
their minimum values for non-trivial knots 
at a (2,3)-torus knot on a torus of revolution?
\end{ques}


\medskip
\subsection{Explosion of the measure of acyclicity for singular knots. } 
\begin{theo} 
The measure of acyclicity 
is an energy functional for knots. 
\end{theo}

\begin{demo}
Fix $\delta$ $(0<\delta\le\frac12)$. For $0<d\le\delta$ put 

$$
\begin{array}{l} 
{{\mathcal K}}(d) ={{\mathcal K}}_{\delta}(d)\, = 
\, \left\{
\begin{array}{l}
K: \mbox{a knot with} \\
\hskip 0.7cm \mbox{length $l(K)$}
\end{array}
\,
\left|
\begin{array}{l}
\mbox{$\exists x,y\in K$ such that}\\
\mbox{i) the shorter arc-length }\\
\hskip 0.4cm \mbox{between $x$ and $y$ is $\delta l(K)$.} \\[2mm]
\mbox{ii) $|x-y|\le dl(K)$.} 
\end{array}
\right\},
\right.\\
\rho (d)\, =\, \displaystyle{\inf _{K\in {{\mathcal K}}(d)}\{{\it mnts}(K)\}}
\end{array}
$$

We remark that $\rho$ is a monotonely decreasing function of $d$ 
because $d_1<d_2$ implies ${{\mathcal K}}(d_1)\subset {{\mathcal K}}(d_2)$, 
and that $\rho (d)=0$ for 
$d\ge (\sin\pi\delta) /\pi$ 
because ${{\mathcal K}}(d)$ contains the standard circle. 

Assume that there is a positive constant $M^{\prime}$ such that 
$$
\sup_{d>0}\rho (d)=\lim _{d\to +0}\rho (d)= M^{\prime}. 
$$
There is a positive constant $M=M(M^{\prime})$ such that if 
the measure of acyclicity of 
$K$ is smaller than $1.1 M^{\prime}$ then 
the ratio of the radii of the two spheres of any concentric 
non-trivial zone for $K$ is smaller than $M$. 

Let $K_d$ be a knot with length 1 in ${{\mathcal K}}(d)$ whose 
measure of acyclicity 
is smaller than $1.1 M^{\prime}$. 
Let $x,y$ be a pair of points on $K_d$ such that 
the shorter arc-length between them is $\delta$ and that 
$d_0=|x-y|\le d$. 
\footnote{We use the axiom of choice here. }
The knot $K$ is divided into two arcs, $A_1$ and $A_2$, by $x$ and $y$. 
Put 
$$
a=\min\left\{\sup_{z\in A_1}\left|z-\frac{x+y}{2}\right|,
\sup_{w\in A_2}\left|w-\frac{x+y}{2}\right|\right\}. 
$$
Let $\Sigma _r$ denote the 2-sphere with center $(x+y)/2$ and radius $r$. 
Since $(\Sigma _a, \Sigma_{d_0/2})$ forms a non-trivial zone, 
$a\le Md_0/2$ and hence $a\le Md/2$. 
Let $\bar K_d$ be the part of $K_d$ which is contained in 
$\Sigma _{Md/2}$. The length $l(\bar K_d)$ of $\bar K_d$ 
is greater than or equal to $\delta$ since $\Sigma _{Md/2}$ 
contains at least one of the arcs $A_1$ and $A_2$. 
Then by Poincar\'e's formula ([Sa] page 259, see also pages 111, 277) 
there holds 
$$
\int_G \#(g(\Sigma _r)\cap \bar K_d)dg
=\frac{O_3O_2O_1\cdot O_0}{O_2O_1} \cdot 
{\rm area}(\Sigma _r)\cdot l(\bar K_d), 
$$
where $O_n$ denotes the volume of the $n$ dimensional unit sphere, 
$G$ denotes the group of the orientation preserving motions of $\RR^3$, 
which is isomorphic to the semidirect product of 
$\RR^3$ and $SO(3)$, 
and $dg$ denotes the kinematic density ([Sa] page 256). 
We have $dg=dP\wedge dK_{[P]}$, where $P\in\RR^3$, 
$dP$ is the volume element of $\RR^3$, and $ dK_{[P]}$ is 
the kinematic density of the group of special rotations around $P$, 
which is isomorphic to $SO(3)$. 
Since 
$$
\int dK_{[P]}=O_2O_1
$$
the Poincar\'e's formula implies 
$$
\int _{\Sigma _r(X,Y,Z)\cap K_d\ne\phi}  \#(\Sigma _r(X,Y,Z)\cap K_d)dXdYdZ
=2\pi r^2 l(\bar K_d), 
$$
where $\Sigma _r(X,Y,Z)$ denotes the sphere with radius $r$ and center 
$(X,Y,Z)$. Suppose $d/2\le r\le Md/2$. 
Then 
the measure of acyclicity of $K_d$ satisfies 
\begin{eqnarray*}
{\it mnts}(K_d)&=&\displaystyle{\int _{NT(K_d)}{}
C^{\#(\Sigma _r(X,Y,Z)\cap K_d)}_2
\cdot\frac1{r^4}dXdYdZdr }\\[2mm]
&\ge&\displaystyle{\int_{d/2}^{Md/2}\frac1{r^4}
\left\{\int _{D^3_{Md/2+r}((x+y)/2)}
\frac{\#(\Sigma _r(X,Y,Z)\cap \bar K_d)-1}2dXdYdZ
\right\}dr, }
\end{eqnarray*}
where\, $D^3_c(R)$ denotes\, the\, 3-ball\, 
with\, center\, $c$\, and \, radius $R$. Thus 
\begin{eqnarray*}
{\it mnts}(K_d)
&\ge&\displaystyle{\frac{1}{2} \int_{d/2}^{Md/2}\frac1{r^4}\left\{
2\pi r^2l(\bar K_d)-\frac{4}{3} \pi \left(\frac{Md}2+r\right)^3
\right\}dr }\\[2mm]
&\ge&\displaystyle{
\frac{1}{2}\left\{2\pi\delta\left[-\frac1r\right]_{d/2}^{Md/2}
-\frac{4}{3} \pi (Md)^3\left[-\frac1{3r^3}\right]_{d/2}^{Md/2}
\right\},}
\end{eqnarray*}
which blows up as $d$ goes down to $0$, which contradicts the assumption that 
\newline
$\lim _{d\to +0}\rho (d)\le M<\infty$. 
\end{demo}
\medskip 
\section{Appendix. }

\subsection{The maximal modulus and cross-ratio of four points. 
\label{appendix1}}

In this sub\-section we give a for\-mula to express the maximal Lorent\-zian modulus 
$\mu _L(x^1, x^2, x^3, x^4)$  of the four points 
in terms of the cross-ratio $(x^2, x^3; x^1, x^4)$. 
It implies that the integrand of $E^{(2)}$ is equal to 
the {\it infinitesimal maximal modulus} 
$\mu _L(x, x+dx, y, y+dy)$. 

But we start with an alternative proof of Lemma \ref{tangencieswith} 
by showing it in lower dimension by one. 
The terms in the following lemma can be defined in a parallel way. 

\begin{lemm} 
Let $x_1 ,x_2 ,x_3$, and $x_4$ be ordered four points on $S^2$. 
Suppose the cross-ratio $(x_2 ,x_3 ;x_1 ,x_4)$ is not a  real number between $0$ and $1$. 

Then 

\noindent (1) The maximal modulus of a separating annulus is attained by 
a separating annulus, which will be called a maximal 
separating annulus. 
Then the pair $x_1 ,x_3$ (or $x_2 ,x_4$) is on a circle of 
the pencil whose base points are the limit points of the pencil defined by 
the boundary circles of the maximal separating annulus. 

\noindent (2) There is exactly one maximal separating annulus. The position of the boundary circles of the annulus is explicitely determined in fonction of the four points.
\label{lemma811}
\end{lemm}

\begin{demo}
\noindent (1) First observe that the two boundary circles $C_1$ and $C_2$ of 
a separating annulus with the maximal modulus should contain, 
one  $x_1$ and $x_3$, the other $x_2$ and $x_4$.
The idea of the proof is the same as the idea of the proof of 
Lemma \ref{tangencieswith}. 
If these points where not on a circle $\Gamma$ through the limit point
 generated by the two boundary circles, say for $x_1 \in C_1$ and $x_3 \in
 C_1$, rotating from $C_1$ in the pencil of base points $x_1$ and $x_3$ will
 define a curve in the quadric $\Lambda$ of circles in $S^2$ transverse at
 $C_1$ to the level $L(C_2 ,C) =L(C_2 ,C_1)$ of $L$. Therefore this curve is transverse to the
 level $m(C_2 ,C) =m(C_2, C_1)$, contradicting the fact that a maximum is a
 critical point.

\noindent (2) Let $\PP \Lambda$ be the set of non-oriented spheres $\PP \Lambda = \Lambda / \pm 1$. The pencils of circles through $x_1 ,x_3$ and through $x_2 ,x_4 $ can be seen as two circles $\gamma_{13}$ and $\gamma_{24}$ in  $\PP \Lambda$. Consider the four circles determined by three of the four points $x_1 ,x_2 ,x_3 ,x_4 $. Two of them bound on $\gamma_{13}$ an open segment: the set of circles contained in the pencil $\gamma_{13}$ separating (in $S^2$)  $x_1 ,x_3$ from $x_2 ,x_4 $.  We have a similar situation on the pencil $\gamma_{24}$. 

Let us prove that the boundary circles of a maximal annulus are the middle points (for the arc-length in  $\PP \Lambda$) of the segments we just defined. This proves the uniqueness of the maximal separating annulus.

The arc-length on a pencil  in $\PP \Lambda$ is is just the angle through a base point of the pencil. 

On the picture \ref{middle_angle} we have represented the pencil $\gamma_{24}$, the circles  $C_{124}$ through $x_1 , x_2 ,x_4 $ and $C_{324}$ through   $x_3 ,x_2 ,x_4$, and the circles $S_1$ and $S_2$, boundary of a maximal separating annulus.

The angles between $C_{124} $ and $S_2$ and the angle between $C_{324}$ and $S_2$ are equal, since the picture is symmetric with respect to the origin and since angles on two sides of a point of intersection of two circles are equal. That completes the proof. 

\begin{figure}[ht]
\centerline{\input{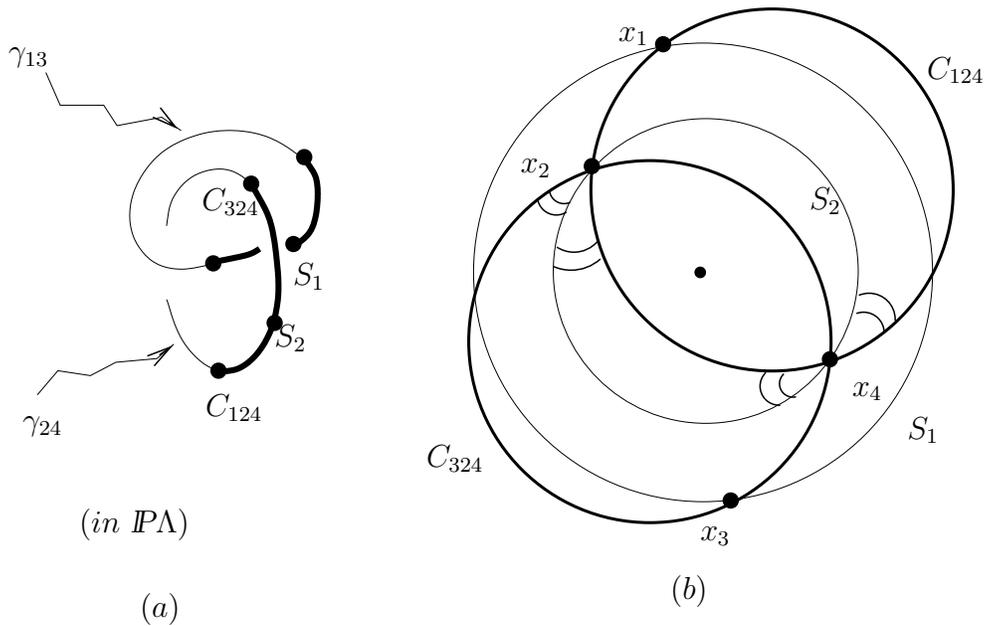}}
\label{middle_angle}
\caption{Characterization of the boundary circles of a maximal separating annulus}
\end{figure}  

\end{demo}

\begin{lemm} 
Let $C_1$ and $C_2$ be a pair of concentric circles on $\RR^2\subset\RR^3$. 
Then among the zones $Z$ with $\partial Z=S_1\cup S_2$ satisfying 
$(S_1\cup S_2)\cap\RR^2=C_1\cup C_2$, 
the maximal modulus is attained by the unique zone 
which has concentric boundary spheres whose center is in the plane $\RR^2$; these spheres are therefore orthogonal to the plane  $\RR^2$. 
\label{lemm812}
\end{lemm}

\begin{demo}
Just notice that, when we presented (in lemma \ref{max_sep_zone}) a maximal  non-trivial zone using two concentric spheres, the four pairwise antipodal points $x_1 ,x_3$ and $x_2 ,x_4$ belong the the same plane, which is othogonal to the two concentric spheres. 
\end{demo}

\begin{coro} Let $x_1 ,x_2 ,x_3$, and $x_4$ be ordered four points 
in $\RR^3$ or $S^3$, and let $\Sigma$ be a 2-sphere passing through them. 
Then the maximal modulus $\mu (x_1, x_2, x_3, x_4)$ of the four points 
is equal to the maximal modulus of a separating annulus in $\Sigma$. 
\end{coro}


\begin{clai} The maximal Lorentzian modulus $\mu _L(x^1, x^2, x^3, x^4)$ 
is expressed in terms of the cross-ratio ${\sl cr}=(x^2, x^3; x^1, x^4)$ as 
$$
\mu _L(x^1, x^2, x^3, x^4)=
\sqrt{1+2|1-{\sl cr}|^2\left\{\left|\frac{{\sl cr}}{1-{\sl cr}}\right|
-\Re\left(\frac{{\sl cr}}{1-{\sl cr}}\right)\right\}}. 
$$
\label{crandmodulus}
\end{clai}

\begin{rema}  Lemma \ref{tangencieswith} implies the above Claim as follows. 
The boundary spheres of the maximal separating zone 
can be mapped into a concentric position by a conformal transformation. 
Then Lemma \ref{tangencieswith} asserts that $x_1$ and $x_3$ are antipodal 
in one of the boundary spheres, and $x_2$ and $x_4$ are antipodal in another. 
We can assume without loss of generality that 
$x_1, x_3=\pm 1$ and $x_2, x_4=\pm z$ in $\CC\subset\RR^3$, and 
the boundary spheres of the maximal separating zone 
are concentric spheres with center the origin and radii 1 and $|z|$. 
Then the maximal Lorentzian modulus $\mu _L(x^1, x^2, x^3, x^4)$ is given by 
$(|z|+1/|z|)/2$ and the cross-ratio ${\sl cr}=(x^2, x^3; x^1, x^4)$ 
is given by $-(z-1)^2/(4z)$, and therefore 
${\sl cr}/(1-{\sl cr})=(z-1)^2/(z+1)^2$. 
One can verify the equality by the direct calculation. 
\end{rema}

\begin{lemm} Let $\Pi$ be a plane and $\Sigma$ be 
a 2-sphere in $\RR^3$ such that the minimum and maximum distance 
between $\Pi$ and $\Sigma$ are $a$ and $b$ $(0<a<b)$ respectively. 
Then the modulus $\rho (\Sigma , \Pi)$ of the zone bounded 
by $\Sigma$ and $\Pi$ is given by 
$$
\rho (\Sigma , \Pi)=\log (\frac{\sqrt b +\sqrt a}{\sqrt b -\sqrt a}). 
$$
\end{lemm}

\begin{demo} 
May assume that $\Pi $ is the $y$-$z$ plane and $\Sigma$ has 
$(\frac{a+b}2, 0, 0)$ as its center. 
Then any inversion with respect to a 2-sphere with center 
$(\pm\sqrt{ab}, 0, 0)$ maps $\Sigma$ and $\Pi$ to a 
concentric position. 
\end{demo}

This Lemma immediately implies the following Lemma. 

\begin{lemm} Let $\Sigma$ be a 2-sphere in $\RR^3$ with center $O$, 
and $A\in\RR^3$ be a point outside $\Sigma$. 
Then among planes $P$ passing through $A$, 
the plane $P_0$ perpendicular to the line $OA$ gives the maximal modulus 
$\rho (\Sigma, P)$. 
\label{maxmodperp}
\end{lemm}

Let $T_0$ be 
an inversion with respect to a 2-sphere whose center lies on 
the line $OA$ that maps $\Sigma$ and $P_0$ 
to a concentric position as is given in the proof of the previous Lemma.  
Then $T_0(A)$ and $T_0(\infty)$ are antipodal in $T_0(P_0)$. 
Therefore this lemma gives an alternative computational proof 
of Lemma \ref{tangencieswith}. 

\begin{lemm} Let $C$ be a circle with center $O$ in a plane $\Pi$, 
and $A\in\Pi$ be a point outside $C$. 
Then among pairs of a sphere $\Sigma$ passing through $C$ 
and a plane $P$ passing through $A$, 
the pair of the sphere $\Sigma _0$ whose center belongs to $\Pi$ 
and the plane $P_0$ perpendicular to the line $OA$ 
gives the maximal modulus. 
\end{lemm}

\begin{demo} Let $r$ be the radius of $C$ and $a=|O-A|$. 
Let $\Sigma (h)$ denote the 2-sphere whose center is apart from 
$\Pi$ by $h$ $(h\ge 0)$. 
Then the preceding two Lemmas imply that the maximum of 
the modulus $\rho (\Sigma (h), P)$ with $P\ni A$ is given by 
$$
\log \left( \sqrt{1+\frac{a^2-r^2}{r^2+h^2}}+\sqrt{\frac{a^2-r^2}{r^2+h^2}} \right), 
$$
which is a decreasing function of $h$. 
\end{demo}

\medskip
\noindent
{\bf Proof of Claim }\ref{crandmodulus}: Without loss of generality 
we can assume that $x^1=\infty$, $x^2=0$, $x^3=1$, and $x^4=z=u+iv$ 
in $\CC\cup\{\infty\}\cong\RR^2\cup\{\infty\}\subset\RR^3\cup\{\infty\}$. 
Then the cross-ratio ${\sl cr}=(0, 1; \infty, z)=1-\frac1z$ and hence 
$\frac{{\sl cr}}{1-{\sl cr}}=z-1.$ 
We assume that $z\not\in\{z\in\RR \, |\, z>1\}$, otherwise 
the both sides of the formula coincide with 1. 
The previous lemma implies that the maximal modulus 
$\mu _L(\infty, 0, 1, z)$ is attained by a zone $Z(\xi)$ 
bounded by $S_1(\xi)$ and $S_2(\xi)$ for some $\xi\in\RR$, 
where $S_2(\xi)$ has center 
$$
O(\xi)=\left(\frac u2-\frac v{\sqrt{u^2+v^2}}\xi, 
\frac v2+\frac u{\sqrt{u^2+v^2}}\xi\right)\in\RR^2, 
$$
and $S_1(\xi)$ is a plane perpendicular to the line joining 
$O(\xi)$ and 1. Then $\rho (S_1(\xi), S_2(\xi))$ is given by 
$$
\rho (Z(\xi))=\log \left( \frac
{\sqrt{1-u+\frac{u^2+v^2}4+\frac{2v}{\sqrt{u^2+v^2}}\xi+\xi^2}
+\sqrt{1-u+\frac{2v}{\sqrt{u^2+v^2}}\xi}}
{\sqrt{\frac{u^2+v^2}4+\xi^2}} \right), 
$$
which takes the maximum value 
$$
\frac{1}{2} \log \left( 1+\frac{2\left\{\sqrt{(u-1)^2+v^2}-(u-1)\right\}}{u^2+v^2} \right). 
$$


\hfill$\square$

As a corollary of Claim \ref{crandmodulus} there holds: 

\begin{coro} The integrand $|\Omega_{CR}|-\Re \Omega_{CR}$ of $E^{(2)}(K)$ can be 
interpreted as the maximal Lorentzian modulus 
$\mu _L(x, x+dx, y, y+dy) -1$ of a pair of infinitesimal segments. 
\end{coro} 


\medskip
\subsection{The circular Gauss map and the inverted open knots.}

Let $C(y,y,x)$  be as in definition  \ref{x_y_z_circle} and $v_y(x)$ be as in subsection \ref{proj_invert}. 
We remark that when $y$ approaches $x$ then 
$C(y,y,x)$ approaches the {\it osculating circle} at $x$, 
which will be denoted by $C(x,x,x)$, 
and $v_y(x)$ approaches the tangent vector at $x$, 
which will be denoted by $v_x(x)$. 

\begin{defi} \rm  
The {\it circular Gauss map} $\stackrel{\circ}{\Phi}_K$ 
of a knot $K$ is defined by 
$$
\stackrel{\circ}{\Phi}_K:S^1\times S^1\ni (s,t)\mapsto v_{f(s)}(f(t))\in S^2.
$$
\end{defi} 

We remark that $\stackrel{\circ}{\Phi}_K$ 
is well-defined at the diagonal set, 
but it is not a symmetric function. 

\begin{rema} 
The circular Gauss map $\stackrel{\circ}{\Phi}_K$ is not conformally 
invariant, namely, for a conformal transformation $T$ there is 
generally no element $g$ in $O(3)$ that satisfies 
$\stackrel{\circ}{\Phi}_{T(K)}(s,t)=g(\stackrel{\circ}{\Phi}_{K}(s,t))$ 
for any $(s,t)\in S^1\times S^1$. 
To see this one can put $T=I_x$ as an extreme case, when 
$\mbox{``$v_{T(x)}(T(y))$''}={I_x}_{\ast}(v_x(y))$ 
is constantly equal to $-v_x(x)$. 
\end{rema}

\begin{clai}
The degree of $\stackrel{\circ}{\Phi}_K$
vanishes for any knot $K$. 
\end{clai}

\begin{demo}
Any knot $K$ can be deformed continuously by an ambient isotopy to a very 
``thin position'' so that $K$ is contained in $xy$-plane except for some 'bridges' 
which are contained in $\{0\le z\le \epsilon\}$ as illustrated in Figure \ref{thinknot}. 
Then the preimage of the vector $(0,0,1)\in S^2$ by $\stackrel{\circ}{\Phi}_K$ 
consists of pairs of points near the crossing points.

\begin{figure}[htb] 
\centerline{\input{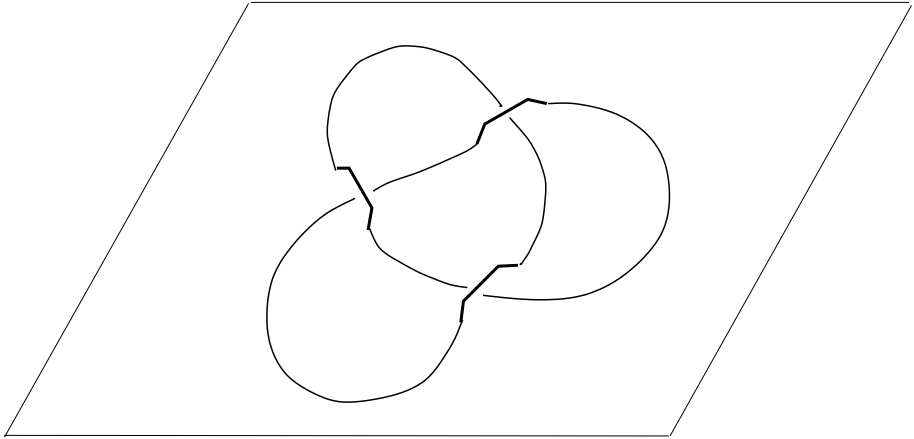}}
\caption{A knot in ``thin position''.}
\label{thinknot}
\end{figure}
\medskip
\noindent
Suppose a knot satisfies 
$$
\left\{ 
\begin{array}{cll} \nonumber 
&f(t)=(t, 0, \epsilon ) &\qquad (-\frac18\le t\le \frac18),  \nonumber\\[2mm]
&f(t)=(0, t-\frac12, 0) &\qquad (\frac12-\frac18\le t\le \frac12+\frac18) 
\nonumber \\
\end{array} \nonumber 
\right. 
$$
for a small $\epsilon$ $(0<\epsilon <\frac18)$. Put 
$$
U=[-\frac18, \frac18]\cup [\frac12-\frac18, \frac12+\frac18]. 
$$
The preimage 
$({\stackrel{\circ}{\Phi}_K|_{U\times U}})^{-1}((0,0,1))\in U\times U$ 
is a pair of points $(\epsilon, \frac12)$ and $(\frac12-\epsilon, 0)$, 
where the signatures of $(\stackrel{\circ}{\Phi}_K|_{U\times U})_{\ast}$ 
are opposite, hence the contribution of $U\times U$ to 
the degree of $\stackrel{\circ}{\Phi}_K$ is 0.

\begin{figure}[ht]
\centerline{\input{degreezero.pstex_t}}
\caption{$({\stackrel{\circ}{\Phi}_K})^{-1}((0, 0, 1))$.}
\end{figure}

\enddemo

Let us introduce two non-trivial functionals derived from 
$\stackrel{\circ}{\Phi}_K$. 
\begin{defi} \rm 
(1) Let $\omega _{S^2}$ be the unit volume 2-form of $S^2$: 
$$
\omega _{S^2}=\frac1{4\pi}
\frac{x_1dx_2\wedge dx_3+x_2dx_3\wedge dx_1+x_3dx_1\wedge dx_2}
{({x_1}^2+{x_2}^2+{x_3}^2)^{\frac32}}.
$$
The {\it absolute circular self-linking number \`a la Gauss}, $|csl|$, 
is defined by 
$$
|csl|(K)=\int_{T^2}|\stackrel{\circ}{\Phi}_K{}^{\ast}(\omega _{S^2})|. 
$$

(2) Let $l_K(x)$ be the length of the curve 
$\cup_{y\in K}v_y(x)=\stackrel{\circ}{\Phi}_K(S^1, s)$ on $S^2$. 
The {\it total variation of the conformally translated tangent vectors}, 
$I_{tv}$, is defined by 
$$
I_{tv}(K)=\int_{S^1}l_K(x)dx.
$$
\end{defi} 
By definition $I_{tv}(K)\ge |csl|(K)$ for any knot $K$. 

Let us give geometric interpretation of $I_{tv}$ using the inverted open knot. 
Let $I_x$, $\tilde{K}_x$, and $\tilde{y}$ be as 
in subsection \ref{cosineformula}. 
Then $I_x(C(y,y,x))$
 is the (oriented) tangent line to $\tilde K_x$ at $\tilde y$. 
Let $\tilde v_{\tilde y}(\tilde y)$ be the unit tangent vector to 
$I_x(C(y,y,x))$
 at $\tilde y$. Then 
$-v_y(x)={I_x}_{\ast}(v_y(x))=\tilde v_{\tilde y}(\tilde y)$. 
Therefore $l_K(x)$ is equal to the length of the curve 
$\cup _{\tilde y\in\tilde{K}_x}\tilde v_{\tilde y}(\tilde y)\subset S^2$. 
Thus $I_{tv}(K)$ can be considered as the 
{\it total variation of the tangent vectors to the inverted open knots}. 

If $T:\RR^3\cup\{\infty\}\to\RR^3\cup\{\infty\}$ is a conformal 
transformation that fixes $\{\infty\}$, then $T|_{\RR^3}$ is a homothety. 
Therefore we have: 
\begin{clai} 
The homothety class of 
the inverted open knot $\tilde K_x$ is 
a conformal invariant of $K$. 
Namely, $I_{T(x)}(T(K))$ is homothetic to $\tilde K_x$ for 
any conformal transformation $T$. 

Therefore the functional $I_{tv}$ is a conformally invariant functional. 
\end{clai}

Although $I_{tv}$ is conformally invariant, it can not be expressed 
in terms of the infinitesimal cross-ratio like $E^{(2)}$ and 
$E_{\sin\theta}$, because $I_{tv}$ is determined by up to the 
second order derivatives of the knot. 

\begin{clai} Neither $|csl|$ nor $I_{tv}$ is 
an energy functional for knots. 
\end{clai} 

\demo
Suffice to show that $I_{tv}$ is not an energy functional 
because $I_{tv}\ge |csl|$ for any knot. 

We show that the contribution to $I_{tv}$ 
of a pair of straight line segments of a fixed length 
with the closest distance $\epsilon $ does not blow up 
even if $\epsilon$ goes down to 0. 
Let $A_1=\{x(s)|a\le s\le b\}$ and 
$A_2=\{y(t)|a^{\prime}\le t\le b^{\prime}\}$ be straight line segments 
with $\min_{x\in A_1, y\in A_2}|x-y|=\epsilon$.  
Let us consider $I_{tv}|_{A_1\cup A_2}$. 
Put 
$$
u=f^{\prime}(s), v=f^{\prime}(t), \> {\rm and}\>\> 
w=\frac{f(t)-f(s)}{|f(t)-f(s)|} 
$$
as before. Then $u$ and $v$ are constant. 
Fix $t$ and put $\gamma (s)=v_x(y)$. 
Then as $\gamma (s)=2(u,w)w-u$ 
$$
\frac{d\gamma}{ds} =2\left(u, \frac{dw}{ds}\right)w+2(u,w)\frac{dw}{ds}, 
$$
hence 
$$
\left|\frac{d\gamma}{ds}\right| \le 2\sqrt 2\left|\frac{dw}{ds}\right|.
$$
Therefore 
$$
l_{A_1}(y) =\int_a^b\left|\frac{d\gamma}{ds}\right|ds
 \le 2\sqrt 2\int_a^b\left|\frac{dw}{ds}\right|ds
 \le 2\sqrt 2\pi,
$$
and hence 
$$
I_{tv}|_{A_1\cup A_2} =\int_{a^{\prime}}^{b^{\prime}}l_{A_1}(y)dt 
+\int_a^bl_{A_2}(x)ds
 \le 2\sqrt 2\pi (b^{\prime}-a^{\prime}+b-a),
$$
which is independent of $\epsilon$. 
\end{demo}

We show that they can detect the unknot. 
\begin{exam}\rm 
Let $K_{\circ}$ be the standard planar circle. 
Then $v_y(x)=f^{\prime}(s)$ for any $y$ and therefore 
$I_{tv}(K_{\circ})= |csl|(K_{\circ})=0$ 
\end{exam}

\begin{theo}
If $K$ is a non-trivial knot then $I_{tv}(K)\ge\pi$.
\end{theo}
\begin{demo}
The unit tangent vector $\tilde v_{\tilde y}(\tilde y)$ of 
the inverted open knot $\tilde K_x$ is asymptotic to 
$-v_x(x)$ as $\tilde y$ goes to $\infty$. 
If the angle between $\tilde v_{\tilde y}(\tilde y)$ and $-v_x(x)$ 
is smaller than $\frac{\pi}2$ for any $\tilde y$ then $\tilde K_x$ 
is unknotted, which contradicts the assumption. Therefore 
$l_K(x)\ge\pi$ for any $x\in K$. 
\end{demo}
\begin{conj} 
There is a positive constant $C$ such that 
if $K$ is a non-trivial knot then $|csl|(K)\ge C$. 
\end{conj}


\subsection{The global radius of curvature. \label{globalcurvature}} 

Let us make a remark on the relation between the 
infimum of the radii of non-trivial spheres 
and the global radius of curvature by 
Gonzalez and Maddocks (\cite{Go-Ma}). 
Let us consider a knot $K$ in $\RR^3$. 

Let $r(C(x,y,z))$ be the radius of the circle $C(x,y,z)$ 
that passes through $x,y$ and $z$ in $K$. 
When two (or three) of them coincide $C(x,y,z)$ means the tangent 
(or respectively, osculating) circle. 
The {\it global radius of curvature} $\rho_K^{(3)}(x)$ of a knot $K$ at 
$x$ is defined by 
$$
\rho_K^{(3)}(x)=\min _{y,z\in K}r(C(x,y,z)). 
$$
This $\rho_K^{(3)}(x)$ is attained by a triple $(x,y,y)$, 
and especially, when it is attained by a triple $(x,x,x)$, 
$\rho_K^{(3)}(x)$ is the radius of curvature in the ordinary sense. 

Let $\sigma (x,y,z,w)$ denote the smallest sphere 
that passes through $x,y,z$ and $w$ in $K$, 
where the points in $\sigma\cap K$ are counted with multiplicity 
according to the order of the tangency. 
Generically $\sigma (x,y,z,w)$ is uniquely determined. 
Let $r(\sigma (x,y,z,w))$ be the radius of $\sigma (x,y,z,w)$. 
Put 
$$
\rho_K^{(4)}(x)=\inf _{(y,z,w)\in K}
r(\sigma (x,y,z,w)). 
$$
%
Since $r(\sigma (x,y,z,w))\ge r(C(x,y,z))$ there holds 
$\rho_K^{(4)}(x)\ge \rho_K^{(3)}(x)$. 
Put 
\begin{eqnarray*}
\triangle (K)&=&\min _{x\in K}\rho_K^{(3)}(x)
=\min _{x,y,z \in K}r(C(x,y,z)), \\[2mm]
\square (K)&=&\inf _{x\in K} \rho_K^{(4)}(x)
=\inf _{(x,y,z,w)\in K}r(\sigma (x,y,z,w)). 
\end{eqnarray*}
Then $\triangle (K)$ is equal to the {\it thickness} 
studied by Buck and Simon 
et al. 
As is shown in \cite{Go-Ma} 
$\triangle (K)$ is either the minimum local radius of curvature or 
the strictly smaller radius of the sphere $\sigma =\partial D^3$ 
with ${\rm Int} D^3\cap K=\phi$ which is 
twice tangent to $K$ at antipodal points $x$ and $y$. 

\begin{clai} {\rm \cite{Ma-Sm} \cite{Go-Ma-Sm} } 
$\rho_K^{(3)}(x)=\rho_K^{(4)}(x)$ for any $x$. Therefore 
$\triangle (K)=\square (K)$. 
\end{clai}
\vskip 3mm 
\subsection{Links}
One can ask similar questions for links ${\mathcal L} = C_1 \cup C_2$

We say that the link is {\it splittable} if there exist two disjoint closed balls $B_1$ and $B_2$ such that $C_1 \subset B_1$ and $C_2 \subset B_2$.

\begin{enumerate}
\item What can be said about the measure of spheres ({\it non-trivial spheres}) intersecting the two components if the link is not splittable?

\item Define a non-trivial zone for a link as a zone crossed by the two components (then at least two stands of one and two strands of the other component cross the zone), and the modulus of a link as the maximum of the modulus of a non-trivial zone for the link. Give a lower bound for the modulus of a non-splittable link. 

\end{enumerate}

A recent result of Langevin and Moniot \cite{La-Mo} is:

\begin{clai}
The modulus of a non-splittable link is bounded below by the modulus of the Hopf Link: 
$${\mathcal H} = C_1 \cup C_2 ; \ C_1 = \{x=0 \} \cap S^3 , \ C_2 = \{ y = 0 \} \cap S^3; \ S^3 \subset \CC$$ 
\end{clai}

As for knots, the claim provides a lower bound for the measure of non-trivial spheres for a non-splittable link.

\newpage

\bibliographystyle{plain}

\end{document}